\documentclass[journal]{IEEEtran}

\pdfoutput=1

\usepackage{cite}
\usepackage{graphicx}
\usepackage[breaklinks=false, allcolors=black, colorlinks=true]{hyperref}

\usepackage{amsfonts}
\usepackage{amsmath} 
\usepackage{amssymb}
\usepackage{amsthm}	
\usepackage{enumerate}
\usepackage{mathrsfs}
\usepackage{subcaption}
\usepackage{epstopdf}
\usepackage{url}
\usepackage{booktabs}
\usepackage{tabularx}

\theoremstyle{definition} 
\newtheorem{theore}{Theorem}
\newtheorem{corollar}{Corollary}
\newtheorem{lemm}[theore]{Lemma}
\newtheorem{definitio}{Definition}
\newtheorem{proble}{Problem}
\newtheorem{assumptio}{Assumption}
\newtheorem{propert}{Property}
\newtheorem{propositio}{Proposition}
\newtheorem{conjectur}{Conjecture}
\newtheorem{exampl}{Example}
\newtheorem{rul}{Rule}

\newtheorem{remar}{Remark}
\newtheorem{fac}{Fact}

\def\nn{\nonumber}
\def\beq{\begin{equation}}
\def\eeq{\end{equation}}
\def\beqs{\begin{equation*}}
\def\eeqs{\end{equation*}}
\def\balign{\begin{align}}
\def\ealign{\end{align}}
\def\bea{\begin{eqnarray}}
\def\eea{\end{eqnarray}}
\def\ba{\begin{array}}
\def\ea{\end{array}}
\def\bcenter{\begin{center}}
\def\ecenter{\end{center}}
\def\bitem{\begin{itemize}}
\def\eitem{\end{itemize}}
\def\benum{\begin{enumerate}}
\def\eenum{\end{enumerate}}
\def\bdoc{\begin{document}}
\def\edoc{\end{document}}
\def\defeq{{\stackrel{\Delta}{=}}}
\def\nn{\nonumber}

\def\eg{{e.g., \/}}
\def\etal{{\it et al. \/}}
\def\viz{{\iet viz.,\ \/}}
\def\ie{{i.e.,\ \/}}
\def\etc{{\it etc.,\ \/}}
\def\cf{{cf.,\ \/}}

\usepackage[usenames,dvipsnames,svgnames,table]{xcolor}

\def\yellow{\textcolor{yellow}}
\def\green{\textcolor{green}}
\def\tcbg{\textcolor{bgrd}}
\def\magenta{\textcolor{magenta}}
\def\white{\textcolor{white}}
\def\gray{\textcolor{gray}}

\newcommand\Quote[1]{\lq\textsl{#1}\rq}
\newcommand\fr[2]{{\textstyle\frac{#1}{#2}}}

\newcommand{\Ambb}{\mathbb{A}}
\newcommand{\Bmbb}{\mathbb{B}}
\newcommand{\Cmbb}{\mathbb{C}}
\newcommand{\Dmbb}{\mathbb{D}}
\newcommand{\Embb}{\mathbb{E}}
\newcommand{\Fmbb}{\mathbb{F}}
\newcommand{\Gmbb}{\mathbb{G}}
\newcommand{\Hmbb}{\mathbb{H}}
\newcommand{\Imbb}{\mathbb{I}}
\newcommand{\Jmbb}{\mathbb{J}}
\newcommand{\Kmbb}{\mathbb{K}}
\newcommand{\Lmbb}{\mathbb{L}}
\newcommand{\Mmbb}{\mathbb{M}}
\newcommand{\Nmbb}{\mathbb{N}}
\newcommand{\Ombb}{\mathbb{O}}
\newcommand{\Pmbb}{\mathbb{P}}
\newcommand{\Qmbb}{\mathbb{Q}}
\newcommand{\Rmbb}{\mathbb{R}}
\newcommand{\Smbb}{\mathbb{S}}
\newcommand{\Tmbb}{\mathbb{T}}
\newcommand{\Umbb}{\mathbb{U}}
\newcommand{\Vmbb}{\mathbb{V}}
\newcommand{\Wmbb}{\mathbb{W}}
\newcommand{\Xmbb}{\mathbb{X}}
\newcommand{\Ymbb}{\mathbb{Y}}
\newcommand{\Zmbb}{\mathbb{Z}}

\newcommand{\Amsc}{\mathcal{A}}
\newcommand{\Bmsc}{\mathcal{B}}
\newcommand{\Cmsc}{\mathcal{C}}
\newcommand{\Dmsc}{\mathcal{D}}
\newcommand{\Emsc}{\mathcal{E}}
\newcommand{\Fmsc}{\mathcal{F}}
\newcommand{\Gmsc}{\mathcal{G}}
\newcommand{\Hmsc}{\mathcal{H}}
\newcommand{\Imsc}{\mathcal{I}}
\newcommand{\Jmsc}{\mathcal{J}}
\newcommand{\Kmsc}{\mathcal{K}}
\newcommand{\Lmsc}{\mathcal{L}}
\newcommand{\Mmsc}{\mathcal{M}}
\newcommand{\Nmsc}{\mathcal{N}}
\newcommand{\Omsc}{\mathcal{O}}
\newcommand{\Pmsc}{\mathcal{P}}
\newcommand{\Qmsc}{\mathcal{Q}}
\newcommand{\Rmsc}{\mathcal{R}}
\newcommand{\Smsc}{\mathcal{S}}
\newcommand{\Tmsc}{\mathcal{T}}
\newcommand{\Umsc}{\mathcal{U}}
\newcommand{\Vmsc}{\mathcal{V}}
\newcommand{\Wmsc}{\mathcal{W}}
\newcommand{\Xmsc}{\mathcal{X}}
\newcommand{\Ymsc}{\mathcal{Y}}
\newcommand{\Zmsc}{\mathcal{Z}}

\newcommand{\Ascr}{\mathscr{A}}
\newcommand{\Bscr}{\mathscr{B}}
\newcommand{\Cscr}{\mathscr{C}}
\newcommand{\Dscr}{\mathscr{D}}
\newcommand{\Escr}{\mathscr{E}}
\newcommand{\Fscr}{\mathscr{F}}
\newcommand{\Gscr}{\mathscr{G}}
\newcommand{\Hscr}{\mathscr{H}}
\newcommand{\Iscr}{\mathscr{I}}
\newcommand{\Jscr}{\mathscr{J}}
\newcommand{\Kscr}{\mathscr{K}}
\newcommand{\Lscr}{\mathscr{L}}
\newcommand{\Mscr}{\mathscr{M}}
\newcommand{\Nscr}{\mathscr{N}}
\newcommand{\Oscr}{\mathscr{O}}
\newcommand{\Pscr}{\mathscr{P}}
\newcommand{\Qscr}{\mathscr{Q}}
\newcommand{\Rscr}{\mathscr{R}}
\newcommand{\Sscr}{\mathscr{S}}
\newcommand{\Tscr}{\mathscr{T}}
\newcommand{\Uscr}{\mathscr{U}}
\newcommand{\Vscr}{\mathscr{V}}
\newcommand{\Wscr}{\mathscr{W}}
\newcommand{\Xscr}{\mathscr{X}}
\newcommand{\Yscr}{\mathscr{Y}}
\newcommand{\Zscr}{\mathscr{Z}}

\def\Atiny{\mbox{\tiny{A}}}
\def\Btiny{\mbox{\tiny{B}}}
\def\Ctiny{\mbox{\tiny{C}}}
\def\Dtiny{\mbox{\tiny{D}}}
\def\Etiny{\mbox{\tiny{E}}}
\def\Ctiny{\mbox{\tiny{C}}}
\def\Ftiny{\mbox{\tiny{F}}}
\def\Gtiny{\mbox{\tiny{G}}}
\def\Htiny{\mbox{\tiny{H}}}
\def\Itiny{\mbox{\tiny{I}}}
\def\Jtiny{\mbox{\tiny{J}}}
\def\Ktiny{\mbox{\tiny{K}}}
\def\Ltiny{\mbox{\tiny{L}}}
\def\Mtiny{\mbox{\tiny{M}}}
\def\Ntiny{\mbox{\tiny{N}}}
\def\Otiny{\mbox{\tiny{O}}}
\def\Ptiny{\mbox{\tiny{P}}}
\def\Qtiny{\mbox{\tiny{Q}}}
\def\Rtiny{\mbox{\tiny{R}}}
\def\Stiny{\mbox{\tiny{S}}}
\def\Ttiny{\mbox{\tiny{T}}}
\def\Utiny{\mbox{\tiny{U}}}
\def\Vtiny{\mbox{\tiny{V}}}
\def\Wtiny{\mbox{\tiny{W}}}
\def\Xtiny{\mbox{\tiny{X}}}
\def\Ytiny{\mbox{\tiny{Y}}}
\def\Ztiny{\mbox{\tiny{Z}}}

\def\alphabf{\hbox{\boldmath$\alpha$\unboldmath}}
\def\betabf{\hbox{\boldmath$\beta$\unboldmath}}
\def\gammabf{\hbox{\boldmath$\gamma$\unboldmath}}
\def\deltabf{\hbox{\boldmath$\delta$\unboldmath}}
\def\epsilonbf{\hbox{\boldmath$\epsilon$\unboldmath}}
\def\zetabf{\hbox{\boldmath$\zeta$\unboldmath}}
\def\etabf{\hbox{\boldmath$\eta$\unboldmath}}
\def\iotabf{\hbox{\boldmath$\iota$\unboldmath}}
\def\kappabf{\hbox{\boldmath$\kappa$\unboldmath}}
\def\lambdabf{\hbox{\boldmath$\lambda$\unboldmath}}
\def\mubf{\hbox{\boldmath$\mu$\unboldmath}}
\def\nubf{\hbox{\boldmath$\nu$\unboldmath}}
\def\xibf{\hbox{\boldmath$\xi$\unboldmath}}
\def\pibf{\hbox{\boldmath$\pi$\unboldmath}}
\def\rhobf{\hbox{\boldmath$\rho$\unboldmath}}
\def\sigmabf{\hbox{\boldmath$\sigma$\unboldmath}}
\def\taubf{\hbox{\boldmath$\tau$\unboldmath}}
\def\upsilonbf{\hbox{\boldmath$\upsilon$\unboldmath}}
\def\phibf{\hbox{\boldmath$\phi$\unboldmath}}
\def\chibf{\hbox{\boldmath$\chi$\unboldmath}}
\def\psibf{\hbox{\boldmath$\psi$\unboldmath}}
\def\omegabf{\hbox{\boldmath$\omega$\unboldmath}}
\def\Sigmabf{\hbox{$\bf \Sigma$}}
\def\Upsilonbf{\hbox{$\bf \Upsilon$}}
\def\Omegabf{\hbox{$\bf \Omega$}}
\def\Deltabf{\hbox{$\bf \Delta$}}
\def\Gammabf{\hbox{$\bf \Gamma$}}
\def\Thetabf{\hbox{$\bf \Theta$}}
\def\Lambdabf{\mbox{$ \bf \Lambda $}}
\def\Xibf{\hbox{\bf$\Xi$}}
\def\Pibf{{\bf \Pi}}
\def\thetabf{{\mbox{\boldmath$\theta$\unboldmath}}}
\def\Upsilonbf{\hbox{\boldmath$\Upsilon$\unboldmath}}
\newcommand{\Phibf}{\mbox{${\bf \Phi}$}}
\newcommand{\Psibf}{\mbox{${\bf \Psi}$}}

\def\abf{{\bf a}}
\def\bbf{{\bf b}}
\def\cbf{{\bf c}}
\def\dbf{{\bf d}}
\def\ebf{{\bf e}}
\def\fbf{{\bf f}}
\def\gbf{{\bf g}}
\def\hbf{{\bf h}}
\def\ibf{{\bf i}}
\def\jbf{{\bf j}}
\def\kbf{{\bf k}}
\def\lbf{{\bf l}}
\def\mbf{{\bf m}}
\def\nbf{{\bf n}}
\def\obf{{\bf o}}
\def\pbf{{\bf p}}
\def\qbf{{\bf q}}
\def\rbf{{\bf r}}
\def\sbf{{\bf s}}
\def\tbf{{\bf t}}
\def\ubf{{\bf u}}
\def\vbf{{\bf v}}
\def\wbf{{\bf w}}
\def\xbf{{\bf x}}
\def\ybf{{\bf y}}
\def\zbf{{\bf z}}
\def\rbf{{\bf r}}
\def\xbf{{\bf x}}
\def\ybf{{\bf y}}
\def\Abf{{\bf A}}
\def\Bbf{{\bf B}}
\def\Cbf{{\bf C}}
\def\Dbf{{\bf D}}
\def\Ebf{{\bf E}}
\def\Fbf{{\bf F}}
\def\Gbf{{\bf G}}
\def\Hbf{{\bf H}}
\def\Ibf{{\bf I}}
\def\Jbf{{\bf J}}
\def\Kbf{{\bf K}}
\def\Lbf{{\bf L}}
\def\Mbf{{\bf M}}
\def\Nbf{{\bf N}}
\def\Obf{{\bf O}}
\def\Pbf{{\bf P}}
\def\Qbf{{\bf Q}}
\def\Rbf{{\bf R}}
\def\Sbf{{\bf S}}
\def\Tbf{{\bf T}}
\def\Ubf{{\bf U}}
\def\Vbf{{\bf V}}
\def\Wbf{{\bf W}}
\def\Xbf{{\bf X}}
\def\Ybf{{\bf Y}}
\def\Zbf{{\bf Z}}
\def\Ac{{\cal A}}
\def\Bc{{\cal B}}
\def\Cc{{\cal C}}
\def\Dc{{\cal D}}
\def\Ec{{\cal E}}
\def\Fc{{\cal F}}
\def\Gc{{\cal G}}
\def\Hc{{\cal H}}
\def\Ic{{\cal I}}
\def\Jc{{\cal J}}
\def\Kc{{\cal K}}
\def\Lc{{\cal L}}
\def\Mc{{\cal M}}
\def\Nc{{\cal N}}
\def\Oc{{\cal O}}
\def\Pc{{\cal P}}
\def\Qc{{\cal Q}}
\def\Rc{{\cal R}}
\def\Sc{{\cal S}}
\def\Tc{{\cal T}}
\def\Uc{{\cal U}}
\def\Vc{{\cal V}}
\def\Wc{{\cal W}}
\def\Xc{{\cal X}}
\def\Yc{{\cal Y}}
\def\Zc{{\cal Z}}

\def\0bf{\mathbf{0}}
\def\1bf{\mathbf{1}}

\def\lambdabf{\boldsymbol{\lambda}}
\def\mubf{\boldsymbol{\mu}}
\def\gammabf{\boldsymbol{\gamma}}
\def\xibf{\boldsymbol{\xi}}
\def\sigmabf{\boldsymbol{\sigma}}
\def\zetabf{\boldsymbol{\zeta}}

\def\Lambdabf{\boldsymbol{\Lambda}}
\def\Mubf{\boldsymbol{\Mu}}
\def\Gammabf{\boldsymbol{\Gamma}}
\def\Xibf{\boldsymbol{\Xi}}
\def\Sigmabf{\boldsymbol{\Sigma}}
\def\Zetabf{\boldsymbol{\Zeta}}

\usepackage{tikz}
\usetikzlibrary{arrows}

\usepackage{makecell}
\setcellgapes{.25ex}

\newcommand\independent{\protect\mathpalette{\protect\independenT}{\perp}}
\def\independenT#1#2{\mathrel{\rlap{$#1#2$}\mkern2mu{#1#2}}}

\newcommand{\bsf}[1]{\textsf{\textbf{#1}}}
\newcommand{\lbsf}[1]{\textsf{\large  \textbf{#1}}}
\newcommand{\Lbsf}[1]{\textsf{\Large  \textbf{#1}}}
\newcommand{\hbsf}[1]{\textsf{\huge  \textbf{#1}}}

\newcommand{\myminipage}[3]{\begin{minipage}[#1]{#2}{#3} \end{minipage}}
\newcommand{\sbs}[4]{\myminipage{c}{#1}{#3} \hfill
\myminipage{c}{#2}{#4}}

\newcommand{\myfig}[2]{\centerline{\psfig{figure=#1,width=#2,silent=}}}
\newcommand{\myfigh}[2]{\centerline{\psfig{figure=#1,height=#2,silent=}}}
\newcommand{\myfigwh}[3]{\centerline{\psfig{figure=#1,width=#2,height=#3,silent=}}}

\newcommand{\beqa}{\begin{eqnarray}}
\newcommand{\eeqa}{\end{eqnarray}}
\newcommand{\beqan}{\begin{eqnarray*}}
\newcommand{\eeqan}{\end{eqnarray*}}
\newcommand{\dst}[1]{\displaystyle{ #1 }}
\newcounter{pf}

\newcommand{\smax}[1] { \bar \sigma \left( #1 \right) }
\newcommand{\Rn}{{\mathbb R}^n}
\newcommand{\R}{{\mathbb R}}
\newcommand{\C}{{\mathbb C}}
\newcommand{\Rm}{\mathbb{R}^m}
\newcommand{\Rmn}{\mathbb{R}^{m \times n}}
\newcommand{\Rpq}{\mathbb{R}^{p \times q}}
\newcommand{\Cn}{\mathbb{C}^n}
\newcommand{\Cm}{\mathbb{C}^m}
\newcommand{\Cnn}{\mathbb{C}^{n \times n}}
\newcommand{\Cmn}{\mathbb{C}^{m \times n}}
\newcommand{\ip}[1]{\left\langle #1 \right\rangle}
\newcommand{\rank}{\mbox{rank}}
\newcommand{\Span}{\mbox{\rm Span }}
\newcommand{\Trace}{\mbox{\rm Tr }}
\newcommand{\trace}[1]{\text{Tr}\left(#1\right)}
\newcommand{\Spec}{\mbox{\rm Spec }}
\newcommand{\vectornorm}[1]{\left\|#1\right\|}

\newcommand{\pd}[2]{\frac{\partial #1}{\partial #2}}
\newcommand{\ppd}[3]{\frac{\partial^2 #1}{\partial #2 \partial #3}}

\newcommand{\thtilde}{\tilde{\theta}}
\newcommand{\thnom}{\theta^\circ}
\newcommand{\thopt}{\theta^{\mbox{\small opt}}}
\newcommand{\thhat}{{\hat{\theta}}}
\newcommand{\Tho}{\Theta^\circ}
\newcommand{\tho}{\theta^\circ}
\newcommand{\np}{{n_p}}

\newcommand{\ii}{{[i]}}
\newcommand{\II}{{[i+1]}}
\newcommand{\iii}{{[ii]}}
\newcommand{\jj}{{[j]}}
\newcommand{\kk}{{[k]}}
\newcommand{\thi}{{\theta^\ii}}
\newcommand{\thI}{{\theta^\II}}
\newcommand{\di}{{d^\ii}}
\newcommand{\gi}{{g^\ii}}
\newcommand{\Hi}{{\HH^\ii}}
\newcommand{\thK}{\theta^{(k+1)}}
\newcommand{\gk}{{g^{(k)}}}
\newcommand{\Hk}{{{\cal H}^{(k)}}}

\newcommand{\bfdelta}{{\bf \Delta}}

\newcommand{\Exp}[1]{\exp \left\{ #1 \right\}} 
\newcommand{\gaussian}[1]{\mathbb{N} \left( #1 \right)}
\newcommand{\uniform}[1]{\mathbb{U} \left[ #1 \right]}
\newcommand{\exponential}[1]{\mathbb{E} \left[ #1 \right]}
\newcommand{\EXP}[1]{\EEXP \left[ #1 \right]} 
\newcommand{\EEXP}{\mbox{\bsf{E}}} 
\newcommand{\Prob}[1]{\mbox{{\sf Pr}} \left(#1 \right)}
\newcommand{\convas}{\stackrel{as}{\longrightarrow}}
\newcommand{\convinp}{\stackrel{p}{\longrightarrow}}
\newcommand{\convind}{\stackrel{d}{\longrightarrow}}
\newcommand{\convqm}{\stackrel{qm}{\longrightarrow}}
\newcommand{\sss}[1]{{_{#1}}}
\newcommand{\density}[2]{p_{_{_{#1}}}\!\!\left(#2 \right)} 
\newcommand{\distro}[2]{P_{_{_{#1}}}\!\!\left(#2 \right)} 
\newcommand{\rxx}[1]{R_{_{#1}}\!} 
\newcommand{\sxx}[1]{S_{_{#1}}} 
\newcommand{\cov}[1]{\Lambda_{_{#1}}} 
\newcommand{\mean}[1]{m_{_{#1}}} 
\newcommand{\LS}[1]{\hat{#1}_{_{LS}}} 
\newcommand{\MV}[1]{\hat{#1}_{_{MV}}} 
\newcommand{\LMV}[1]{\hat{#1}_{_{LMV}}} 
\newcommand{\ML}[1]{\hat{#1}_{_{ML}}} 

\renewcommand{\arraystretch}{0.9}
\newcommand{\bmat}[1]{ \begin{bmatrix} #1 \end{bmatrix}}
\newcommand{\pmat}[1]{ \begin{pmatrix} #1 \end{pmatrix}}
\newcommand{\mat}[1]{ \left[ \begin{array}{cccccccc} #1 \end{array}
\right] }
\newcommand{\smallmat}[1]{\small{\mat{#1}}}
\newcommand{\sysblk}[4]{\begin{array}{c|cccc}#1&#2\\ \hline#3&#4
\end{array}}
\newcommand{\sysmat}[4]{\left[\sysblk{#1}{#2}{#3}{#4}\right]}
\newcommand{\SGeq}{\succ}
\newcommand{\SLeq}{\prec}
\newcommand{\Geq}{\succeq}
\newcommand{\Leq}{\preceq}

\newcommand{\Bset}{\mathbb{B}}
\newcommand{\Cset}{\mathbb{C}}
\newcommand{\Fset}{\mathbb{F}}
\newcommand{\Mset}{\mathbb{M}}
\newcommand{\Nset}{\mathbb{N}}
\newcommand{\Qset}{\mathbb{Q}}
\newcommand{\Rset}{\mathbb{R}}
\newcommand{\Sset}{\mathbb{S}}
\newcommand{\Tset}{\mathbb{T}}
\newcommand{\Uset}{\mathbb{U}}
\newcommand{\Vset}{\mathbb{V}}
\newcommand{\Wset}{\mathbb{W}}
\newcommand{\Zset}{\mathbb{Z}}

\newcommand{\Ical}{{\cal I}}
\newcommand{\Acal}{{\cal A}}
\newcommand{\Bcal}{{\cal B}}
\newcommand{\Ccal}{{\cal C}}
\newcommand{\Dcal}{{\cal D}}
\newcommand{\Ecal}{{\cal E}}
\newcommand{\Fcal}{{\cal F}}
\newcommand{\Gcal}{{\cal G}}
\newcommand{\Hcal}{{\cal H}}
\newcommand{\Kcal}{{\cal K}}
\newcommand{\Lcal}{{\cal L}}
\newcommand{\Mcal}{{\cal M}}
\newcommand{\Ncal}{{\cal N}}
\newcommand{\Pcal}{{\cal P}}
\newcommand{\Qcal}{{\cal Q}}
\newcommand{\Rcal}{{\cal R}}
\newcommand{\Scal}{{\cal S}}
\newcommand{\Tcal}{{\cal T}}
\newcommand{\Wcal}{{\cal W}}
\newcommand{\Ucal}{{\cal U}}
\newcommand{\Vcal}{{\cal V}}
\newcommand{\Xcal}{{\cal X}}
\newcommand{\Zcal}{{\cal Z}}

\newcommand{\EE}{{\bf E}}
\newcommand{\FF}{{\bf F}}
\newcommand{\GG}{{\bf G}}
\newcommand{\HH}{{\bf H}}
\newcommand{\LL}{{\bf L}}
\newcommand{\NN}{{\bf N}}
\newcommand{\MM}{{\bf M}}
\newcommand{\PP}{{\bf P}}
\newcommand{\QQ}{{\bf Q}}
\newcommand{\RR}{{\bf R}}
\renewcommand{\SS}{{\bf S}}
\newcommand{\TT}{{\bf T}}
\newcommand{\VV}{{\bf V}}
\newcommand{\WW}{{\bf W}}

\newcommand{\thk}{\theta^{(k)}}
\newcommand{\thb}{\theta^{\rm opt}}
\newcommand{\alb}{\alpha^{\rm opt}}
\newcommand{\dk}{d^{(k)}}
\newcommand{\Hinf}{{\cal H}_\infty}
\newcommand{\Htwo}{{\cal H}_2}

\renewcommand{\arraystretch}{1.1}

\newcommand{\red}[1]{{\color{red} #1}}
\newcommand{\blue}[1]{{\color{Blue} #1}}
\newcommand{\black}[1]{{\color{Black} #1}}

\newcounter{l1}
\newcounter{l2}
\newcounter{l3}
\newcommand{\bdotlist}{\begin{list}{$\bullet$}{}}
\newcommand{\bboxlist}{\begin{list}{$\Box$}{}}
\newcommand{\bbboxlist}{\begin{list}{\raisebox{.005in}{{\tiny
$\blacksquare$ \ \ }}}{}}
\newcommand{\bdashlist}{\begin{list}{$-$}{} }
\newcommand{\blist}{\begin{list}{}{} }
\newcommand{\barablist}{\begin{list}{\arabic{l1}}{\usecounter{l1}}}
\newcommand{\balphlist}{\begin{list}{(\alph{l2})}{\usecounter{l2}}}
\newcommand{\bAlphlist}{\begin{list}{\Alph{l2}.}{\usecounter{l2}}}
\newcommand{\bdiamlist}{\begin{list}{$\diamond$}{}}
\newcommand{\bromalist}{\begin{list}{(\roman{l3})}{\usecounter{l3}}}

\newcommand{\thm}[1]{\noindent \begin{theorem} #1   \end{theorem}}
\newcommand{\prop}[1]{\begin{proposition} #1 \end{proposition}}
\newcommand{\lem}[1]{\begin{lemma} #1  \hfill $\blacksquare$ \end{lemma}}
\newcommand{\ex}[1]{\begin{example} {\rm #1} \end{example}}
\newcommand{\prf}[1]{ \noindent {\em Proof:} \, #1 \hfill $\blacksquare$}
\newcommand{\rem}[1]{\begin{remark} {\rm #1} \hfill $\Box$ \end{remark}}
\newcommand{\defn}[1]{\begin{definition} {\rm #1 } \end{definition}}
\newcommand{\prob}[1]{\begin{exercise} {\rm  #1 } \end{exercise}}
\newcommand{\cor}[1]{\begin{corollary}   #1  \end{corollary}}

\newcommand{\argmin}{\mathop{\rm argmin}}
\newcommand{\argmax}{\mathop{\rm argmax}}
\newcommand{\diag}{\mathop{\mathrm{diag}}}
\newcommand{\tr}{\mathop{\rm Tr}}
\newcommand{\conv}{\mathop{\rm conv}}
\newcommand{\var}{\mathop{\rm Var}}
\renewcommand{\b}[1]{\ensuremath{\boldsymbol{\mathrm{#1}}}}
\newcommand{\ms}{{\rm MS}}
\newcommand{\tcs}{{\rm TCS}}
\newcommand{\scs}{{\rm SCS}}

\newcommand{\E}[1]{\b{\mu}_{{#1}}}
\newcommand{\Var}[1]{{\Sigma_{#1}}}

\newcommand{\bone}{\mathbf{1}}
\newcommand{\bi}{{\bf i}}

\def\u{u}	
\def\U{U}	
\def\v{v}	
\def\V{V}	
\def\s{s}	
\def\S{S}	
\def\f{q}	
\def\q{q}	
\def\Q{Q}	
\def\c{c}	
\def\z{z}	
\def\xcap{b}
\def\ccap{\mathbf{\c}}
\def\rentftr{\Phi}
\def\rentfsr{\Sigma}
\def\rentDer{\Delta}
\def\lambdao{\lambda^{o}}

\def \add [#1]{\blue{#1}}
\def \remove [#1]{\red{#1}}
\def \replace [#1]#2{\red{#1} \blue{#2}}

\def \dm [#1]{\magenta{(\textbf{Daniel says:} #1)}}
\def \eb [#1]{\red{(\textbf{Eilyan says:} #1)}}

\def \rone{\black}

\graphicspath{{figures/}}

\IEEEoverridecommandlockouts

\begin{document}

\makeatletter
\renewcommand{\@maketitle}{%
    \newpage\null%
    \begin{center}%
    		\vspace{-.25in}
        \let\footnote\thanks{\LARGE\@title\par}%
        \vskip1.5em{\normalsize\lineskip0.5em\begin{tabular}[t]{c}\@author\end{tabular}\par }%
        \vspace{-.2in}
    \end{center}%
    \par}
\makeatother

\title{\Huge Decentralized Stochastic Control of \\ Distributed Energy Resources
}

\author{\vspace{.12in} Weixuan Lin   \qquad Eilyan Bitar 
\thanks{Weixuan Lin ({\tt\small wl476@cornell.edu}) and Eilyan Bitar ({\tt\small eyb5@cornell.edu}) are with the School of Electrical and Computer Engineering, Cornell University, Ithaca, NY, 14853, USA. }
}

\maketitle

\begin{abstract}
We consider the decentralized control of radial distribution systems with controllable  photovoltaic inverters and energy storage resources. For such systems, we investigate the problem of designing fully decentralized controllers that minimize the expected cost of balancing demand, while guaranteeing the satisfaction of  individual resource and distribution system voltage constraints.
Employing a linear approximation of the branch flow model,  we formulate this problem as the design of a decentralized disturbance-feedback controller that minimizes the expected value of a convex quadratic cost function, subject to robust convex quadratic constraints on the system state and input. As such problems are, in general, computationally intractable, we derive a tractable inner approximation to this decentralized control problem, which enables the efficient  computation of an affine control policy via the solution of a finite-dimensional conic program.
As affine policies are, in general, suboptimal for the family of systems considered, we provide an efficient method to bound their suboptimality via the optimal solution of another finite-dimensional conic program. 
A case study of a 12 kV radial distribution system demonstrates that decentralized affine controllers can perform close to optimal.
\end{abstract}



\section{Introudction}

The increasing penetration of distributed and renewable energy resources introduces challenges to the operation of distribution systems, including rapid fluctuations in bus voltage magnitudes, reverse power flows at distribution substations, and deteriorated power quality due to the intermittency of  supply from renewables. These challenges are exasperated by the fact that traditional techniques for distribution system management, including the deployment of on-load tap changing (OLTC) transformers and  shunt capacitors, cannot effectively deal with the rapid variation in the power supplied from renewable resources \cite{cheng2015mitigating}.  In this paper, we aim to address such challenges by developing a systematic approach to the design of decentralized feedback controllers for distribution networks with a high penetration of distributed solar and energy storage resources, in order to minimize the expected cost of meeting demand over a finite horizon, while respecting  network and  resource constraints.

\emph{Related Work:} \ 
Although current industry standards require that photovoltaic (PV) inverters operate at a unity power factor \cite{IEEE1547}, 
the latent reactive power capacity of PV inverters can be utilized to regulate voltage profiles  \cite{turitsyn2011options, farivar2013equilibrium, ferreira2013distributed, li2014real, robbins2013two, zhu2015fast}, and reduce active power losses \cite{bolognani2013distributed,  dall2013distributed, dall2014decentralized, sulc2014optimal, zhang2015optimal, farivar2012optimal, dall2014optimal, dall2015optimal, bazrafshan2014decentralized, kekatos2015stochastic} in distribution networks.
A large swath of the literature on the reactive power management of PV inverters prescribes the solution of an optimal power flow (OPF) problem to determine the reactive power injections of PV inverters in real time
\cite{turitsyn2011options, li2014real, farivar2013equilibrium, ferreira2013distributed, robbins2013two, zhu2015fast, farivar2012optimal, bolognani2013distributed, dall2013distributed, dall2014optimal, dall2014decentralized, sulc2014optimal, zhang2015optimal}. 
The resulting OPF problem must be repeatedly solved over fast time scales (e.g., every minute) to accommodate the  rapid fluctuations in the active power supplied from the PV resources. 
In the presence of a large number of PV resources, the sheer size of the resulting OPF problem that needs to be solved, and the communication requirements it entails, gives rise to the need for distributed optimization methods
\cite{li2014real, ferreira2013distributed, robbins2013two, zhu2015fast, bolognani2013distributed, dall2013distributed, dall2014decentralized, sulc2014optimal, zhang2015optimal}. In particular, there has emerged a recent stream of literature developing distributed optimization methods, which enable the real-time control of reactive power injections of  PV inverters using only local measurements of bus voltage magnitudes \cite{li2014real, ferreira2013distributed, robbins2013two, zhu2015fast}. Under the assumption that the 
underlying OPF problem being solved is time-invariant,  such methods are guaranteed to asymptotically converge to the globally optimal reactive power injection profile. There is, however, no guarantee on the performance or constraint-satisfaction of these methods in finite time. The aforementioned methods can be interpreted as being fully decentralized, in that they do not require the explicit exchange of information between local controllers. Instead, the local controllers can be interpreted as 
communicating implicitly through the distribution network, which couples them physically. There exists another class of distributed optimization methods, which rely on the explicit exchange of information between neighboring controllers through a digital communication network \cite{bolognani2013distributed, dall2013distributed, dall2014decentralized, sulc2014optimal, zhang2015optimal}. Additionally,  there exists a related stream of literature, which aims to explicitly treat uncertainty in renewable supply and demand by leveraging on methods grounded in stochastic optimization \cite{dall2015optimal, bazrafshan2014decentralized, kekatos2015stochastic}.

In addition to the reactive power control of PV inverters, one can imagine a future power system in which a broader class of distributed energy resources with storage capability (e.g., electric vehicles, standalone battery packs)
are  actively controlled to  mitigate  voltage fluctuations and distribution system losses
\cite{traube2013mitigation, marra2013improvement, cheng2015mitigating, caramanis2009management, clement2009coordinated, deilami2011real, gan2013optimal, loukarakis2016decentralized, vrakopoulou2014stochastic}. As the set of feasible power injections to and withdrawals from energy storage systems are naturally coupled across time, 
the problem of managing their operation amounts to a multi-period, constrained stochastic control problem \cite{cheng2015mitigating, caramanis2009management, clement2009coordinated, deilami2011real, gan2013optimal, loukarakis2016decentralized, vrakopoulou2014stochastic}. In the presence of network constraints and uncertainty in demand and renewable supply,  the calculation of the optimal control policy is, in general, computationally intractable. 
These computational difficulties in control design are underscored by a recent report from the U.S. Department of Energy pointing to an apparent lack of effective control methods capable of  ``integrating [PV] inverter controls with control of other DERs or the management of uncertainty from intermittent generation'' \cite{DoE2016}[p. 31]. The development of computational methods to enable the tractable calculation of feasible control policies with computable bounds on their suboptimality is therefore desired, and stands as the primary subject of this paper.

\emph{Contribution:} \ 
The setting we consider entails the decentralized control of distributed energy resources spread throughout a radial distribution network,
subject to uncertainty in demand and renewable supply. 
The power flow equations over the radial network are described according to a linearized branch flow model.
Our objective is to minimize the expected amount of active power supplied at the substation required to meet demand, while guaranteeing the satisfaction of network and individual resource constraints. For the setting considered, this is technically equivalent to minimizing the expected active power losses plus the terminal energy stored in the distribution network.
The determination of an optimal decentralized control policy for such problems is, in general, computationally intractable, due to the presence of stochastic disturbances  and hard constraints on the system state and and input.
Our primary contributions are two-fold. First, we develop a convex programming approach to the design  of decentralized, affine disturbance-feedback controllers.
Second, as such control policies are,  in general, suboptimal, we provide a technique to bound their suboptimality through the solution of another convex program. We verify that the decentralized affine  policies we derive are close to optimal for the problem instance considered in our case study.

\emph{Organization:} \  
The remainder of this paper is organized as follows. 
Section \ref{sec:model} describes our models of the distribution network and the distributed energy resources.
Section \ref{sec:decent_control} formally states the decentralized control design problem. 
Section \ref{sec:affine} describes an approach to the computation of decentralized affine control policies via a finite-dimensional conic program. 
Section \ref{sec:LB} describes an approach to the tractable calculation of guaranteed bounds on the suboptimality incurred by these affine control policies via another finite-dimensional conic program.
Section \ref{sec:case_study} demonstrates the proposed techniques with a numerical study of a 12 kV radial distribution network. Section \ref{sec:Conclusions} concludes the paper.

\emph{Notation:} \ Let $\RR$ denote the set of real numbers. Let $e_i$ be the $i^{\rm th}$ real standard basis vector,
of dimension appropriate to the context in which it is used.
 Denote by $x'$ the transpose of vector $x \in \RR^n$.
For any pair of vectors $x = (x_1, \dots, x_n) \in \RR^n$ and $y = (y_1, \dots, y_m) \in \RR^m$, we define their concatenation as $(x,y) = (x_1, \dots, x_n, y_1, \dots, y_m) \in \RR^{n+m}$. Given a   process $\{x(t)\}$ indexed by  $t = 0, \cdots, T-1$, we denote by $x^t = (x(0), \cdots, x(t))$ its history up until and including time $t$. 
We denote by $I_n$ the $n$-by-$n$ identity matrix, and by $0_{m \times n}$ the $m$-by-$n$ zero matrix. Subscripts are omitted when the underlying matrix dimension is clear from the context.  
We denote the trace of a square matrix $A$ by $\trace{A}$.
Finally, we denote by $\Kcal$ a proper cone (i.e., convex, closed, and pointed with a nonempty interior). Let $\Kcal^*$ denote its dual cone. We write $x \succeq _\Kcal y$ to indicate that $x - y \in \Kcal$. For a matrix $A$ of appropriate dimension, 
$A \succeq_\Kcal 0$ denotes its columnwise inclusion in $\Kcal$.
\section{Network and Resource Models} \label{sec:model}

\subsection{Branch Flow Model}
Consider a radial distribution network whose topology is described  by a \emph{rooted tree} $\Gcal = (\Vcal, \Ecal)$, where $\Vcal = \{0, 1, .., n\}$ denotes its set of (nodes) buses, and $\Ecal$ its set of (directed edges) distribution lines. In particular, bus 0 is defined as the root of the network, and represents the substation that connects to the external power system. Each directed distribution line admits the natural orientation, i.e., away from the root. For each distribution line $(i, j) \in \Ecal$, we denote by $ r_{ij} + \bi x_{ij}$ its \emph{impedance}. In addition, define $I_{ij}$ as the \emph{complex current} flowing from bus $i$ to $j$, and $  p_{ij} + \bi q_{ij}$ as the \emph{complex power} flowing from bus $i$ to $j$. For each bus $i \in \Vcal$, let $V_i$ denote its \emph{complex voltage}, and $ p_i + \bi q_i$  the \emph{complex power injection} at this bus. We assume that the complex voltage $V_0$ at the substation is fixed and known.

We employ the \emph{branch flow model} proposed in \cite{baran1989optimal,baran1989optimal_sizing} to describe the steady-state, single-phase AC power flow equations associated with this radial distribution network. In particular, for each bus $j = 1, \dots, n$, and its unique \emph{parent} $i \in \Vcal$, we have
\begin{align}
-p_j \ &= \ p_{ij} - r_{ij} \ell_{ij} - \sum_{k: (j,k) \in \Ecal} p_{jk} , \label{eq:branch_quad1}\\
-q_j \ &= \ q_{ij} - x_{ij} \ell_{ij} - \sum_{k: (j,k) \in \Ecal} q_{jk} , \\ 
v_j^2 \ &= \ v_i^2 - 2 (r_{ij} p_{ij} + x_{ij} q_{ij}) + (r_{ij}^2 + x_{ij}^2) \ell_{ij} , \label{eq:branch_quad3}\\
\ell_{ij} \ &= \ (p_{ij}^2 + q_{ij}^2)/v_i^2 , \label{eq:branch_quad4}
\end{align}
where $\ell_{ij} = |I_{ij}|^2$ and $v_i = |V_i| $. We note that the branch flow model is well defined only for radial distribution networks, as we require that each bus $j$ (excluding the substation) have a unique parent $i \in \Vcal$.

For the remainder of the paper, we consider a linear approximation of the branch flow model \eqref{eq:branch_quad1}-\eqref{eq:branch_quad4} based on the Simplified Distflow method developed in \cite{baran1989network}.
The derivation of this approximation relies on the assumption that  $\ell_{ij} = 0$ for all $(i, j) \in \Ecal$, as the active and reactive power losses on distribution lines are considered small relative to the power flows. 
According to \cite{IEC2007efficient, farivar2013equilibrium}, such an approximation tends to introduces a relative model error of  1-5\% for practical distribution networks. 
Under this assumption, Eq. \eqref{eq:branch_quad1}-\eqref{eq:branch_quad3} can be simplified to 
\begin{align}
-p_j \ &= \ p_{ij}  - \sum_{k: (j,k) \in \Ecal} p_{jk} , \label{eq:branch_lin1}\\
-q_j \ &= \ q_{ij}  - \sum_{k: (j,k) \in \Ecal} q_{jk} , \label{eq:branch_lin2}\\ 
v_j^2 \ &= \ v_i^2 -  2(r_{ij} p_{ij} + x_{ij} q_{ij})  . \label{eq:branch_lin3}
\end{align}
The linearized branch flow Eq. \eqref{eq:branch_lin1}-\eqref{eq:branch_lin3} can be written more compactly as
\begin{align}
v^2 = Rp + Xq + v_0^2\mathbf{1}. \label{eq:v}
\end{align}
Here, $v^2= (v_1^2, .., v_n^2)$, $p = (p_1, .., p_n)$, and $q = (q_1, .., q_n)$ denote the vectors of squared bus voltage magnitudes, real power injections, and reactive power injections, respectively, and $\mathbf{1} = (1, .., 1)$ is a vector of all ones in $\RR^n$.
The matrices $R, X \in \RR^{n \times n}$ are defined according to
\begin{align*}
R_{ij} =  2 \hspace{-1em}  \sum_{(h, k) \in \Pcal_i \cap \Pcal_j} \hspace{-.15in} r_{hk}  , \\
X_{ij} =  2 \hspace{-1em} \sum_{(h, k) \in \Pcal_i \cap \Pcal_j} \hspace{-.15in}  x_{hk} , 
\end{align*}
where $\Pcal_i \subset \Ecal$ is defined as the set of edges on the unique path from bus 0 to $i$. 

In the sequel, we will consider the control of the distribution system over discrete time periods indexed by $t = 0, \cdots, T-1$.
Each discrete time period $t$ is defined over a time interval of length $\Delta$.
We require the vector of bus voltage magnitudes $v(t) = (v_1(t), .., v_n(t)) \in \RR^n$ at each time period $t$ to satisfy
\begin{align}
\underline{v} \leq \v(t) \leq \overline{v}, \label{eq:v_con}
\end{align}
where the allowable range of voltage magnitudes is defined by $\underline{v}, \overline{v} \in \RR^n$.

\subsection{Energy Storage Model} \label{sec:storage}
We consider a distribution system consisting of $n$ perfectly efficient energy storage devices, where each bus $i$ (excluding the substation) is assumed to have an energy storage capacity of $b_i \in \RR$.
The dynamic evolution of each energy storage device $i$ is described according to the state equation
\begin{align}
x_i (t+1) = x_i (t) - \Delta p^S_i (t), \quad t=0, \dots, T-1, \label{eq:storage}
\end{align}
where the state $x_i (t) \in \RR$ denotes the amount of energy stored in storage device $i$ just preceding period $t$, and $p^S_i (t) \in \RR$ denotes the active power extracted from device $i$ during period $t$.  \rone{For ease of exposition, we assume that the initial condition $x_i(0)$ of each storage device is fixed and known.}\footnote{\rone{
We emphasize that all  results presented in this paper are easily generalized to the setting in which the initial condition $x_i (0)$ is modeled as a random variable with known distribution. In particular, one can treat the initial condition as an additive disturbance to the state equation at time period $t=0$. We refer the readers to \cite{Lin2016} for a detailed  treatment of such systems.}} 
We impose state and input constraints of the form
\begin{align}
0 \leq x_i(t) \leq b_i, \quad &t = 0, \dots, T \label{eq:SoC_con}\\
\underline{p}^S_i \leq p^S_i (t) \leq \overline{p}^S_i, \quad &t = 0, \dots, T-1. \label{eq:pS_con}
\end{align}
for $i = 1, .., n$. The interval $[\underline{p}^S_i, \; \overline{p}^S_i] \subset \RR$ defines the range of allowable inputs for storage device $i$ at each time period $t$.

\subsection{Photovoltaic Inverter Model}

We assume that, in addition to energy storage capacity, each bus $i$ (excluding the substation) has a photovoltaic (PV) inverter whose reactive power injection can be actively controlled. 
We denote by $\xi_i^I (t) \in \RR$ the active power injection,  and by $q_i^I (t) \in \RR$ the reactive power injection from the PV inverter at bus $i$ and time $t$.
Due to  the intermittency  of solar irradiance, we will model $  \xi^I_i (t) $ as a discrete-time stochastic process, whose precise specification is presented in Section \ref{sec:uncert}. 
Additionally, we require that the reactive power injections  respect capacity constraints of the form
\begin{align}
\left| q_i^I (t) \right|   \leq \sqrt{\left. s_i^I \right.^2 -  \left.\xi_i^I (t) \right. ^2},  \quad i = 1, \dots, n, \label{eq:q_I}
\end{align}
for $t = 0, \dots, T-1$. Here, $s_i^I \in \RR$ denotes the apparent power capacity of PV inverter $i$.
Clearly, it must hold that $\xi_i^I (t) \leq s_i^I$.

\subsection{Load Model}
Each bus in the distribution network is assumed to have a constant power load, which we will treat as a discrete-time stochastic process.
Accordingly, we denote by $\xi^p_i(t) \in \RR$ and $\xi^q_i(t)  \in \RR$ the active and reactive power demand, respectively,  at bus $i$ and time $t$. It follows that the nodal active and reactive power balance equations can be  expressed as
\begin{align}
p_i(t) \ & = \ p^S_i (t) + \xi^I_i (t) -  \xi^p_i (t),  \label{eq:p}  \\
q_i(t) \ & = \ q^I_i (t) - \xi^q_i (t),  \label{eq:q}
\end{align}
where $p_i(t) \in \RR$ and $q_i(t) \in \RR$ denote the net active and reactive power injections, respectively, at each bus $i = 1, \dots,n$ and time period $t=0,\dots, T-1$.

\subsection{Uncertainty Model} \label{sec:uncert}
As indicated earlier, we model the  active power demand,  reactive power demand, and PV active power supply as discrete-time stochastic processes. 
Accordingly, we associate with each bus $i$ a \emph{disturbance process} defined as $\xi_i (t) = (\xi^p_i (t) ,\xi^q_i (t) , \xi^I_i (t) ) \in \RR^3$. 
We define the \emph{full disturbance trajectory} as 
\begin{align}
\xi &= (1, \xi(0), \dots, \xi (T - 1)) \in \RR^{N_\xi}, \label{eq:xi}
\end{align}
where $N_\xi = 1 + 3n T$ and  $\xi (t) = (\xi_1 (t), \dots, \xi_n (t) ) \in \RR^{3n}$ for each time period  $t$.
Note that, in our specification of the disturbance trajectory $\xi$, we have included a constant scalar as its initial component. Such notational convention will prove useful in simplifying the specification of affine control policies in the sequel. \footnote{\rone{
The inclusion of a constant as the initial component of $\xi$ is for notational convenience, as it allows one to represent any affine function of $\xi$ as a \emph{linear} function of $\xi$.}} 

We assume that 
the disturbance trajectory $\xi$ has support $\Xi$ that is a nonempty and  compact subset of $\RR^{N_\xi}$, representable by
\begin{align*}
\Xi = \{ \xi \in \RR^{N_\xi} \ | \ \xi_1 = 1 \ \text{and} \  W \xi \succeq_{\Kcal} 0  \},
\end{align*}
where the matrix $W \in \RR^{\ell \times N_\xi}$ is known. 
It follows from the compactness of $\Xi$ that the second-order moment matrix 
$$\rone{M = \EE\left(\xi \xi'\right),}$$
is finite-valued. We assume, without loss of generality, that $M$ is a positive definite matrix. 
We emphasize that our specification of the disturbance trajectory $\xi$ captures a large family of disturbance processes, including those whose support can be described as the intersection of polytopes and ellipsoids.

\section{Decentralized Control Design} \label{sec:decent_control}
\subsection{State Space Description}
In what follows, we build on the individual resource models developed in Section \ref{sec:model} to develop a discrete-time state space model describing the collective dynamics of  the distribution network. 
We partition the system into $n$ subsystems, where each subsystem $i  \in  \{1, \dots , n\}$ encapsulates the  dynamics of resources connected to bus $i$.
For each subsystem $i$, we let the energy storage state $x_i (t)$ be its \emph{state} at time $t$, and define its \emph{input} according to
\begin{align*}
u_i (t) = \bmat{p^S_i (t) \\ q^I_i (t) }.
\end{align*}
The corresponding state equation for each subsystem $i$ is therefore given by Eq. \eqref{eq:storage}. 
We define the full system state and input at time $t$ by $x(t) = (x_1 (t), .. , x_n(t)) \in \RR^{n}$ and $u(t) = (u_1 (t), .., u_n(t)) \in \RR^{2n}$, respectively. 
The full system state equation admits the following representation
$$x(t+1) = x(t) + B u (t).$$
Here, the matrix $B$ is given by
\begin{align}
\rone{B = I_n \otimes \bmat{-\Delta & 0},} \label{eq:B_mat}
\end{align}
where $\otimes$ denotes the Kronecker product operator. 
The  initial condition\footnote{Recall that the initial condition $x(0)$ is assumed fixed and known.} and  system  trajectories are related according to
\begin{align}
x = \Ambb x(0) + \Bmbb u , \label{eq:x_full}
\end{align}
where $x$ and $u$ represent the \emph{state} and \emph{input trajectories}, respectively. They are given by
\begin{align*}
x &= (x(0), \dots, x(T)) \in \RR^{N_x} , && N_x = n (T+1) ,\\
u &= (u(0), \dots, u(T-1)) \in \RR^{N_u}, && N_u = 2nT.
\end{align*}
We refer the reader to  Appendix \ref{app:matrix} for a precise specification of the matrices $(\Ambb, \Bmbb)$.

\subsection{Decentralized Control Design}
The controller information structure considered in this paper is such that  each subsystem is required to determine its local control input using only  its local measurements. We therefore restrict ourselves to fully decentralized disturbance-feedback control policies.\footnote{For simplicity of exposition, it is assumed in this paper that  each subsystem can  perfectly observe its local disturbance process. We note, however, that all of the results presented in this paper can be immediately  generalized to the setting in which each subsystem has only partial linear observations of its local disturbance process. We refer the readers to \cite{Lin2016} for a  detailed treatment of such systems.} That is to say, at each  time  $t$, the control input to each subsystem $i$ is restricted to be of the form
\begin{align*}
u_i (t) = \gamma_i (\xi_i^t, t),
\end{align*}
where $\gamma_i (\cdot, t)$ is a causal measurable function of the local disturbance history. 
We define the \emph{local control policy} for subsystem $i$ as $\gamma_i = ( \gamma_i ( \cdot, 0), .., \gamma_i (\cdot, T-1) )$; and   refer to the collection of local control policies $\gamma = (\gamma_1, .., \gamma_n)$ as the \emph{decentralized control policy} for the system. Finally,  we define $\Gamma$ to be the \emph{set of all  admissible decentralized control policies}.

In this paper, we consider the objective of minimizing the expected amount of active power required to meet demand over the distribution network. For the setting considered, this is technically equivalent to minimizing the expected active power losses plus the terminal energy stored in the distribution network. In a similar spirit to \cite{baran1989network, turitsyn2011options}, we approximate the active power loss on line $(i, j) \in \Ecal$ at time period $t$ as\footnote{Implicit in this approximation is the assumption that the bus voltage magnitudes are uniform across the network.}
\begin{align*}
\delta p_{ij}(t) = r_{ij} \left( \frac{p_{ij} (t)^2 + q_{ij} (t)^2}{v_0 (t)^2} \right). 
\end{align*}
By a direct substitution of the linearized branch flow Eqs. \eqref{eq:branch_lin1}-\eqref{eq:branch_lin2} into the above approximation, one can represent the total active power losses as a convex quadratic function in the input trajectory $u$ and disturbance  trajectory $\xi$. Namely, one can construct matrices $L_u \in \RR^{2nT \times N_u}$, $L_\xi \in \RR^{2nT \times N_\xi}$, and a positive definite diagonal matrix $\Sigma \in \RR^{2nT \times 2nT}$, such that the total active power losses can be written as
\begin{align}
  \sum_{t = 0}^{T-1} \sum_{(i,j) \in \Ecal} \delta p_{ij} (t) = (L_u u + L_\xi \xi) ' \Sigma (L_u u + L_\xi \xi). \label{eq:loss}
\end{align}
In addition, the sum of the terminal energy storage states across the network can be written as a linear function of the state trajectory $x$. Namely, it is straightforward to construct a vector $c \in \RR^{N_x}$, such that
\begin{align}
 \sum_{i=1}^n x_i (T) = c'x. \label{eq:terminal_storage}
\end{align}
We refer the reader to Appendix \ref{app:matrix} for an exact specification of  $c$, $L_u, \ L_\xi$, and $\Sigma$. Henceforth, we define the expected cost associated with a decentralized control policy $\gamma \in \Gamma$ according to 
\begin{align}
J (\gamma) = \EE^\gamma \left[ c'x + (L_u u + L_\xi \xi) ' \Sigma (L_u u + L_\xi \xi) \right]. \label{eq:cost}
\end{align}
Here, expectation is taken with respect to the joint distribution on $(x, u, \xi)$ induced by the control  policy $\gamma$.

We define the \emph{\rone{decentralized control design}} problem as
\begin{equation}
\begin{alignedat}{8}
&\text{minimize}  \ \ &&  J (\gamma) \\
&\text{subject to}  \ \ &&  \gamma \in \Gamma \\
&&&\hspace{-.105in} \left. \begin{array}{l}
x  \in \Xcal, \ u  \in \Ucal (\xi ) \\ 
x = \Ambb x(0) + \Bmbb u \\
u = \gamma (\xi) 
\end{array}  \right\} \forall \xi \in \Xi,
\end{alignedat} \label{opt:decent_full}
\end{equation}
\rone{where the decision variable is the decentralized control policy $\gamma \in \Gamma$.}
The set of feasible states $\Xcal$ is defined according to inequality \eqref{eq:SoC_con}. The set of feasible control inputs $\Ucal (\xi)$ is defined according to inequalities \eqref{eq:v_con}, \eqref{eq:pS_con}, and \eqref{eq:q_I}. We let $J^{*}$ denote the \emph{optimal value} of problem \eqref{opt:decent_full}.

\section{Design of Affine  Controllers} \label{sec:affine}

The decentralized control design problem \eqref{opt:decent_full} amounts to an infinite-dimensional convex program, and is, in general, computationally intractable. We therefore resort to approximation by restricting the space of admissible decentralized control policies to be causal affine functions of the measured disturbance process. In addition, we approximate the feasible region of problem \eqref{opt:decent_full} from within by a polyhedral set. The combination of these two approximations enables the  computation of a decentralized control policy, which is guaranteed to be feasible for problem \eqref{opt:decent_full}, through solution of a finite-dimensional conic program.

\subsection{Polyhedral Inner Approximation of Constraints}
\begin{figure}[htb]
\centering

\newcommand*{\ra}{0.8}

\begin{subfigure}[t]{0.32 \linewidth}
\centering
\begin{tikzpicture}

\tikzstyle{blank} = [circle, inner sep = 0pt, minimum size = 0mm]
\tikzstyle{Line} = [line width = 0.35pt]
\tikzstyle{Line_dashed} = [line width = 0.35pt, densely dashed]
\tikzstyle{Arrows} = [-angle 60, line width = 0.35pt]
\tikzstyle{circ} = [circle, draw, line width = 0.35pt, inner sep = 0pt, minimum size = 1.5cm]
\tikzstyle{textbox} = [below, text centered]

\draw [Line] (-2,{-1 * \ra}) arc[radius = \ra, start angle=-90, end angle=90];

\node (center1) at (-2, 0) [blank] {};

\node (left1) at ([xshift = -8mm, yshift = 0mm]center1) [blank] {};
\node (right1) at ([xshift = 12mm, yshift = 0mm]center1) [blank] {};

\node at ([xshift = 3mm, yshift = 1mm]right1) [textbox] {\footnotesize{$\xi_i^I (t)$}};

\node (down1) at ([xshift = 0mm, yshift = -16mm]center1) [blank] {};
\node (up1) at ([xshift = 0mm, yshift = 16mm]center1) [circle, inner sep = 0pt, minimum size = 0mm, label = east: \footnotesize{$q_i^I (t)$}] {};

\filldraw[fill = gray!20, draw = black, line width = 0.35pt] (-2,{-1 * \ra}) arc[radius = \ra, start angle=-90, end angle=-30] -- ({-2 + \ra * cos(30)}, {\ra * sin(30)}) arc[radius = \ra, start angle = 30, end angle = 90] -- cycle;

\draw [Arrows] (left1) -- (right1); 
\draw [Arrows] (down1) -- (up1);

\node at (-2, {\ra * sin(90)}) [circle, inner sep = 0pt, minimum size = 0mm, label = west: \footnotesize{$s_i^I$}] {};

\node at (-2, {- \ra * sin(90)}) [circle, inner sep = 0pt, minimum size = 0mm, label = west: \footnotesize{$-s_i^I$}] {};

\end{tikzpicture}
\caption{}
\label{fig:original}
\end{subfigure}
\begin{subfigure}[t]{0.32 \linewidth}
\centering
\begin{tikzpicture}

\tikzstyle{blank} = [circle, inner sep = 0pt, minimum size = 0mm]
\tikzstyle{Line} = [line width = 0.35pt]
\tikzstyle{Line_dashed} = [line width = 0.35pt, densely dashed]
\tikzstyle{Arrows} = [-angle 60, line width = 0.35pt]
\tikzstyle{circ} = [circle, draw, line width = 0.35pt, inner sep = 0pt, minimum size = 1.5cm]
\tikzstyle{textbox} = [below, text centered]

\draw [Line] (-2,-\ra) arc[radius = \ra, start angle=-90, end angle=90];

\node (center1) at (-2, 0) [blank] {};

\node (left1) at ([xshift = -8mm, yshift = 0mm]center1) [blank] {};
\node (right1) at ([xshift = 12mm, yshift = 0mm]center1) [blank] {};

\node at ([xshift = 3mm, yshift = 1mm]right1) [textbox] {\footnotesize{$\xi_i^I (t)$}};

\node (down1) at ([xshift = 0mm, yshift = -16mm]center1) [blank] {};
\node (up1) at ([xshift = 0mm, yshift = 16mm]center1) [circle, inner sep = 0pt, minimum size = 0mm, label = east: \footnotesize{$q_i^I (t)$}] {};

\filldraw[fill = gray!20, draw = black, line width = 0.35pt] (-2,{-\ra * sin(30)}) -- ({-2 + \ra * cos(30)}, {-\ra * sin(30)}) -- ({-2 + \ra * cos(30)}, {\ra * sin(30)}) -- (-2,{\ra * sin(30)}) -- cycle;

\draw [Arrows] (left1) -- (right1); 
\draw [Arrows] (down1) -- (up1);

\node at (-2, {\ra * sin(30)}) [circle, inner sep = 0pt, minimum size = 0mm, label = west: \footnotesize{$\overline{q}_i^I (t)$}] {};

\node at (-2, {-\ra * sin(30)}) [circle, inner sep = 0pt, minimum size = 0mm, label = west: \footnotesize{$-\overline{q}_i^I (t)$}] {};

\end{tikzpicture}
\caption{}
\label{fig:inner_approx}
\end{subfigure}
\begin{subfigure}[t]{0.32 \linewidth}
\centering
\begin{tikzpicture}

\tikzstyle{blank} = [circle, inner sep = 0pt, minimum size = 0mm]
\tikzstyle{Line} = [line width = 0.35pt]
\tikzstyle{Line_dashed} = [line width = 0.35pt, densely dashed]
\tikzstyle{Arrows} = [-angle 60, line width = 0.35pt]
\tikzstyle{circ} = [circle, draw, line width = 0.35pt, inner sep = 0pt, minimum size = 1.5cm]
\tikzstyle{textbox} = [below, text centered]

\node (center1) at (-2, 0) [blank] {};

\node (left1) at ([xshift = -8mm, yshift = 0mm]center1) [blank] {};
\node (right1) at ([xshift = 15mm, yshift = 0mm]center1) [blank] {};

\node at ([xshift = 3mm, yshift = 1mm]right1) [textbox] {\footnotesize{$\xi_i^I (t)$}};

\node (down1) at ([xshift = 0mm, yshift = -16mm]center1) [blank] {};
\node (up1) at ([xshift = 0mm, yshift = 16mm]center1) [circle, inner sep = 0pt, minimum size = 0mm, label = east: \footnotesize{$q_i^I (t)$}] {};

\filldraw[fill = gray!20, draw = black, line width = 0.35pt] (-2,{-\ra * sqrt(2)}) -- ( {-2+\ra * cos(30)}, {( \ra * cos(30)) - \ra * sqrt(2)}) -- ({-2+\ra * cos(30)}, {(- \ra *  cos(30)) + \ra * sqrt(2)}) -- (-2, {\ra * sqrt(2)}) -- cycle;

\draw [Line_dashed] ( {-2+ \ra * cos(30)}, {( \ra * cos(30)) - \ra * sqrt(2)}) -- ({-2 + \ra * sqrt(2)}, 0) -- ({-2+\ra * cos(30)}, {(- \ra * cos(30)) + \ra * sqrt(2)});

\draw [Arrows] (left1) -- (right1); 
\draw [Arrows] (down1) -- (up1);

\draw [Line] (-2,-\ra) arc[radius = \ra, start angle=-90, end angle=90];

\newcommand*{\eps}{0.05}
\newcommand*{\epsy}{0}

\node at ({-2 + \eps} , {\ra * sqrt(2) - \epsy}) [circle, inner sep = 0pt, minimum size = 0mm, label = west: \footnotesize{$\sqrt{2} s_i^I$}] {};

\node at ({-2 + \eps} , {-\ra * sqrt(2) + \epsy}) [circle, inner sep = 0pt, minimum size = 0mm, label = west: \footnotesize{$-\sqrt{2} s_i^I$}] {};

\node at ([xshift = -3mm, yshift = 7.5mm]right1) [textbox] {\footnotesize{$\sqrt{2} s_i^I$}};
\end{tikzpicture}
\caption{} \label{fig:outer_approx}
\end{subfigure}

\caption{\rone{The above plots depict an inverter's  range of feasible reactive power injections (in gray) at a particular time period $t$ as specified by (a) the \emph{original} quadratic constraints \eqref{eq:q_I}, (b) the \emph{inner} linear constraints \eqref{eq:q_I_inner}, and (c) the \emph{outer} linear constraints  \eqref{eq:q_I_outer}.}}
\label{fig:approx}
\end{figure}
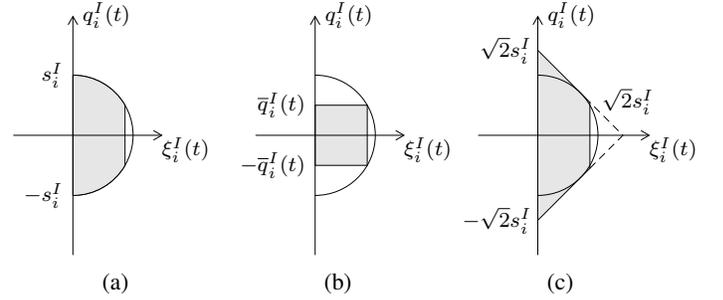

The feasible state space $\Xcal$ is clearly polyhedral. The feasible input space $\Ucal(\xi)$ is not. It can, however, be approximated from within by a polyhedral set  by replacing the quadratic constraint in \eqref{eq:q_I} with the following pair of  linear constraints:
\begin{align}
\left| q_i^I (t) \right| \leq  \overline{q}_i^I (t). \label{eq:q_I_inner}
\end{align}
Here, the deterministic constant $\overline{q}_i^I (t)$ is defined according to
\begin{align*}
\overline{q}_i^I (t) = \inf \left\{ \left. \sqrt{ \left. s_i^I \right. ^2 -  \left.\xi_i^I (t) \right. ^2} \; \right| \; \xi \in \Xi \right\}. 
\end{align*}
Essentially, $\overline{q}_i^I (t)$ denotes the minimum reactive power capacity that is guaranteed to be available at inverter $i$ at time $t$. \rone{We provide a graphical illustration of this polyhedral inner approximation in Figure \ref{fig:inner_approx}.}

Although an inner approximation of this form may appear conservative at first glance, several recent studies \cite{turitsyn2011options, kekatos2015stochastic} have observed such approximations to result in a small loss of performance, as measured by the objective function considered in this paper. We corroborate these claims in Section \ref{sec:LB} by developing a technique to bound the loss of optimality incurred by this inner approximation. 
In particular, the suboptimality incurred by such an approximation is  shown to be  small for the case study considered in this paper.

Inequalities \eqref{eq:v_con}, \eqref{eq:SoC_con}, \eqref{eq:pS_con}, and \eqref{eq:q_I_inner} define a collection of $m = 8nT$ linear constraints on the state, input, and disturbance trajectories. We represent them more succinctly as
\begin{align*}
\underline{F_x} x + \underline{F_u} u + \underline{F_\xi} \xi \leq 0, \quad \forall  \ \xi \in \Xi,
\end{align*}
where it is straightforward to construct the  matrices $\underline{F_x} \in \RR^{m \times N_x}, \underline{F_u} \in \RR^{m \times N_u}$, and $\underline{F_\xi} \in \RR^{m \times N_\xi}$ using the given problem data.
The following optimization problem is an inner approximation to the original decentralized control design problem \eqref{opt:decent_full}:
\begin{align}
\nonumber \text{minimize}  \ \ &  \EE^\gamma \left[ c' x + (L_u u + L_\xi \xi) ' \Sigma (L_u u + L_\xi \xi) \right]\\
\nonumber \text{subject to}  \ \ &  \gamma \in \Gamma \\
&  \hspace{-.105in} \left. \begin{array}{l}
\underline{F_x} x + \underline{F_u} u + \underline{F_\xi} \xi \leq 0 \\
 x  = \Ambb x(0) + \Bmbb u  \\ 
 u = \gamma(\xi) \\ 
\end{array} \right\} \forall \  \xi \in \Xi, \label{opt:decent_linear_con}
\end{align}
\rone{where the decision variable is given by $\gamma$.}
Although convex,
problem \eqref{opt:decent_linear_con} is an infinite-dimensional program, and is therefore computationally intractable, in general.
In what follows, we  refine this approximation by further restricting the space of admissible controllers to be \emph{affine} functions of the disturbance trajectory.

\subsection{Affine Control Design via Conic Programming}

\rone{We restrict our attention  to decentralized affine control policies of the form
\begin{align}
u_i (t) = \overline{u}_i (t) + \sum_{s=0}^t Q_{i} (t,s ) \xi_i (s) \label{eq:affine_policy}
\end{align}
for each subsystem $i = 1, \dots, n$ and time $t = 0, \dots, T-1$. Here, $\overline{u}_i (t) \in \RR^{2}$ denotes the open loop component of the local control, and $(Q_{i}(t,0), \dots, Q_{i}(t,t))$ the collection of feedback control gains at time $t$.
One can write the decentralized affine control policy in \eqref{eq:affine_policy} more compactly as
\begin{align*} 
u = Q \xi.
\end{align*}
We enforce the desired information structure in $Q$ by requiring that  $Q \in S$, where $S$ denotes the subspace of matrices that respect the information structure associated with the set of admissible decentralized control policies $\Gamma$. Specifically, 
\begin{align*}
S = \left\{ \left. Q \in \RR^{N_u \times N_\xi} \right| Q \in \Gamma \right\}.
\end{align*}}

The restriction to decentralized affine control policies gives rise to the following semi-infinite program,\footnote{A semi-infinite program is an optimization problem involving finitely many decision variables, and an infinite number of constraints.} which stands as a more conservative inner approximation to the original decentralized control design problem \eqref{opt:decent_full}.
\begin{align}
\nonumber \text{minimize}  \ \ &  \EE  \left[ c' x + (L_u u + L_\xi \xi) ' \Sigma (L_u u + L_\xi \xi) \right]\\
\nonumber \text{subject to}  \ \ &  Q \in S \\
&  \hspace{-.105in} \left. \begin{array}{l}
\underline{F_x} x + \underline{F_u} u + \underline{F_\xi} \xi \leq 0 \\
 x  = \Ambb x(0) + \Bmbb u  \\ 
 u = Q \xi \\ 
\end{array} \right\} \forall \  \xi \in \Xi, \label{opt:decent_aff}
\end{align}
\rone{where the decision variable is given by $Q$.}
Given our assumption that the uncertainty set $\Xi$ has a conic representation, one can directly apply a previous result from \cite{Lin2016} to equivalently reformulate the semi-infinite program \eqref{opt:decent_aff} as a finite-dimensional conic program.\footnote{We note that a direct application of Proposition 3 from \cite{Lin2016} also requires that the information structure of the underlying decentralized control problem  be \emph{partially nested}. This condition requiring partial nestedness of the information structure is trivially satisfied for the decentralized control design problem \eqref{opt:decent_full} under consideration in this paper.} In Proposition \ref{prop:primal_affine_conic}, we present the finite-dimensional conic reformulation of the semi-infinite program \eqref{opt:decent_aff} implied by \cite{Lin2016}[Prop. 3].

\begin{propositio} \label{prop:primal_affine_conic}
The semi-infinite program \eqref{opt:decent_aff} 
admits an equivalent reformulation as the following finite-dimensional
conic  program
\begin{equation}
\begin{alignedat}{8}
&\text{minimize} && {\rm Tr} \bigg( \Big(  Q' L_u' \Sigma  L_u Q + \left( 2 L_\xi' \Sigma L_u + e_1 c' \Bmbb \right) Q  \\
&&& \qquad +  L_\xi' \Sigma L_\xi  \Big) M \bigg) + c' \Ambb x(0)\\
&\text{subject to} \quad && Q \in S  \\
&&& Z \in \RR^{m \times N_\xi},  \  \ \Pi \in \RR^{ \ell \times m}, \ \  \nu \in \RR^m_+ \\
&&& (\underline{F_u} + \underline{F_x} \Bmbb ) Q + \underline{F_x} \Ambb x(0) e_1' +  \underline{F_\xi} + Z = 0, \\
&&& Z = \nu e_1' +  \Pi' W, \\
&&& \Pi \succeq_{\Kcal^*}  0,
\end{alignedat} \label{eqn:primal_affine_PN_conic}
\end{equation}
\rone{where the decision variables are  given by $Q$, $Z$,  $\Pi$, and $\nu$.}
Let $J^{\rm{in}}$  denote the optimal value of the above program. It stands as an \emph{upper bound} on the optimal value of the original decentralized control problem \eqref{opt:decent_full}, i.e., $J^{*} \leq J^{\rm in}$. 
\end{propositio}

\rone{Several comments are in order. First, the  specification of the conic program \eqref{eqn:primal_affine_PN_conic}  relies on the probability distribution of the disturbance  $\xi$ only through its support $\Xi$ and second-order moment matrix $M$. Second, the conic program can be efficiently solved for a variety of cones $\Kcal$,  including polyhedral and second-order cones.
For such cones,  problem \eqref{eqn:primal_affine_PN_conic} amounts to a conic program with $O(n^2 T^2)$ decision variables and $O(n^2 T^2)$ constraints. It can thus be solved in time that is polynomial in the control horizon $T$ and the number of subsystems $n$. Finally,  assuming that the decentralized affine controller $Q^*$ is computed at a central location, the decentralized implementation of the controller will require  the communication of each local control policy to its corresponding subsystem. This entails the transmission of $3T^2 + 5T$ real numbers to each subsystem.}

\rone{
\begin{remar}[Fast Time-Scale  Implementation]
In practice, the active power generated by a photovoltaic resource may fluctuate over time-scales (e.g., seconds to minutes) that are substantially shorter than the time-scale being used for control design (e.g., hourly).  In Appendix \ref{app:continuous_time}, we propose a method to enable the implementation of controllers designed according to Proposition \ref{prop:primal_affine_conic} over more finely grained time-scales. Under a mild assumption on the quasi-stationarity of the support of the underlying disturbance trajectory, the proposed fast time-scale implementation of the controller is shown to yield state, input, and voltage magnitude trajectories, which are guaranteed to be feasible on this more finely grained time-scale. 
\end{remar}
}

\section{Design of Performance Bounds} \label{sec:LB}
The restriction to affine policies computed according to Proposition \ref{prop:primal_affine_conic} may result in the  loss of optimality with respect to the original decentralized control design problem. In this section, we develop a tractable method to bound this loss of optimality via the solution of a conic programming relaxation -- the optimal value of which is guaranteed to stand as a lower bound on the optimal value of the  original decentralized control design problem \eqref{opt:decent_full}. 
With such  a lower bound in hand, one can estimate the suboptimality incurred by any feasible decentralized control policy.

\subsection{Polyhedral Outer Approximation of Constraints}
As an initial step in the derivation of this relaxation, we construct a polyhedral outer approximation of the feasible region of problem \eqref{opt:decent_full}. Specifically, the quadratic constraint  in \eqref{eq:q_I} can be relaxed to the following pair of  linear constraints:
\begin{align}
\left| q_i^I (t) \right| \leq \sqrt{2} s_i^I - \xi_i^I (t). \label{eq:q_I_outer}
\end{align}
We provide a graphical illustration of this polyhedral outer approximation in Figure \ref{fig:outer_approx}.

\rone{Inequalities \eqref{eq:v_con}, \eqref{eq:SoC_con}, \eqref{eq:pS_con}, and \eqref{eq:q_I_outer}} 
define a collection of $m$ linear constraints on the state, input, and disturbance trajectories. We represent them more succinctly as
\begin{align*}
\overline{F_x} x + \overline{F_u} u + \overline{F_\xi} \xi \leq 0, \quad \forall \ \xi \in \Xi,
\end{align*}
where it is straightforward to construct the matrices $\overline{F_x} \in \RR^{m \times N_x}, \overline{F_u} \in \RR^{m \times N_u}$, and  $\overline{F_\xi} \in \RR^{m \times N_\xi}$  using the given problem data. 
\rone{
The following optimization problem is an outer approximation to the original decentralized control design problem \eqref{opt:decent_full}:
\begin{align}
\nonumber \text{minimize}  \ \ &  \EE^\gamma \left[ c' x + (L_u u + L_\xi \xi) ' \Sigma (L_u u + L_\xi \xi) \right]\\
\nonumber \text{subject to}  \ \ &  \gamma \in \Gamma \\
&  \hspace{-.105in} \left. \begin{array}{l}
\overline{F_x} x + \overline{F_u} u + \overline{F_\xi} \xi \leq 0 \\
 x  = \Ambb x(0) + \Bmbb u  \\ 
 u = \gamma(\xi) \\ 
\end{array} \right\} \forall \  \xi \in \Xi, \label{opt:decent_linear_con_outer}
\end{align}
\rone{where the decision variable is given by $\gamma$.}
}

\subsection{Lower Bounds via Conic Programming}
\rone{Problem \eqref{opt:decent_linear_con_outer} is, in general, computationally intractable due to the infinite-dimensionality of its decision space. It is possible, however, to further relax problem \eqref{opt:decent_linear_con_outer} to a finite-dimensional conic program, under an additional technical assumption on the probability distribution of the disturbance trajectory $\xi$.}

\begin{assumptio}[Disturbance Process] \label{ass:noise} There exist matrices $H_i^t \in \RR^{3n(t + 1) \times (1 + 3(t+1))}$ and $H^t \in \RR^{N_{\xi} \times (1 + 3n(t + 1) )}$  such that
\begin{align*}
\EE \left[ \xi^t \left| \xi_i^t \right. \right] = H_i^t \bmat{1 \\ \xi_i^t}  \ \  \text{and } \ \  \EE \left[ \xi \left| \xi^t \right. \right] = H^t \bmat{1 \\ \xi^t} 
\end{align*}
almost surely, for all time periods $t=0, \dots, T-1$ and subsystems $i = 1, \dots, n$.
\end{assumptio}

Although Assumption \ref{ass:noise} may appear restrictive, it was shown in \cite{Hadjiyiannis2011} to hold for a large family of distributions.  In particular, Assumption \ref{ass:noise} holds for all disturbance processes, which possess elliptically contoured  distributions. These include the multivariate Gaussian distribution,  the multivariate $t$-distribution, their truncated versions, and uniform distributions over ellipsoids \cite{Cambanis1981}.  It is also straightforward to show that Assumption \ref{ass:noise} is satisfied by any disturbance process for which the random vectors $\xi_i(t)$  ($i = 1, \dots, n$, $t = 0, \dots, T-1$) are mutually independent.

With Assumption \ref{ass:noise} in hand, a direct application of Proposition 4 in  \cite{Lin2016} yields a conic programming relaxation\footnote{\rone{The conic programming relaxation \eqref{eqn:dual_affine_PN} is constructed by approximating the \emph{Lagrangian dual problem} of  \eqref{opt:decent_linear_con_outer}. We refer the readers to \cite{Lin2016, Hadjiyiannis2011} for the technical details of this derivation.}} of problem \eqref{opt:decent_linear_con_outer}. Its optimal value stands as a lower bound 
on the optimal value of the original decentralized control design problem \eqref{opt:decent_full}.

\begin{propositio} \label{prop:dual_affine_PN}  Consider the following finite-dimensional conic program:
\begin{equation}
\begin{alignedat}{8}
&\text{minimize} && {\rm Tr} \bigg( \Big(  Q' L_u' \Sigma  L_u Q + \left( 2 L_\xi' \Sigma L_u + e_1 c' \Bmbb \right) Q  \\
&&& \qquad +  L_\xi' \Sigma L_\xi  \Big) M \bigg) + c' \Ambb x(0) \\
&\text{subject to} \quad && Q \in S , \ \ Z \in \RR^{m \times N_\xi }\\
&&& (\overline{F_u} + \overline{F_x} \Bmbb) Q  + \overline{F_x} \Ambb x(0) e_1'  + \overline{F_\xi} + Z = 0, \\
&&& W M Z' \succeq_{\Kcal} 0,  \\
&&& e_1' M Z' \geq 0,
\end{alignedat} \label{eqn:dual_affine_PN}
\end{equation}
where the decision variables are given by $Q$ and $Z$. Let  $J^{\rm out}$  denote the optimal value of the above program. If Assumption \ref{ass:noise} holds, then  $J^{\rm out} \leq J^*$. 
\end{propositio}

\rone{Given Assumption \ref{ass:noise},  the conic program \eqref{eqn:dual_affine_PN} can be used to evaluate the performance of any feasible control policy. Namely, a policy $\gamma \in \Gamma$ is close to optimal (for a given problem instance) if $J(\gamma)$ is close to $J^{\rm out}$. Additionally, Propositions \ref{prop:primal_affine_conic} and \ref{prop:dual_affine_PN} imply that the optimal value of the original decentralized control problem \eqref{opt:decent_full} satisfies $$J^{\rm out} \leq J^* \leq J^{\rm in}.$$
Therefore, a small gap between $J^{\rm in}$ and $J^{\rm out}$  implies that decentralized affine control policies are  close to optimal for the underlying problem instance.
Finally, we note that the conic program  \eqref{eqn:dual_affine_PN} can be efficiently solved for a variety of cones $\Kcal$, including polyhedral and second-order cones. 
For such cones,  problem \eqref{eqn:dual_affine_PN} amounts to a conic program with $O(nT^2)$ decision variables and $O(nT)$ constraints. It can thus be solved in time that is polynomial in the control horizon $T$ and the number of subsystems $n$.}

\section{Case Study} \label{sec:case_study}

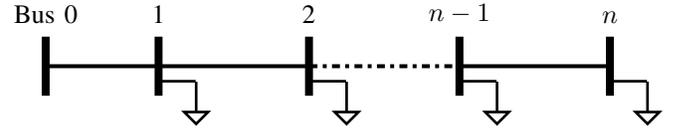
\begin{figure}[http]
\centering

\tikzstyle{Bus} = [line width = 3pt]
\tikzstyle{Line} = [line width = 1.5pt]
\tikzstyle{Load} = [-open triangle 90, line width = 1pt]
\tikzstyle{Lines} = [line width = 1.5pt, dashdotted]
\tikzstyle{blank} = [circle, inner sep = 0pt, minimum size = 0mm]

\begin{tikzpicture}

\node (b0_up) at (-3.5, -0.4) [blank] {};
\node (b0_down) at (-3.5, 0.4) [blank] {};
\node (b0_mid) at (-3.5, 0) [blank] {};

\node (b1_up) at (-2, -0.4) [blank] {};
\node (b1_down) at (-2, 0.4) [blank] {};
\node (b1_mid) at (-2, 0) [blank] {};

\node (b2_up) at (0, -0.4) [blank] {};
\node (b2_down) at (0, 0.4) [blank] {};
\node (b2_mid) at (0, 0) [blank] {};

\node (bn-1_up) at (2, -0.4) [blank] {};
\node (bn-1_down) at (2, 0.4) [blank] {};
\node (bn-1_mid) at (2, 0) [blank] {};

\node (bn_up) at (4, -0.4) [blank] {};
\node (bn_down) at (4, 0.4) [blank] {};
\node (bn_mid) at (4, 0) [blank] {};

\node (load1_1) at (-2, -0.2) [blank] {};
\node (load1_2) at (-1.5, -0.2) [blank] {};
\node (load1_3) at (-1.5, -0.8) [blank] {};

\node (load2_1) at (0, -0.2) [blank] {};
\node (load2_2) at (0.5, -0.2) [blank] {};
\node (load2_3) at (0.5, -0.8) [blank] {};

\node (loadn-1_1) at (2, -0.2) [blank] {};
\node (loadn-1_2) at (2.5, -0.2) [blank] {};
\node (loadn-1_3) at (2.5, -0.8) [blank] {};

\node (loadn_1) at (4, -0.2) [blank] {};
\node (loadn_2) at (4.5, -0.2) [blank] {};
\node (loadn_3) at (4.5, -0.8) [blank] {};

\draw [Bus] (b0_up) -- (b0_down) 
	node [above, text centered] {Bus 0};
\draw [Bus] (b1_up) -- (b1_down)
	node [above, text centered] {1};
\draw [Bus] (b2_up) -- (b2_down)
	node [above, text centered] {2};
\draw [Bus] (bn-1_up) -- (bn-1_down)
	node [above, text centered] {$n-1$};
\draw [Bus] (bn_up) -- (bn_down)
	node [above, text centered] {$n$};

\draw [Line] (b0_mid) -- (b1_mid);
\draw [Line] (b1_mid) -- (b2_mid);
\draw [Lines] (b2_mid) -- (bn-1_mid);
\draw [Line] (bn-1_mid) -- (bn_mid);

\draw [Load] (load1_1) -- (load1_2) -- (load1_3);
\draw [Load] (load2_1) -- (load2_2) -- (load2_3);
\draw [Load] (loadn-1_1) -- (loadn-1_2) -- (loadn-1_3);
\draw [Load] (loadn_1) -- (loadn_2) -- (loadn_3);

\end{tikzpicture}
\caption{Schematic diagram of a 12 kV radial distribution feeder with $n+1$ buses.}
\label{fig:feeder}
\end{figure}

We consider the control of distributed energy resources in a 12 kV radial distribution feeder depicted in Fig. \ref{fig:feeder}. The distribution feeder considered in this paper is similar in structure to the network considered in \cite{zhu2015fast}. Apart from the substation, the distribution feeder consists of $n = 14$ buses. We operate the system over a finite time horizon of $T = 24$ hours, beginning at  twelve o'clock (midnight).

\subsection{System Description}
\begin{figure}[htbp]
\centering
\includegraphics[width = 0.85 \linewidth]{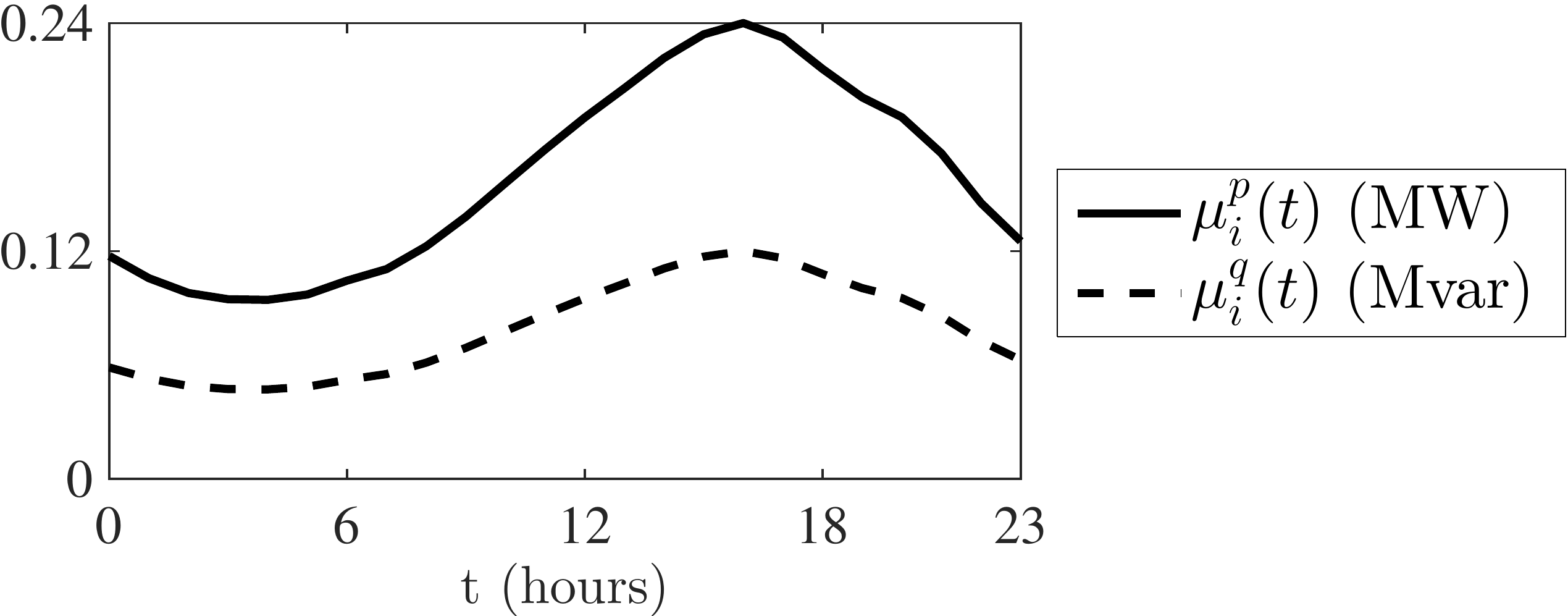}
\caption{Buses $i$ = 3, 4, 5, 13, and  14 are assumed to be identical in terms of their mean load trajectories. The above figure depicts the mean active power and reactive power demand trajectories at these buses. Both trajectories are scaled versions of the load profile DOM-S/M on 07/01/2016 from Southern California Edison \cite{SCE2016}.}
\label{fig:load}
\end{figure}

We assume that only buses 4 and 8 have storage devices and PV inverters installed. All PV inverters are assumed to have an identical active power capacity, which we denote by
by $\theta$ (MW). 
As for demand,  we assume that only buses 3, 4, 5, 13, and 14 have loads; and these loads are assumed to have identical distributions. 
We specify their mean active and reactive power trajectories  in Fig. \ref{fig:load}. 
In order to ensure that Assumption  \ref{ass:noise} is satisfied,  we assume that the random vectors $\xi_i(t)$  ($i = 1, \dots, n$, $t = 0, \dots, T-1$) are mutually independent. In addition, we assume that the random variables $\xi_i^p (t)$, $\xi_i^q(t)$, and $\xi_i^I (t)$ are mutually independent for each bus $i$ and time $t$. \rone{Recall that Assumption \ref{ass:noise} is necessary only for the calculation of the performance bound  specified in Proposition \ref{prop:dual_affine_PN}.}

In Table \ref{tab:notations}, we present additional notation pertinent to this section.
In Table \ref{tab:problem_data}, we specify the parameter values of the distribution network, storage devices, PV inverters, and load.

\newcommand{\tabincell}[2]{\begin{tabular}{@{}#1@{}}#2\end{tabular}}
\begin{table}[htb]
\centering 
\caption{Additional notation.}
\label{tab:notations}
\makegapedcells
\begin{tabular}{ll}
\toprule
\textbf{Notation} & \textbf{Description} \\
\midrule
$\theta$ & Active power capacity of each PV inverter. \\
$\mu_i^p (t)$ & Mean active power demand at bus $i$ and time $t$  . \\
$\mu_i^q (t)$ & Mean reactive power demand at bus $i$ and time $t$. \\
$\mu_i^I (t)$ & Mean active power supply from PV inverter $i$ at time $t$. \\
$\text{Uni} \, [a,b]$ & Uniform distribution on $[a, b]$. \\
\bottomrule
\end{tabular}
\end{table}

\begin{table}[htb]
\centering
\setlength{\tabcolsep}{1pt}
\caption{Specification of system data.}
\label{tab:problem_data}
\makegapedcells
\begin{tabularx}{\linewidth}{ll}
\toprule
\textbf{Distribution network} \\
\midrule
Base voltage magnitude & 12 kV \\
Substation voltage magnitude  & $v_0 = 1$ (per-unit) \\
Impedance on line $(i, j ) \in \Ecal$ &  $r_{ij} = 0.466$, \   $x_{ij} = 0.733$ ($\Omega$) \\
Voltage magnitude constraints & $\underline{v} = 0.95 \cdot  \bone$, \ $\overline{v} = 1.05 \cdot \bone$ (per-unit) \\
\toprule
\textbf{Storage at bus} $i \in \{4, 8\}$ \\
\midrule
Energy capacity & $b_i = 0.5$ (MWh) \\
Power capacity & $\underline{p}^S_i = -0.2$, \ $\overline{p}^S_i = 0.2$ (MW) \\
Initial condition & $x_i(0) = 0$ (MWh)\\
\toprule
\textbf{PV inverter at bus} $i \in \{4, 8\}$ \\
\midrule
Apparent power capacity & $s_i^I = 1.25 \theta$ (MVA) \\
Active power supply & $\xi_i^I (t) \sim \text{Uni} \, [0, 2 \mu_i^I (t) ]$  (MW) \\
Mean active power supply & $\mu_i^I (t) = \theta \cdot \max \left\{ 0.5 \sin \left( \frac{t - 6}{12} \pi \right), 0 \right\}$ \\
\toprule
\textbf{Load at bus} $i \in \{3, 4, 5, 13, 14\}$ \\
\midrule
Active power demand & $\xi_i^p (t) \sim  \text{Uni} \, [ 0.7 \mu_i^p (t), \; 1.3 \mu_i^p (t)  ]$   (MW) \\
Reactive power demand  & $\xi_i^q (t) \sim  \text{Uni} \, [ 0.7 \mu_i^q (t), \; 1.3 \mu_i^q (t) ]$ (Mvar) \\
\bottomrule
\end{tabularx}
\end{table}

\subsection{Discussion}
\begin{figure}[htb]
\centering
\includegraphics[width = 0.9\linewidth]{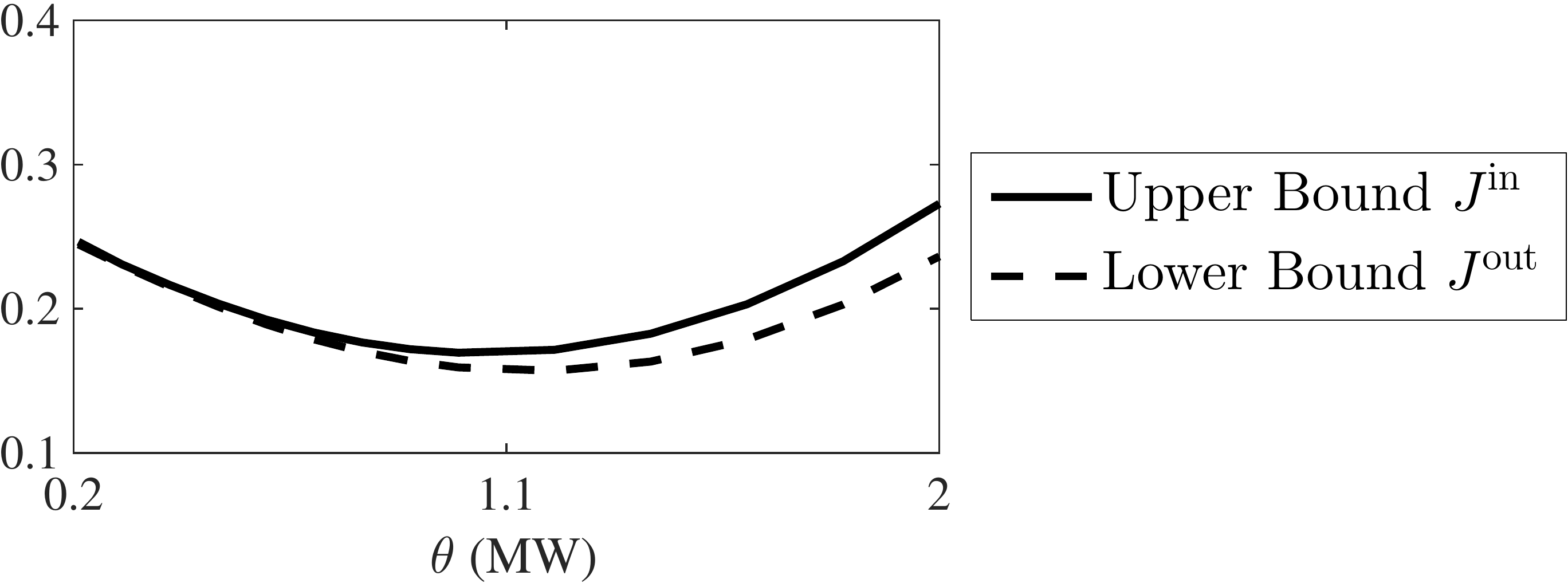}
\caption{This figure depicts the upper  and lower bounds, $J^{\rm in}$ and  $J^{\rm out}$, respectively, on the optimal value of the decentralized control design problem $J^*$ (measured in MWh) as a function of the PV inverter active power capacity $\theta$.
}
\label{fig:para_PV_cap}
\end{figure}

We begin by examining the performance of the decentralized controller proposed in this paper. In Fig. \ref{fig:para_PV_cap}, we plot both the upper and lower  bounds on the optimal value $J^*$ of the decentralized control design problem \eqref{opt:decent_full}, as a function of the PV inverter active power capacity $\theta$. Recall that $J^{\rm in}$ measures the cost incurred by the decentralized affine control policy computed according to Proposition \ref{prop:primal_affine_conic}. Notice that, at low  PV penetration levels (i.e., for low values of $\theta$), the upper and lower bounds nearly coincide. This indicates that the decentralized affine control policy is nearly optimal for the original decentralized control design problem.
More interestingly, at high PV penetration levels (i.e., for high values of $\theta$), the gap between the upper and lower bounds remains small. This reveals that decentralized affine control policies persist in being close to optimal for the system considered, despite the presence of large and unpredictable fluctuations  in PV active power generation. Therefore, for the
system under consideration, there is little additional value to be had in the design of more sophisticated (nonlinear) control
policies.

\begin{figure*}[htbp]

\centering

\begin{subfigure}[t]{0.4557 \linewidth}
\centering
\begin{minipage}[t]{0.49 \linewidth}
\centering
$\xi_4^I (t)$ (MW)

\includegraphics[width = \linewidth]{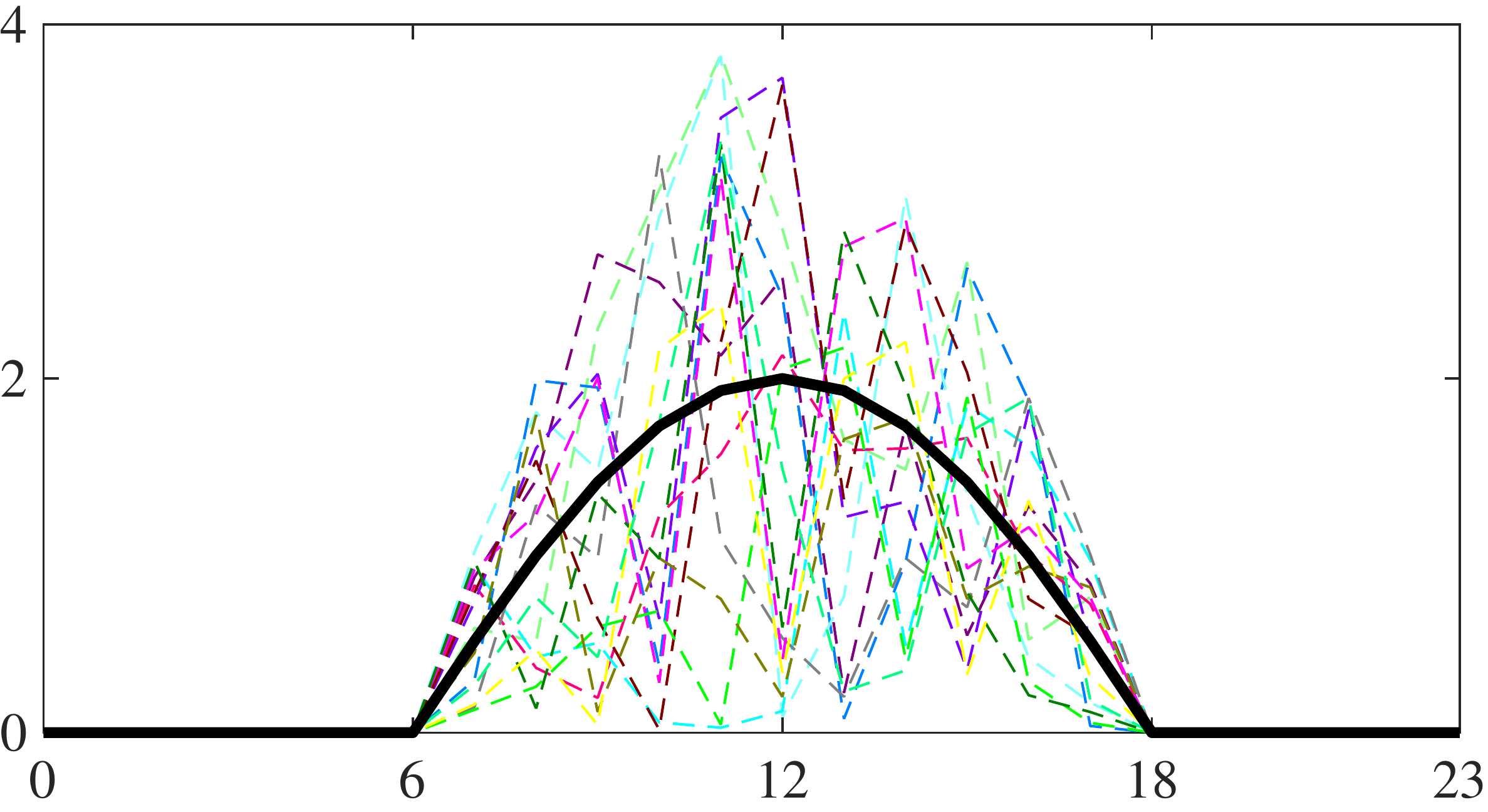}
\end{minipage}
\begin{minipage}[t]{0.49 \linewidth}
\centering
$\xi_4^I (t)$ (MW)

\includegraphics[width = \linewidth]{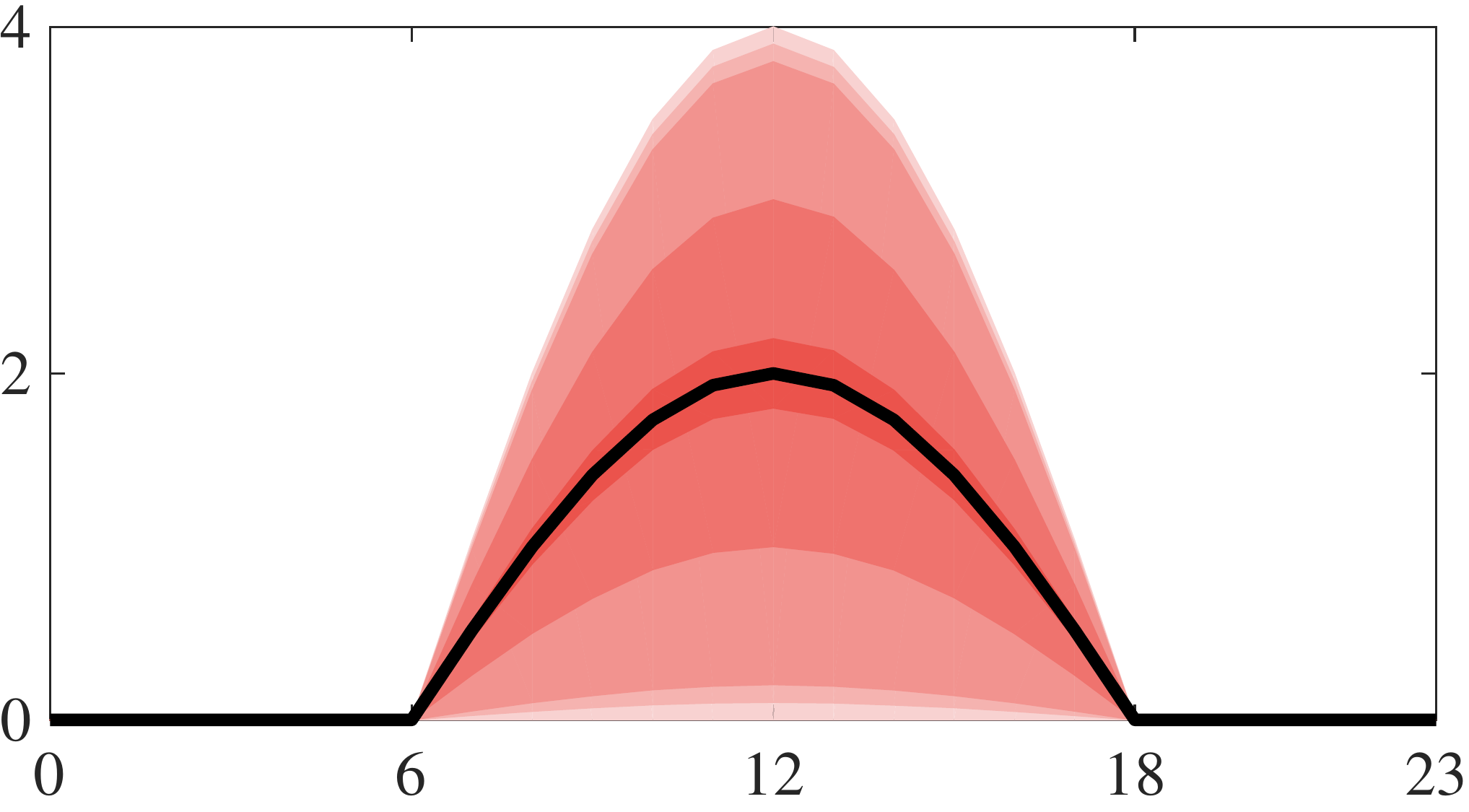}
\end{minipage}
\vspace{.005in}

\begin{minipage}[t]{0.49 \linewidth}
\centering
$q_4^I (t)$ (Mvar)

\includegraphics[width = \linewidth]{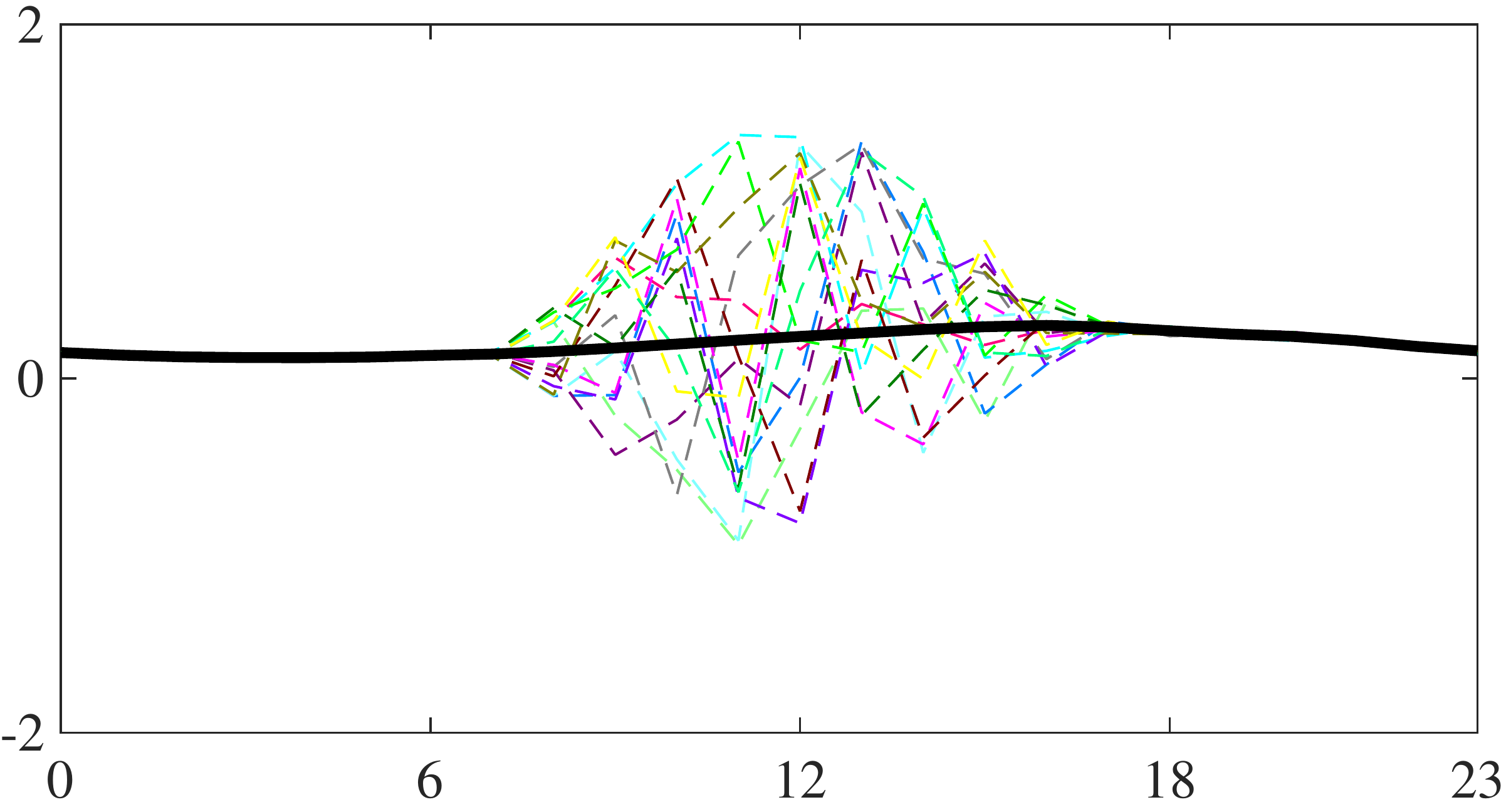}
\end{minipage}
\begin{minipage}[t]{0.49 \linewidth}
\centering
$q_4^I (t)$ (Mvar)

\includegraphics[width = \linewidth]{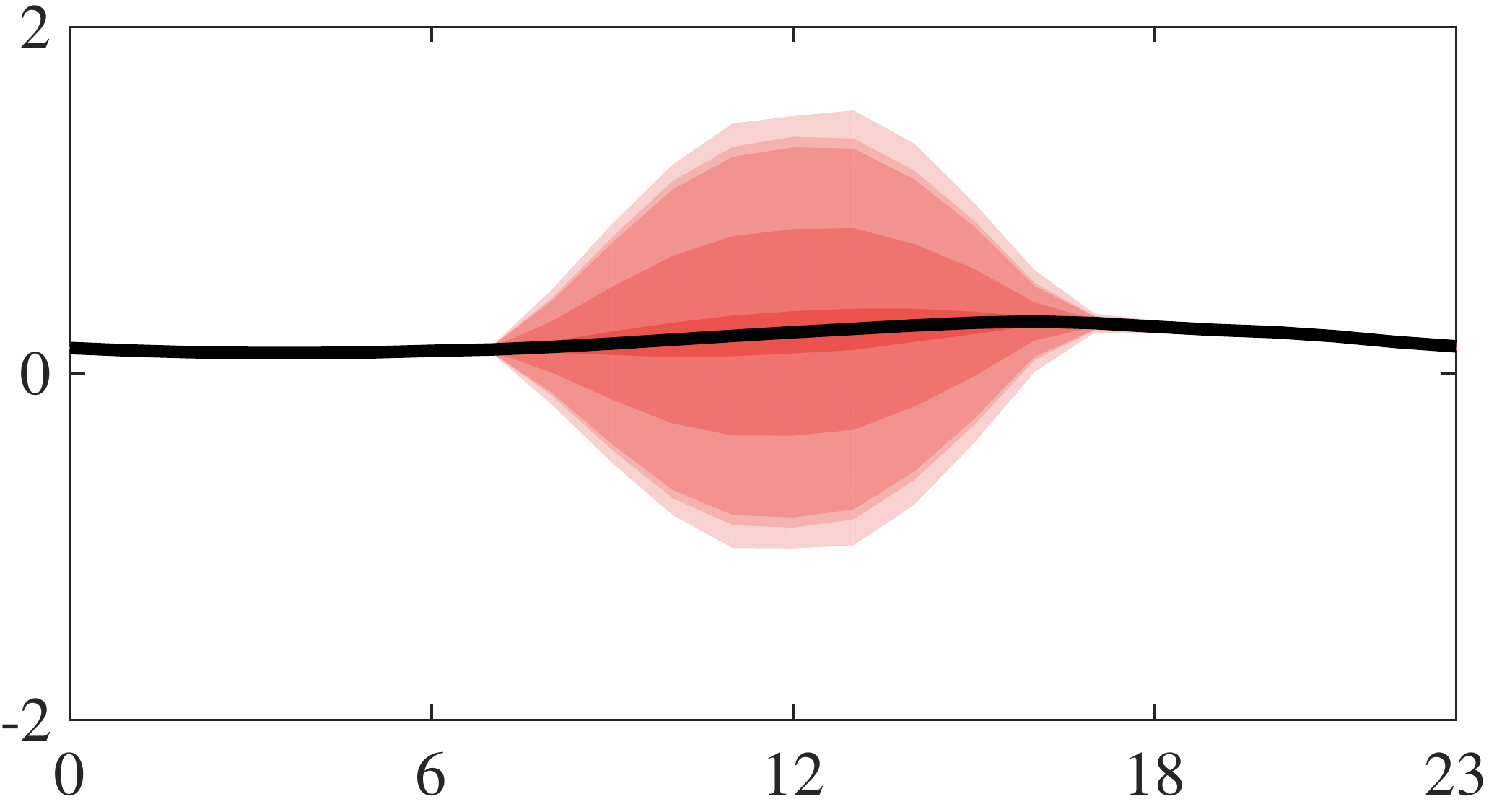}
\end{minipage}
\vspace{.005in}

\begin{minipage}[t]{0.49 \linewidth}
\centering
$p_4^S (t)$ (MW)

\includegraphics[width = \linewidth]{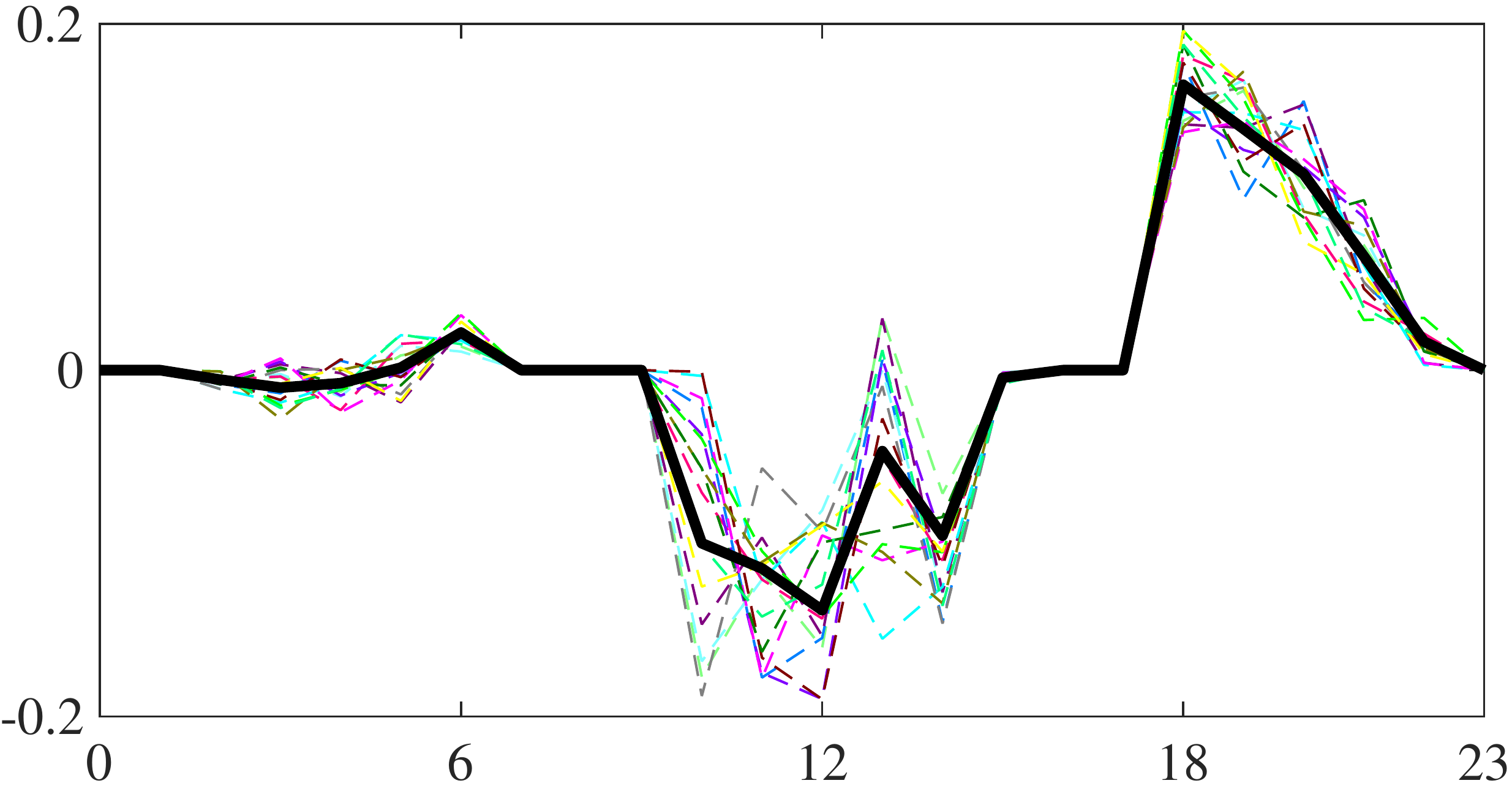}
\end{minipage}
\begin{minipage}[t]{0.49 \linewidth}
\centering
$p_4^S (t)$ (MW)

\includegraphics[width = \linewidth]{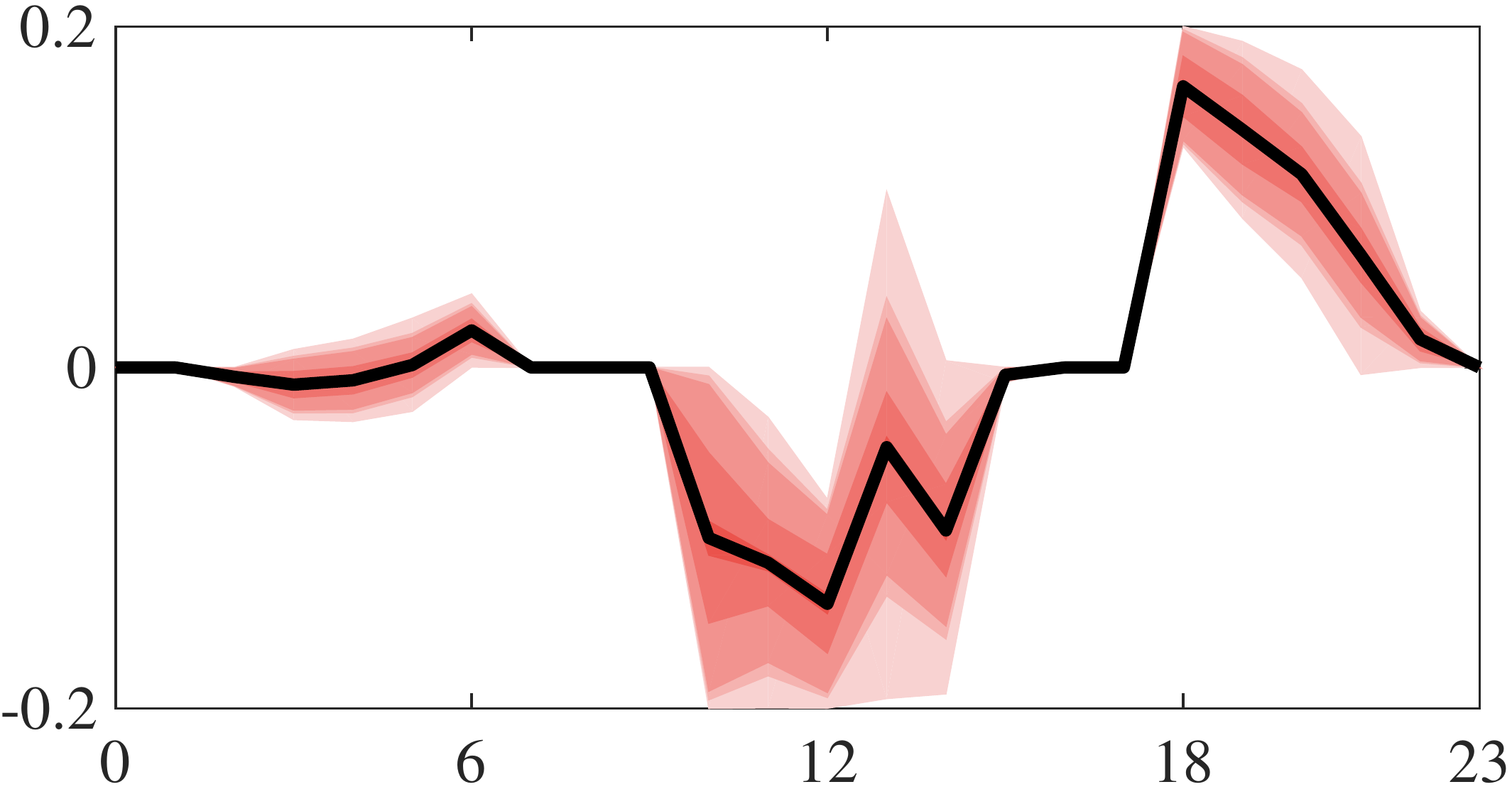}
\end{minipage}
\vspace{.005in}

\begin{minipage}[t]{0.49 \linewidth}
\centering
$x_4 (t)$ (MWh)

\includegraphics[width = \linewidth]{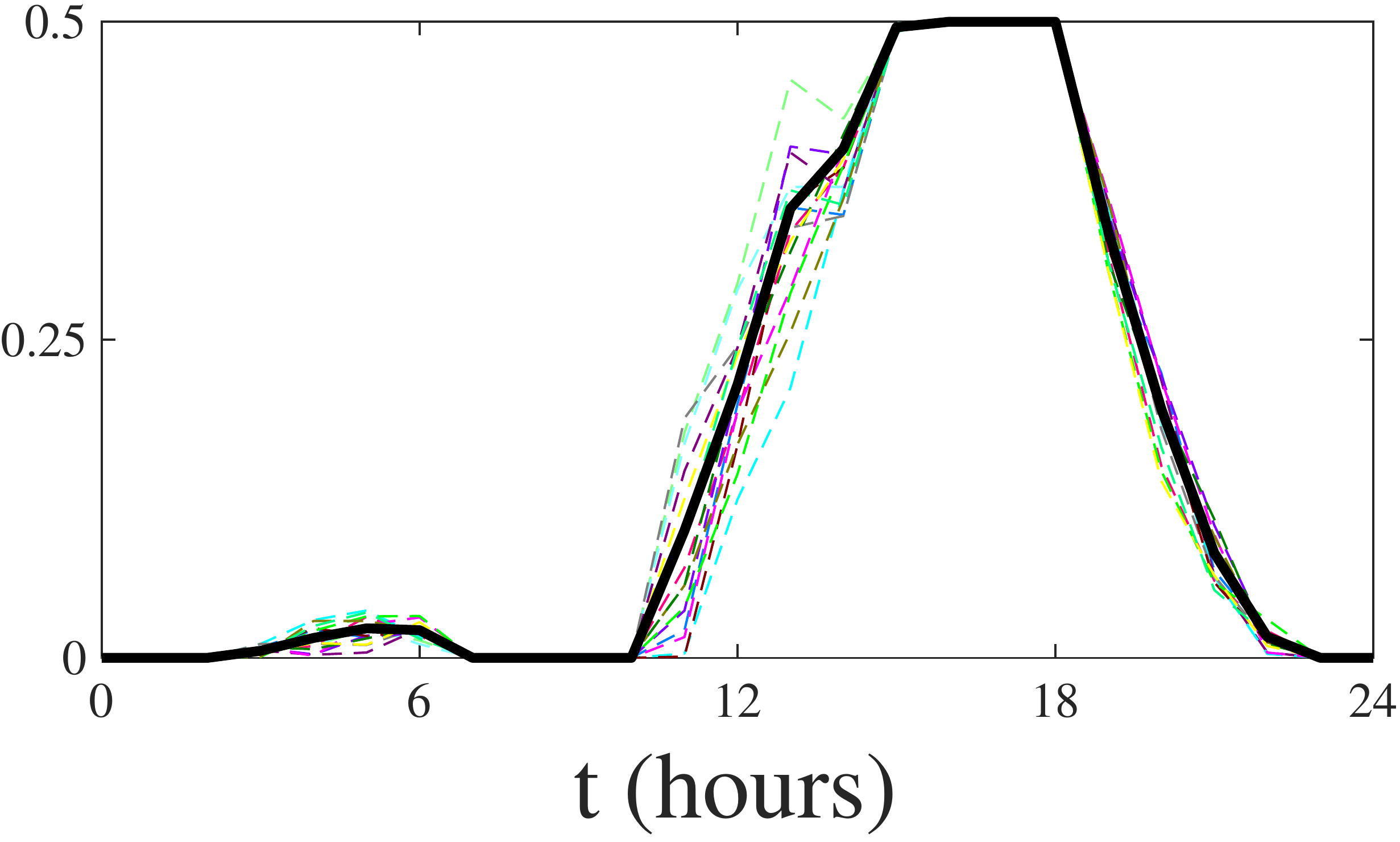}
\end{minipage}
\begin{minipage}[t]{0.49 \linewidth}
\centering
$x_4 (t)$ (MWh)

\includegraphics[width = \linewidth]{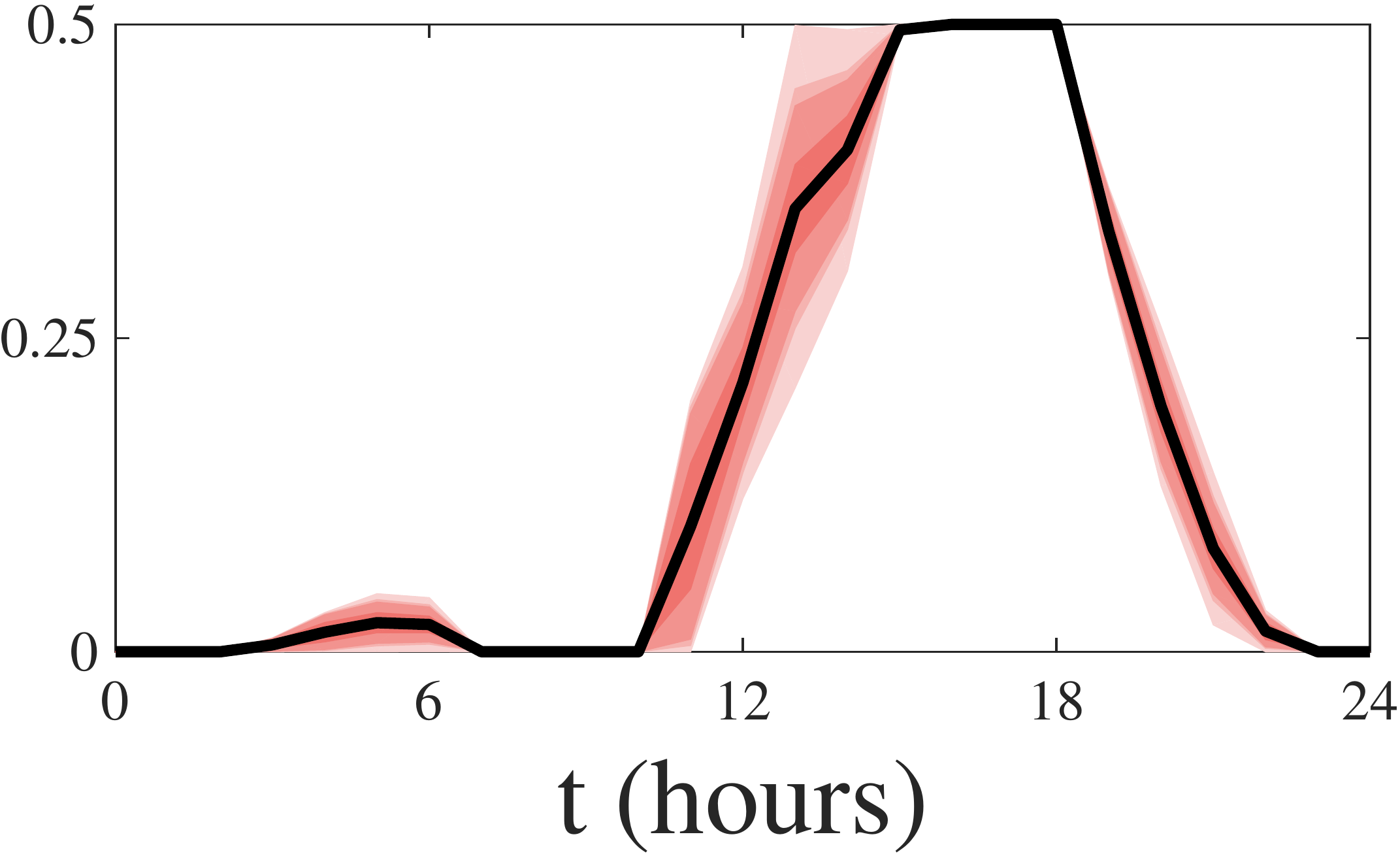}

\end{minipage}
\caption{System trajectories and their confidence intervals at bus 4.} \label{fig:control_traj_Bus4}
\end{subfigure}
\begin{subfigure}[t]{0.49784 \linewidth}
\centering
\begin{minipage}[t]{0.448 \linewidth}
\centering
$\xi_8^I (t)$ (MW)

\includegraphics[width = \linewidth]{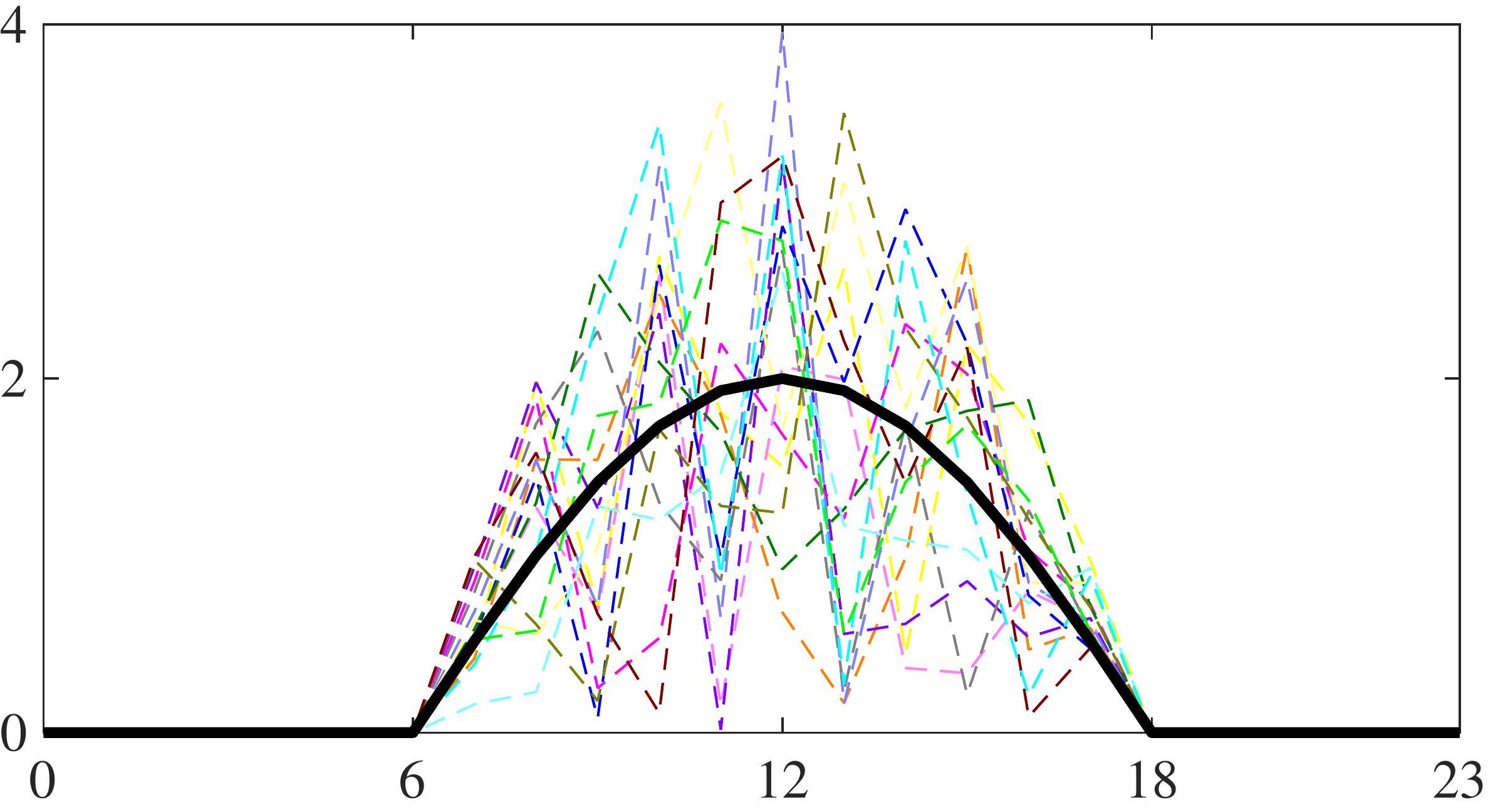}
\end{minipage}
\begin{minipage}[t]{0.532 \linewidth}
\centering
$\xi_8^I (t)$ (MW)

\includegraphics[width = \linewidth]{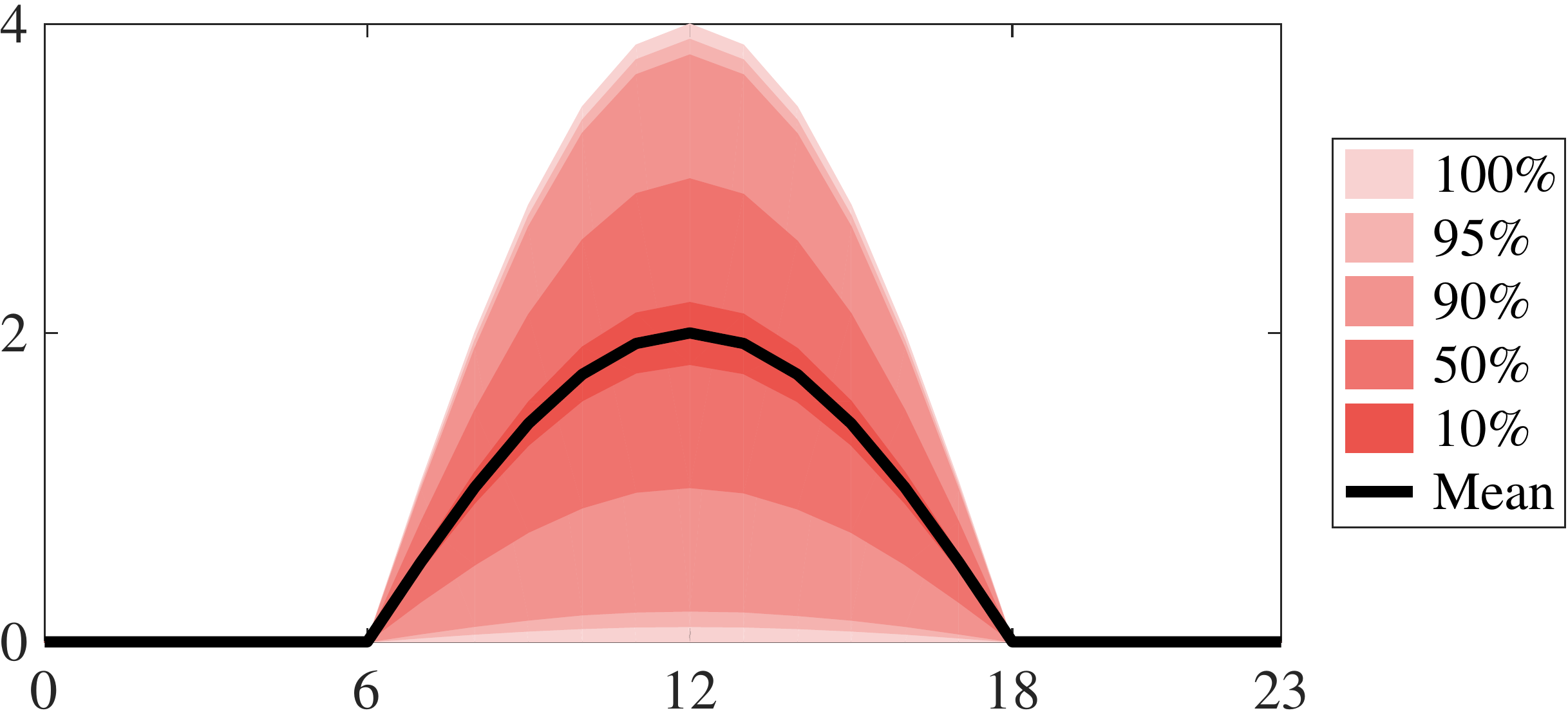}
\end{minipage}
\vspace{.005in}

\begin{minipage}[t]{0.448 \linewidth}
\centering
$q_8^I (t)$ (Mvar)

\includegraphics[width = \linewidth]{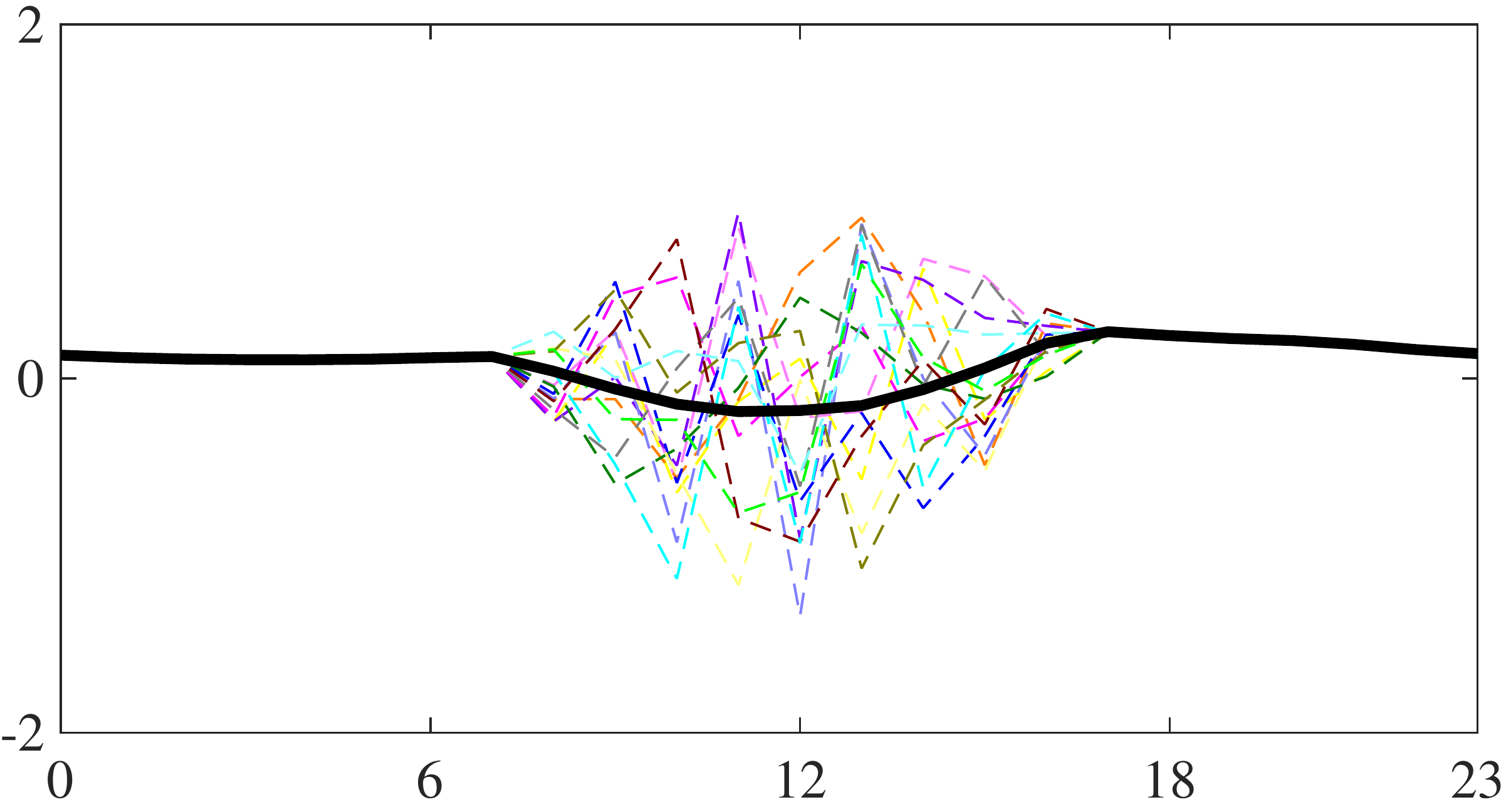}
\end{minipage}
\begin{minipage}[t]{0.532 \linewidth}
\centering
$q_8^I (t)$ (Mvar)

\includegraphics[width = \linewidth]{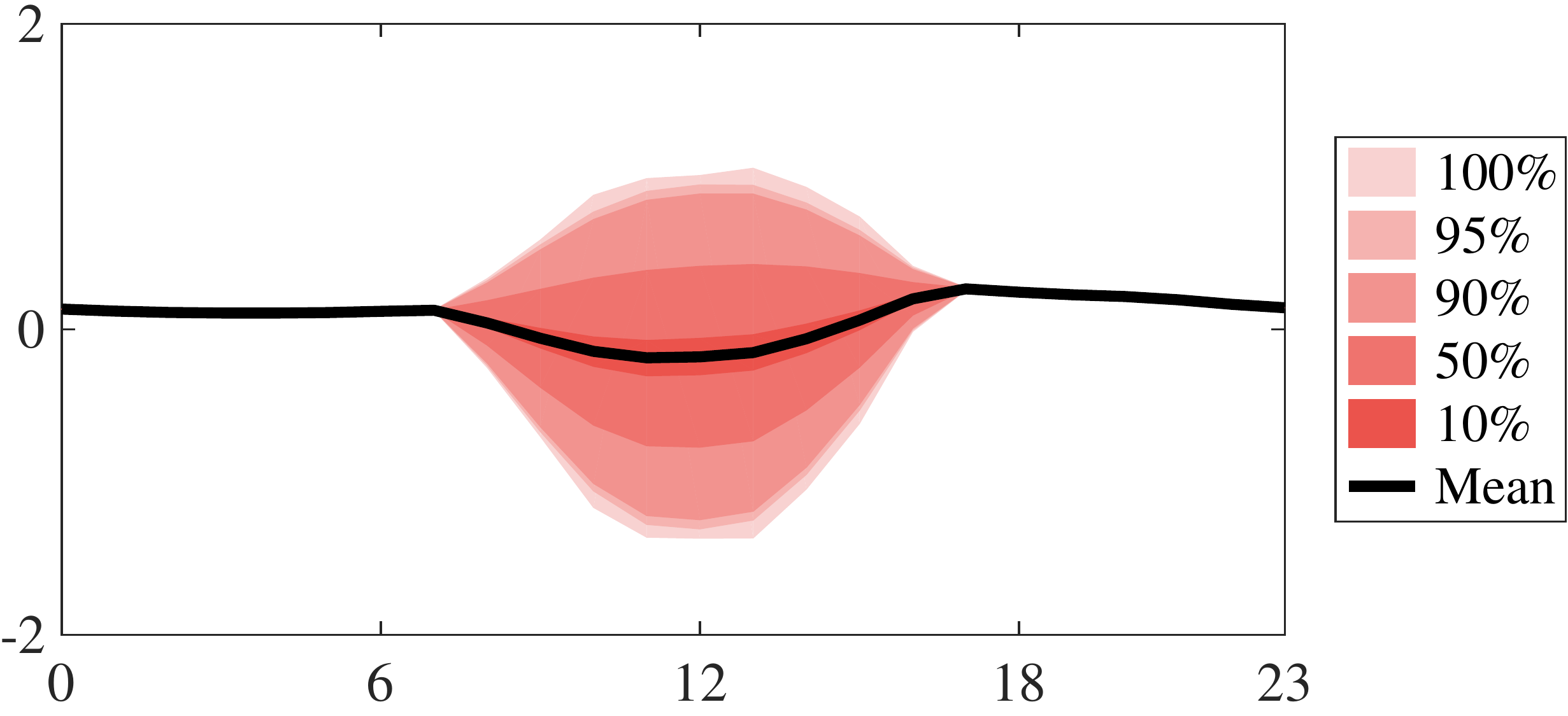}
\end{minipage}
\vspace{.005in}

\begin{minipage}[t]{0.448 \linewidth}
\centering
$p_8^S (t)$ (MW)

\includegraphics[width = \linewidth]{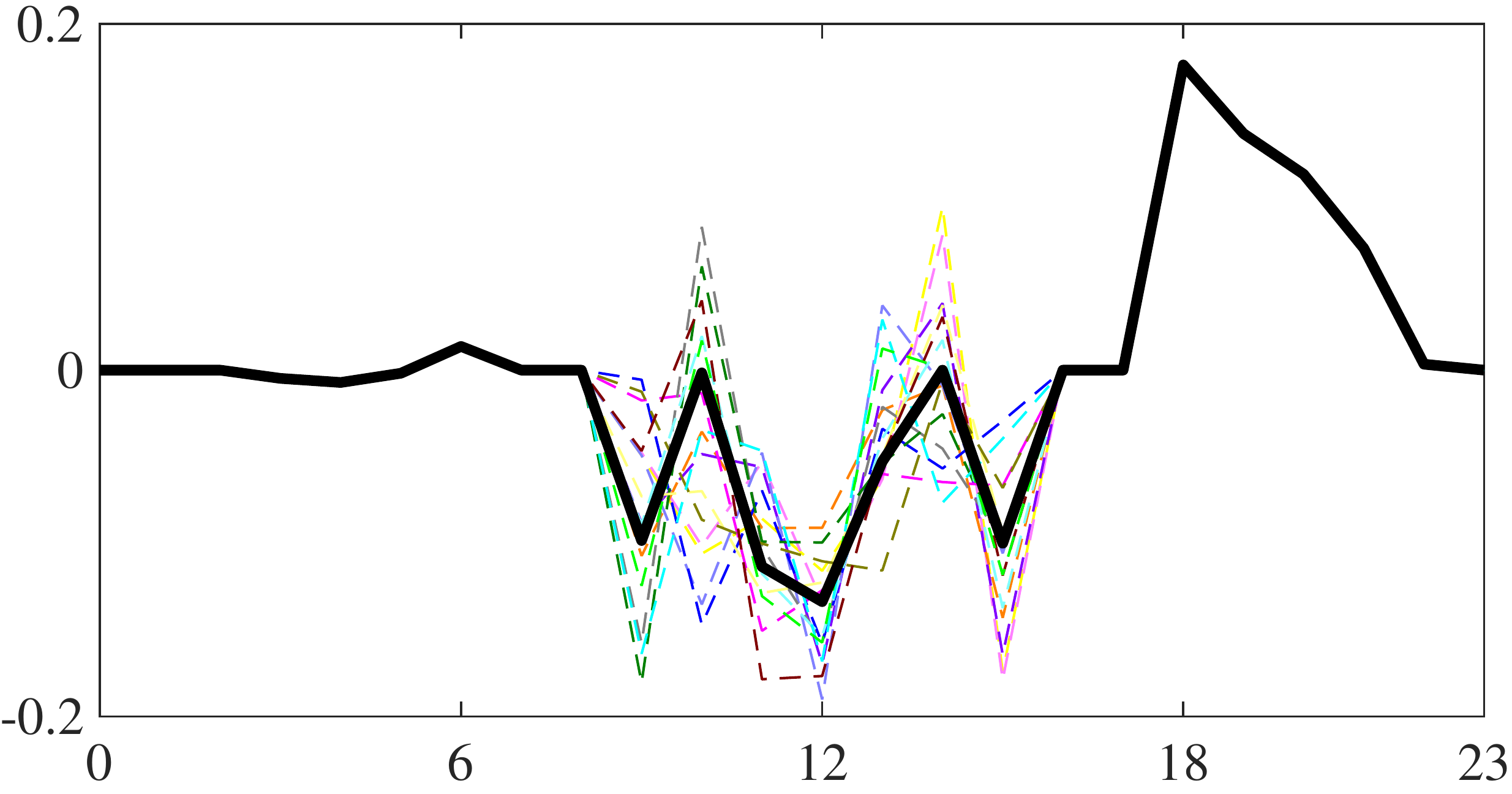}
\end{minipage}
\begin{minipage}[t]{0.532 \linewidth}
\centering
$p_8^S (t)$ (MW)

\includegraphics[width = \linewidth]{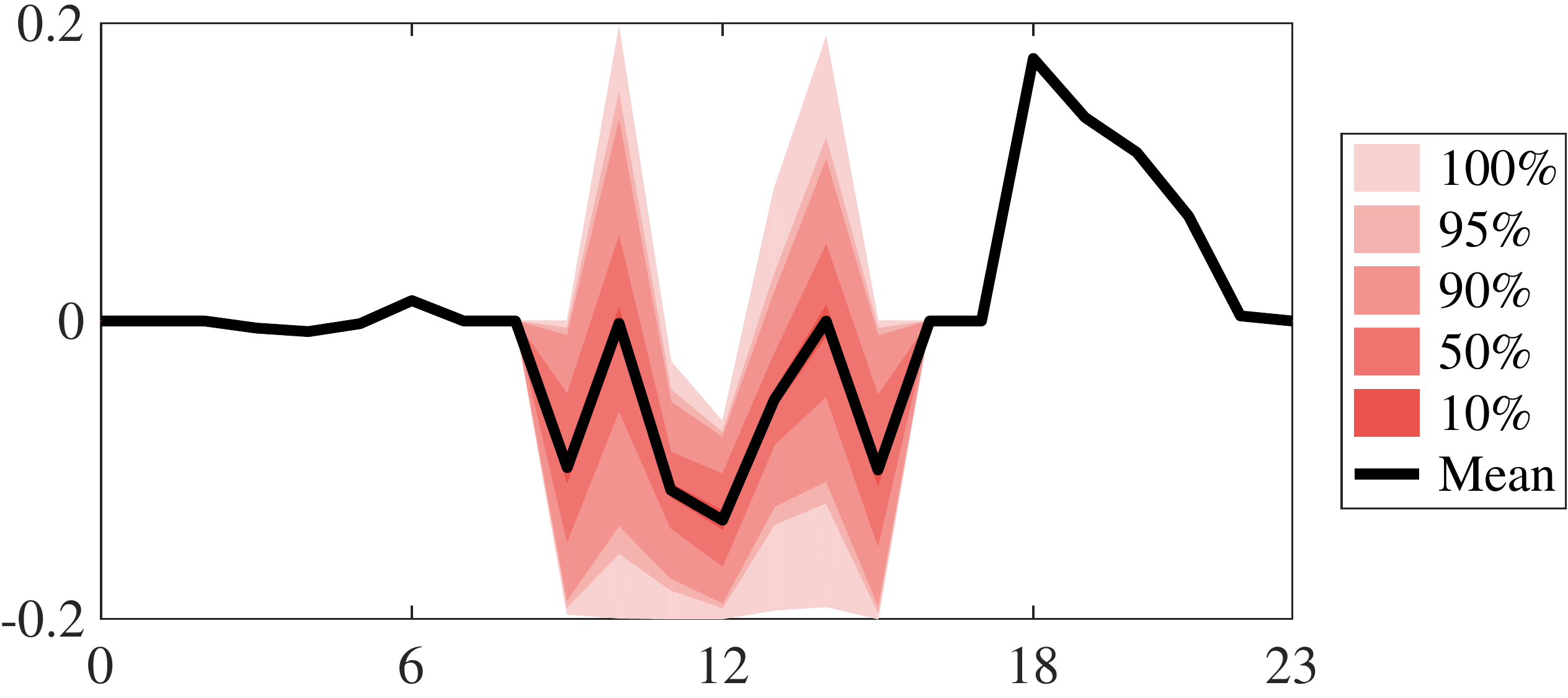}
\end{minipage}
\vspace{.005in}

\begin{minipage}[t]{0.448 \linewidth}
\centering
$x_8 (t)$ (MWh)

\includegraphics[width = \linewidth]{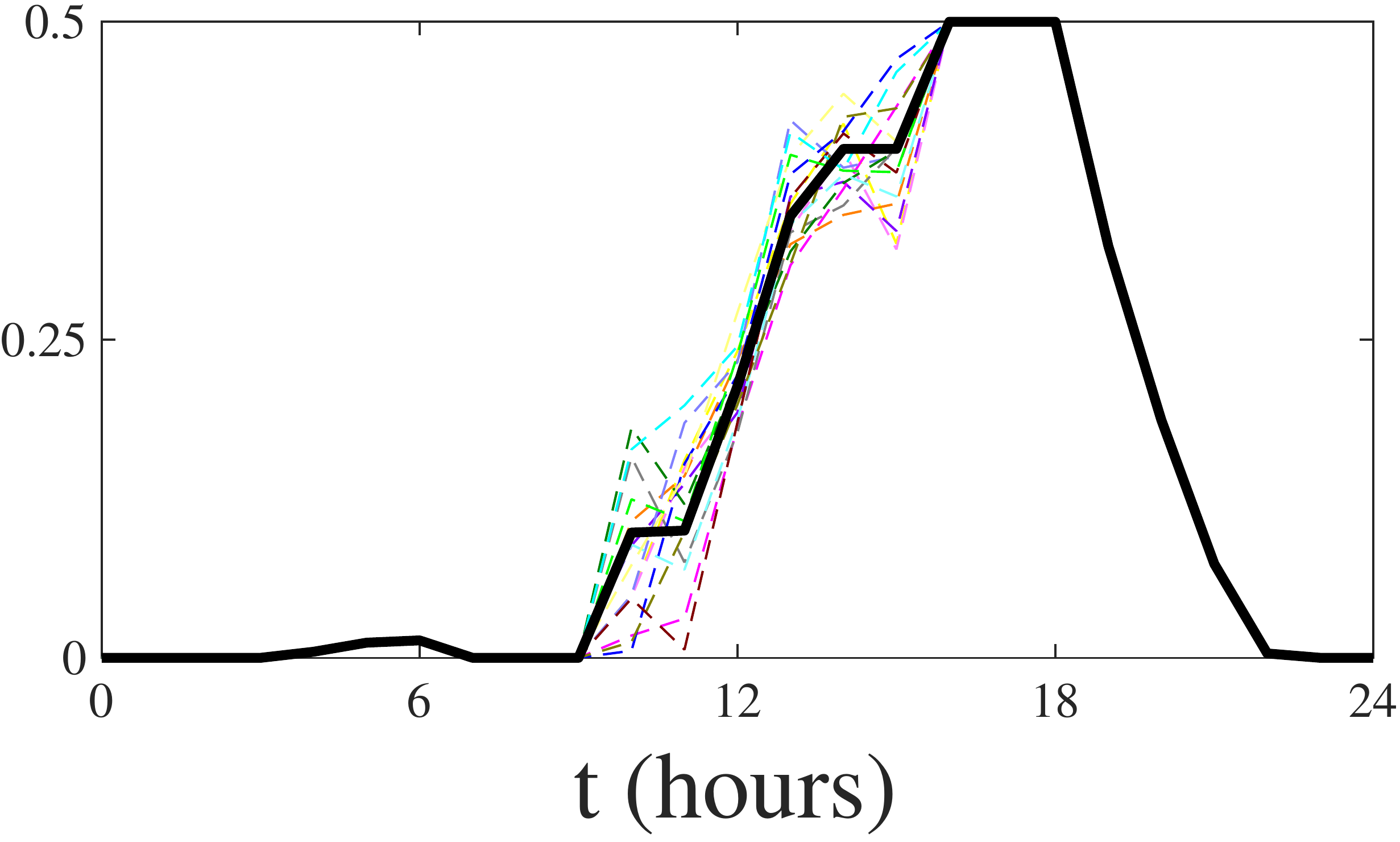}
\end{minipage}
\begin{minipage}[t]{0.532 \linewidth}
\centering
$x_8 (t)$ (MWh)

\includegraphics[width = \linewidth]{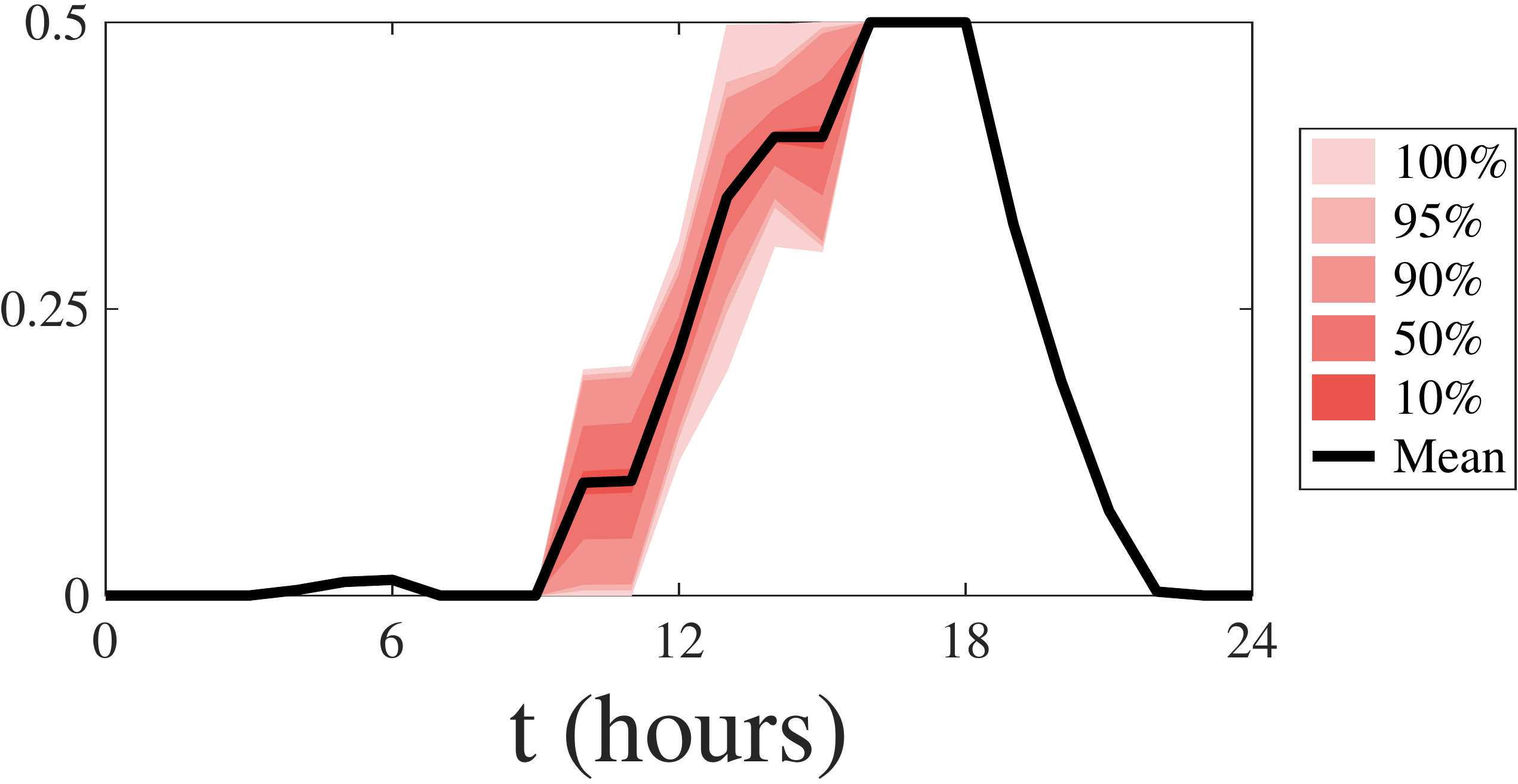}

\end{minipage}
\caption{System trajectories and their confidence intervals at bus 8.} \label{fig:control_traj_Bus8}
\end{subfigure}
\caption{The figures in the first and third columns plot independent realizations of disturbance, input, and state trajectories associated with bus 4 and 8, respectively. The  dashed colored lines represent the trajectory realizations, while the  solid black lines denote the mean trajectories.  The figures in the second and fourth columns depict the empirical confidence intervals associated with these trajectories. They were estimated using $3\times10^5$ independent realizations of the disturbance trajectories.}
\label{fig:control_traj}
\end{figure*}

\begin{figure*}[htbp]

\centering

\begin{subfigure}[t]{0.4557 \linewidth}
\centering
\begin{minipage}[t]{0.49 \linewidth}
\centering

$v_4(t)$ (pu)

\includegraphics[width = \linewidth]{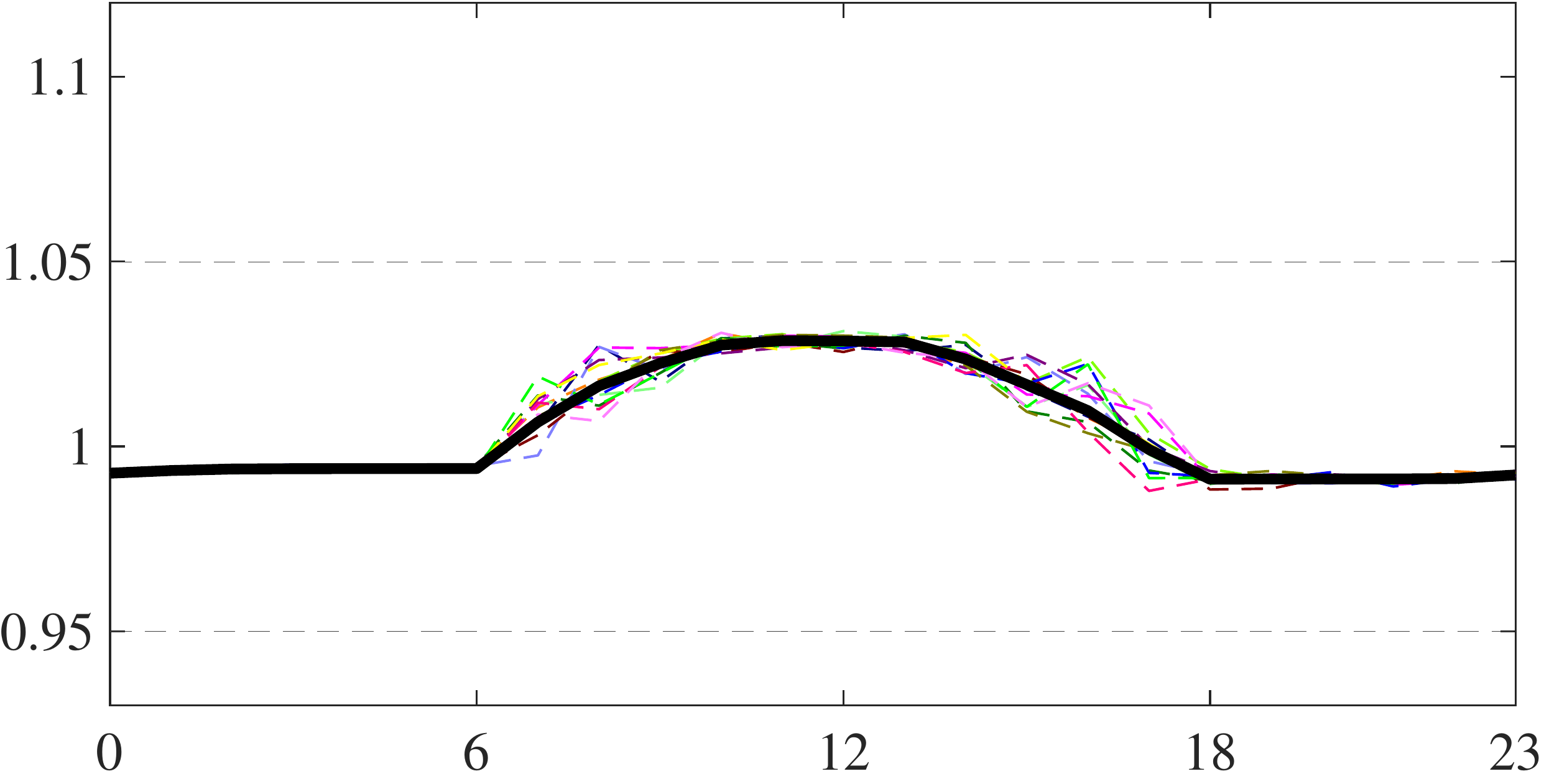}
\end{minipage}
\begin{minipage}[t]{0.49 \linewidth}
\centering

$v_4(t)$ (pu)

\includegraphics[width = \linewidth]{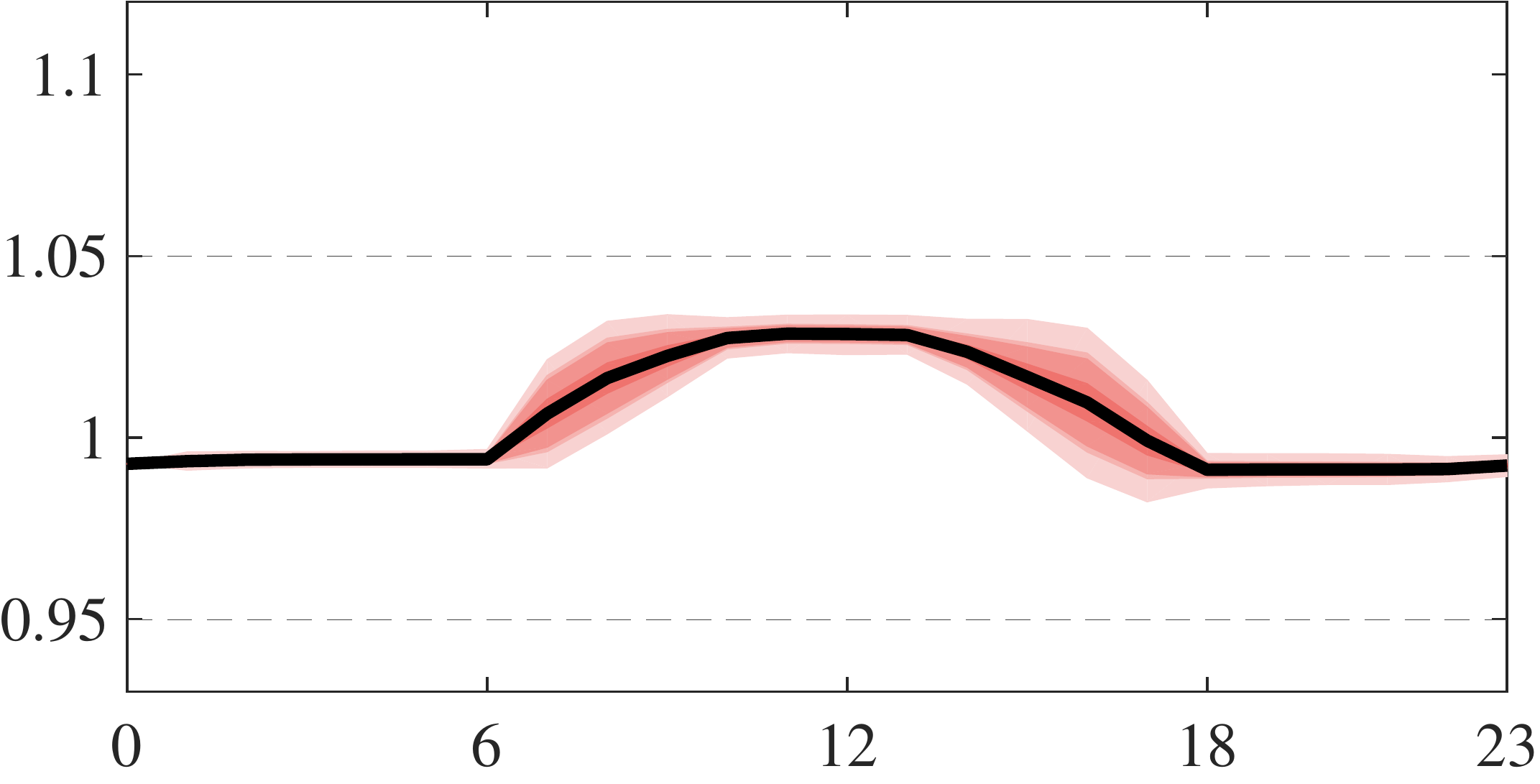}
\end{minipage}
\vspace{.005in}

\begin{minipage}[t]{0.49 \linewidth}
\centering

$v_8(t)$ (pu)

\includegraphics[width = \linewidth]{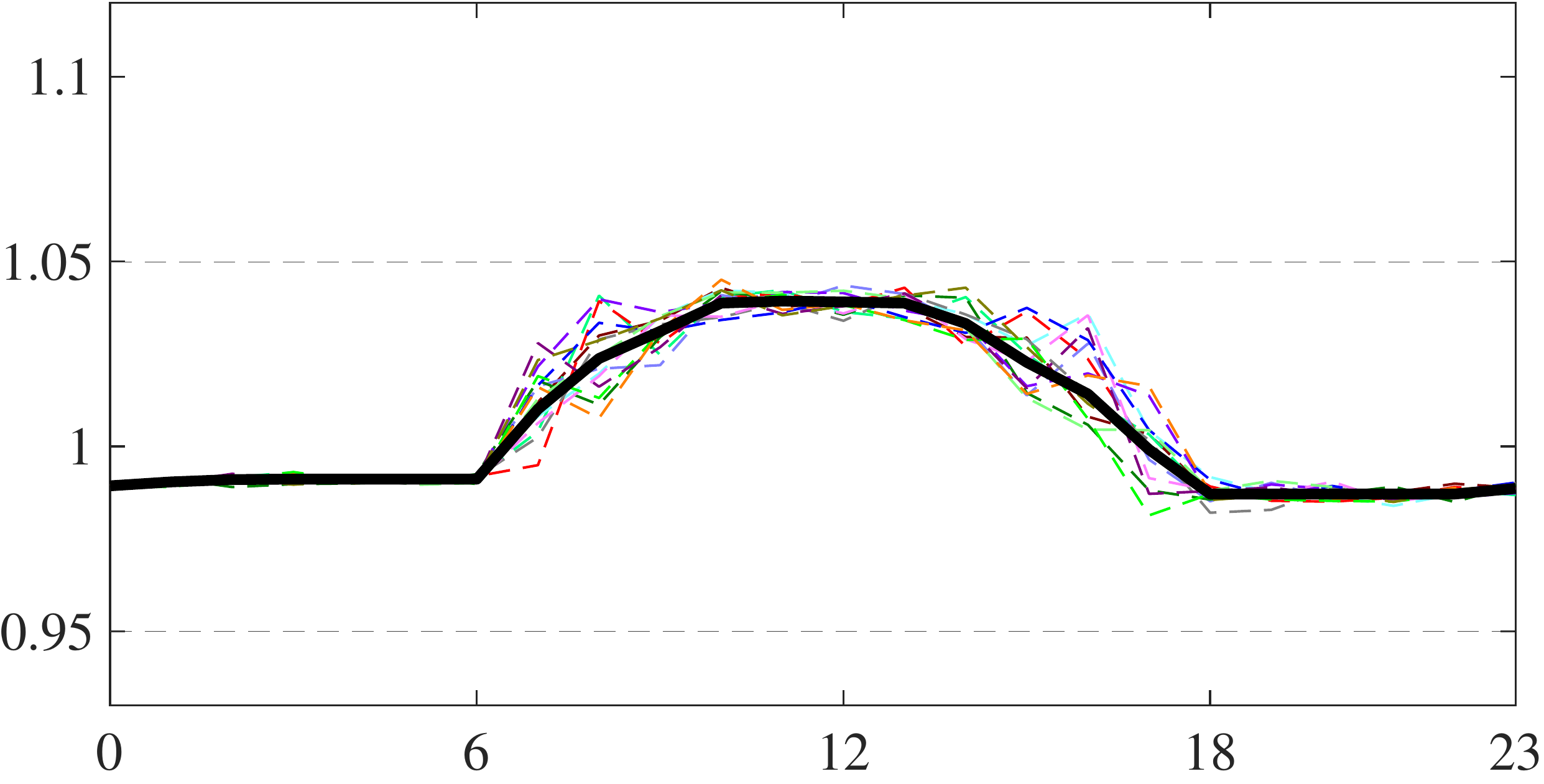}
\end{minipage}
\begin{minipage}[t]{0.49 \linewidth}
\centering
$v_8(t)$ (pu)

\includegraphics[width = \linewidth]{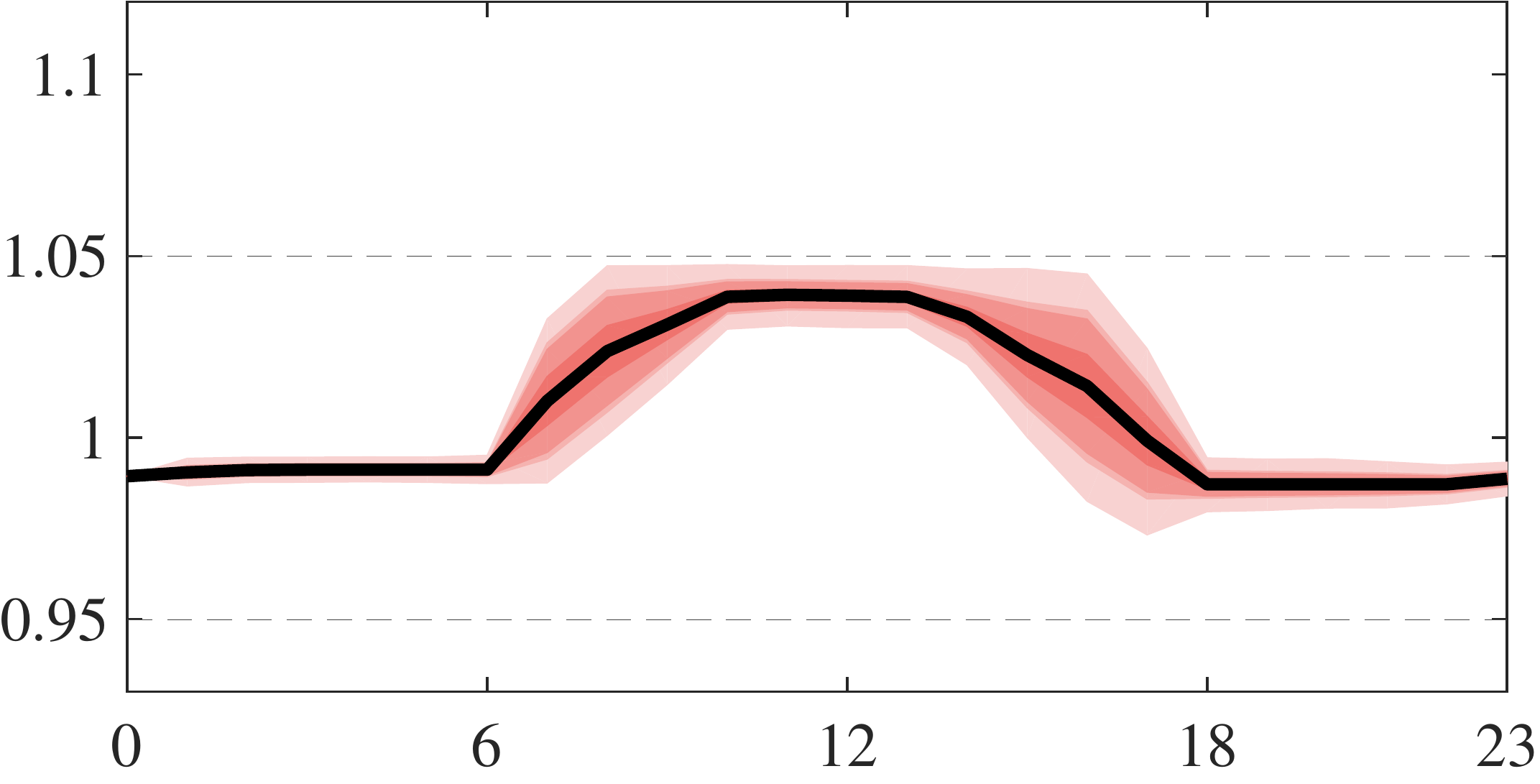}
\end{minipage}
\vspace{.005in}

\begin{minipage}[t]{0.49 \linewidth}
\centering
$v_{14}(t)$ (pu)

\includegraphics[width = \linewidth]{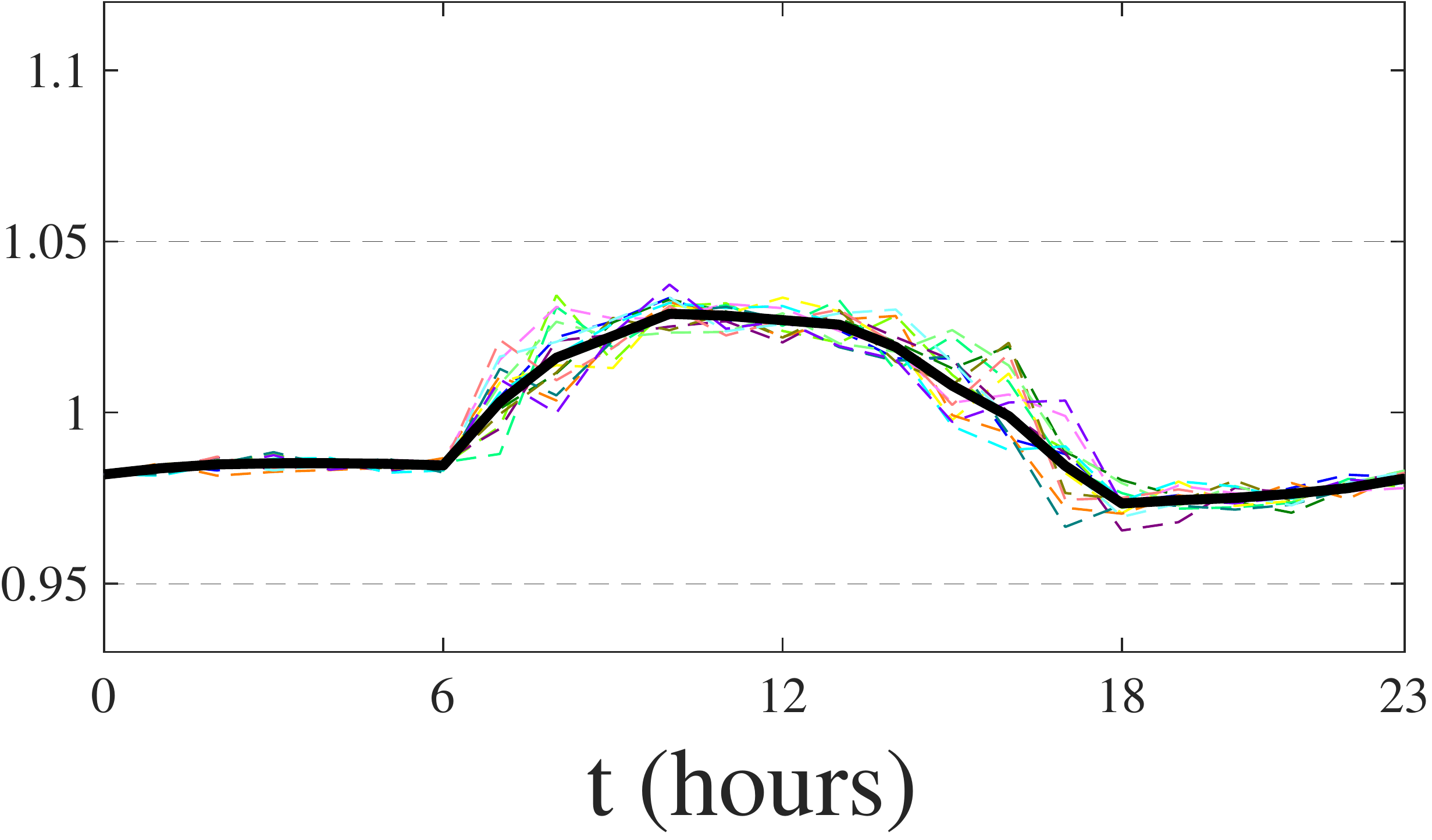}
\end{minipage}
\begin{minipage}[t]{0.49 \linewidth}
\centering
$v_{14}(t)$ (pu)

\includegraphics[width = \linewidth]{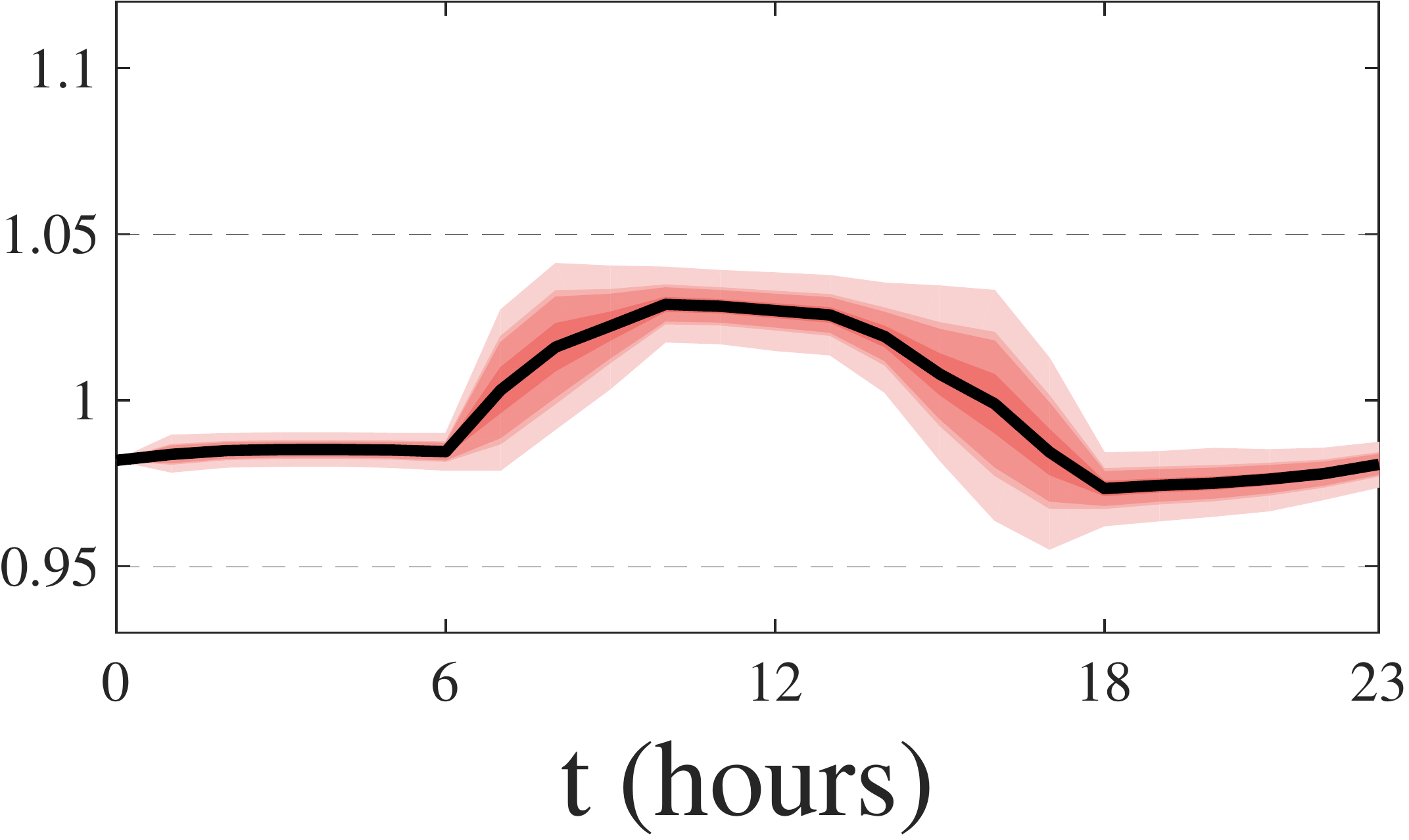}
\end{minipage}

\caption{Bus voltages in a \emph{controlled distribution system}.} \label{fig:volt_aff}
\end{subfigure}
\begin{subfigure}[t]{0.49784 \linewidth}
\centering
\begin{minipage}[t]{0.448 \linewidth}
\centering
$v_4(t)$ (pu)

\includegraphics[width = \linewidth]{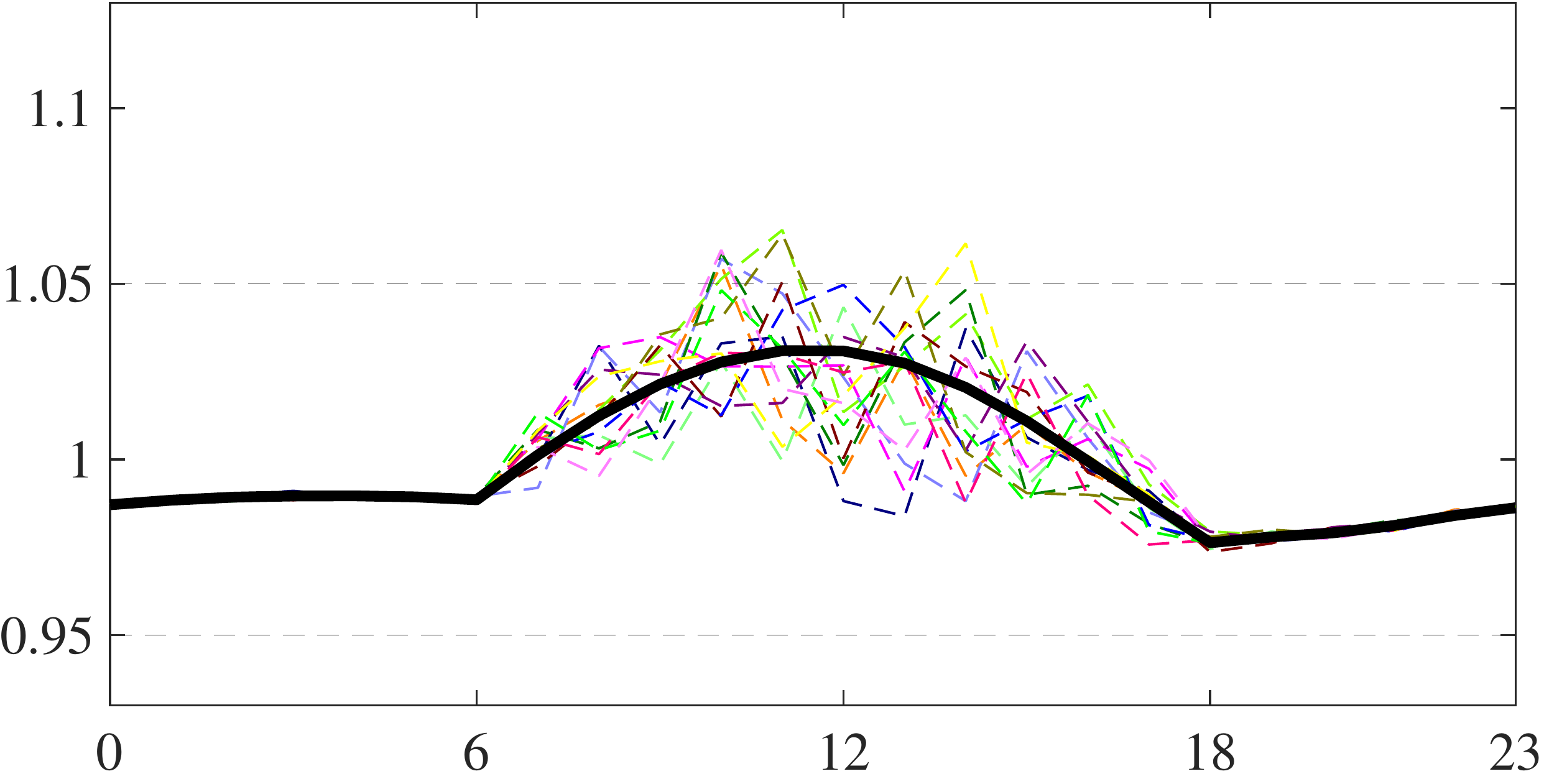}
\end{minipage}
\begin{minipage}[t]{0.532 \linewidth}
\centering
$v_4(t)$ (pu)

\includegraphics[width = \linewidth]{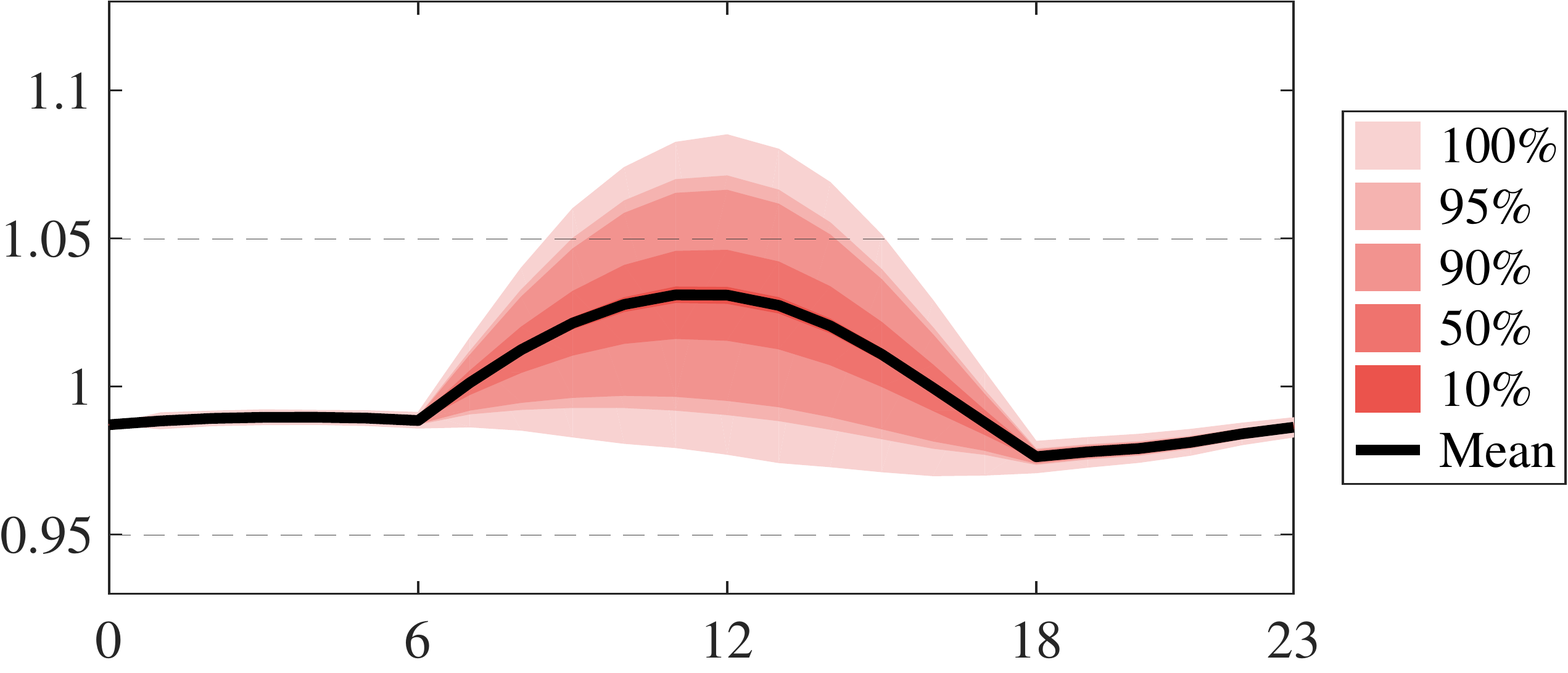}
\end{minipage}
\vspace{.005in}

\begin{minipage}[t]{0.448 \linewidth}
\centering
$v_8(t)$ (pu)

\includegraphics[width = \linewidth]{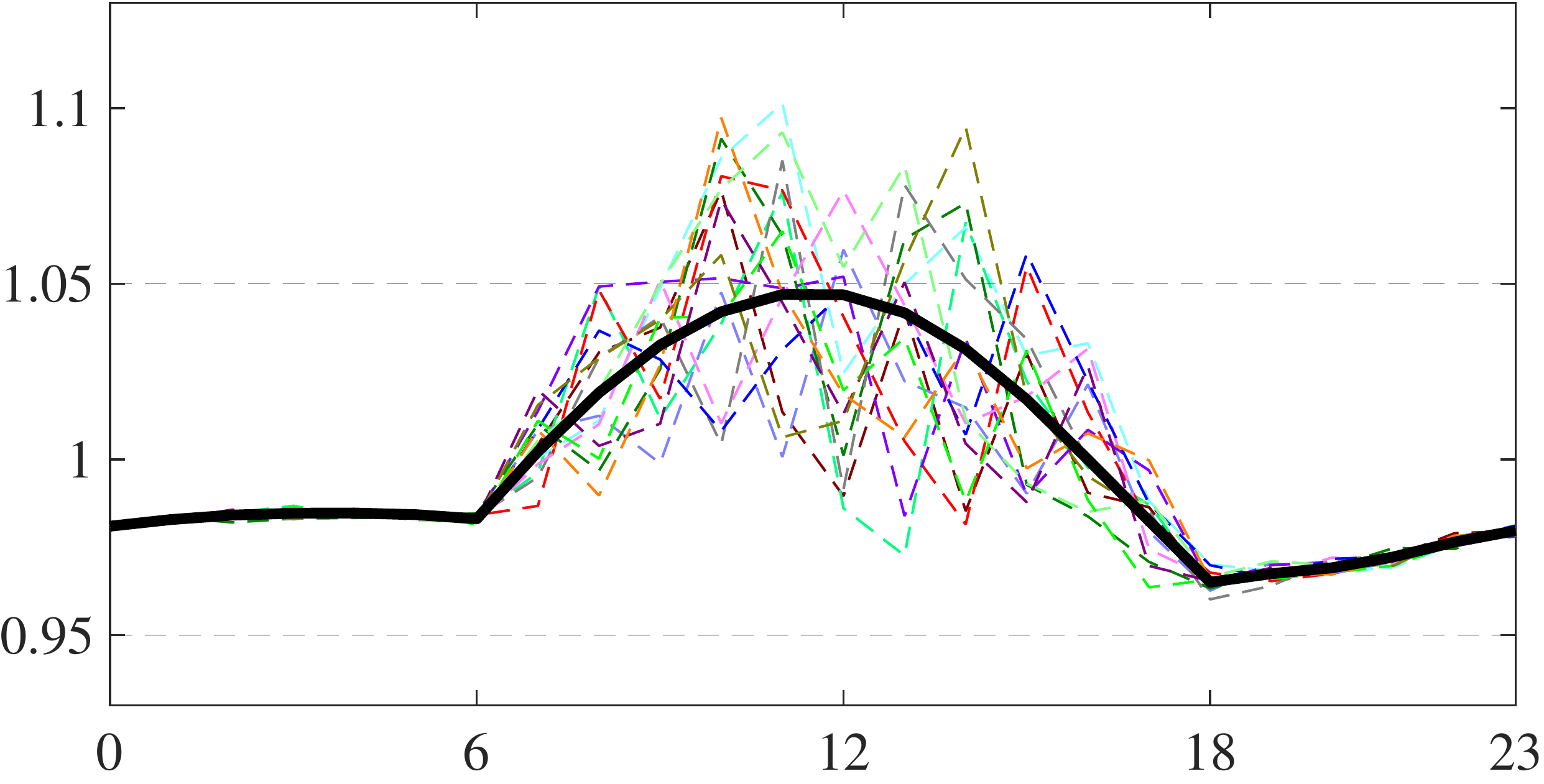}
\end{minipage}
\begin{minipage}[t]{0.532 \linewidth}
\centering
$v_8(t)$ (pu)

\includegraphics[width = \linewidth]{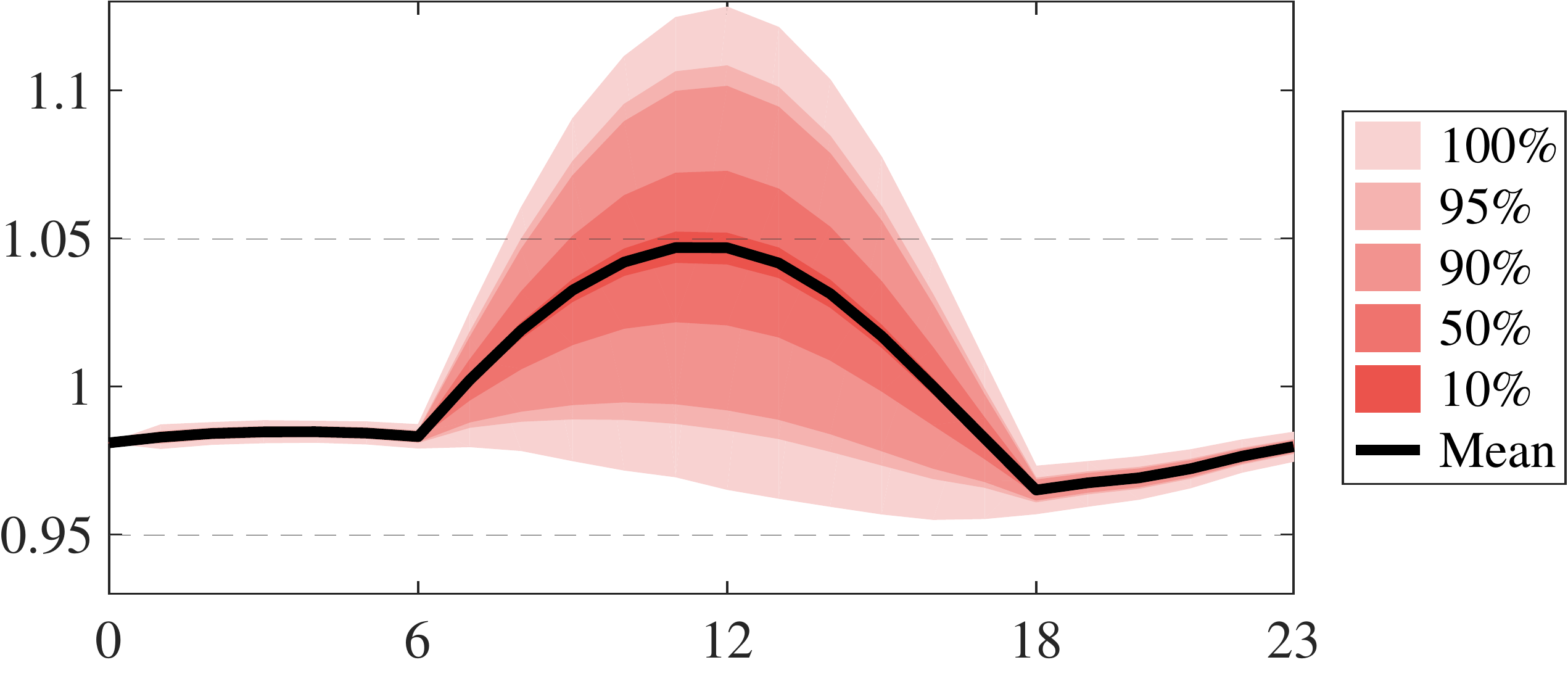}
\end{minipage}
\vspace{.005in}

\begin{minipage}[t]{0.448 \linewidth}
\centering
$v_{14}(t)$ (pu)

\includegraphics[width = \linewidth]{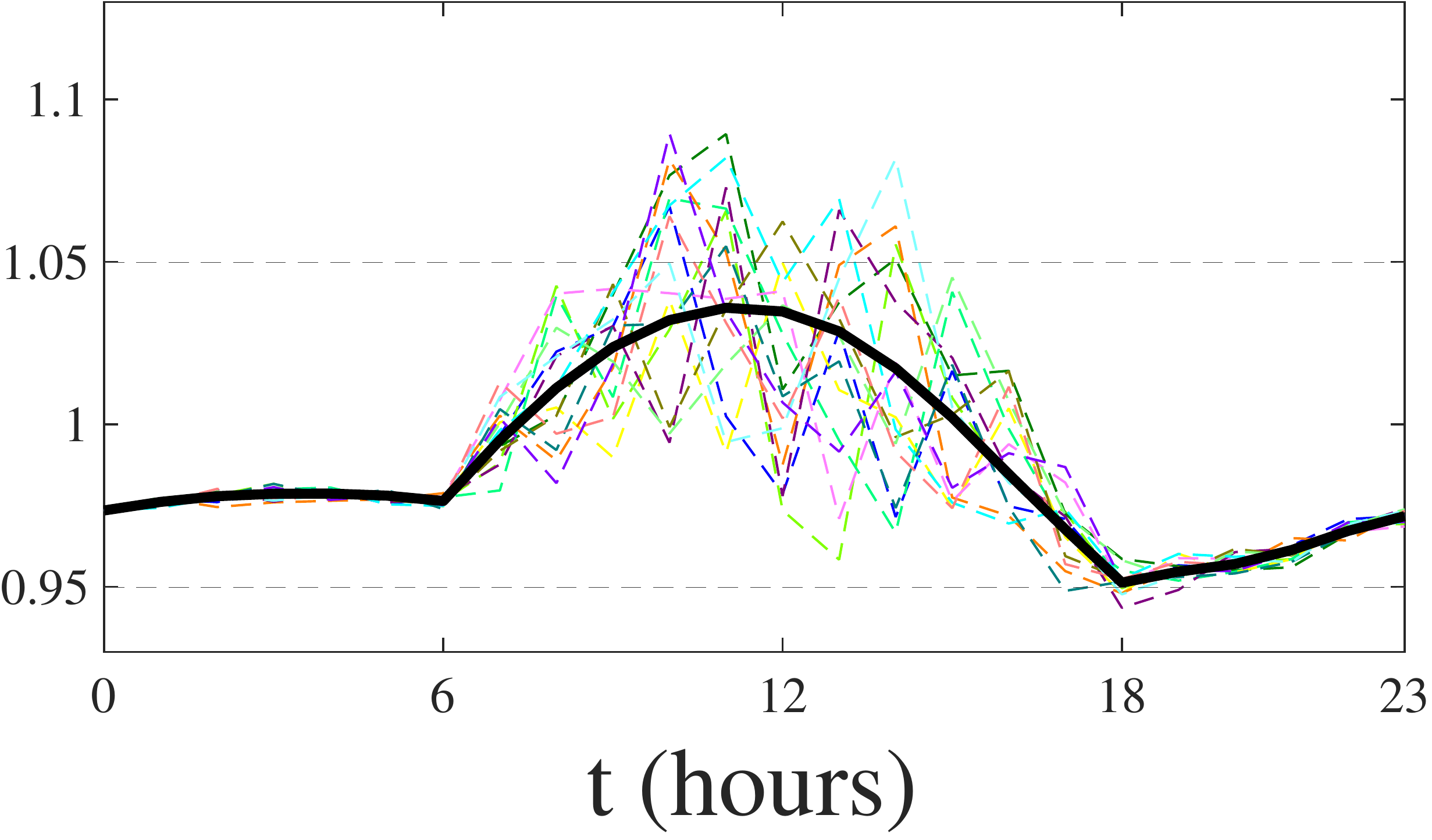}
\end{minipage}
\begin{minipage}[t]{0.532 \linewidth}
\centering
$v_{14}(t)$ (pu)

\includegraphics[width = \linewidth]{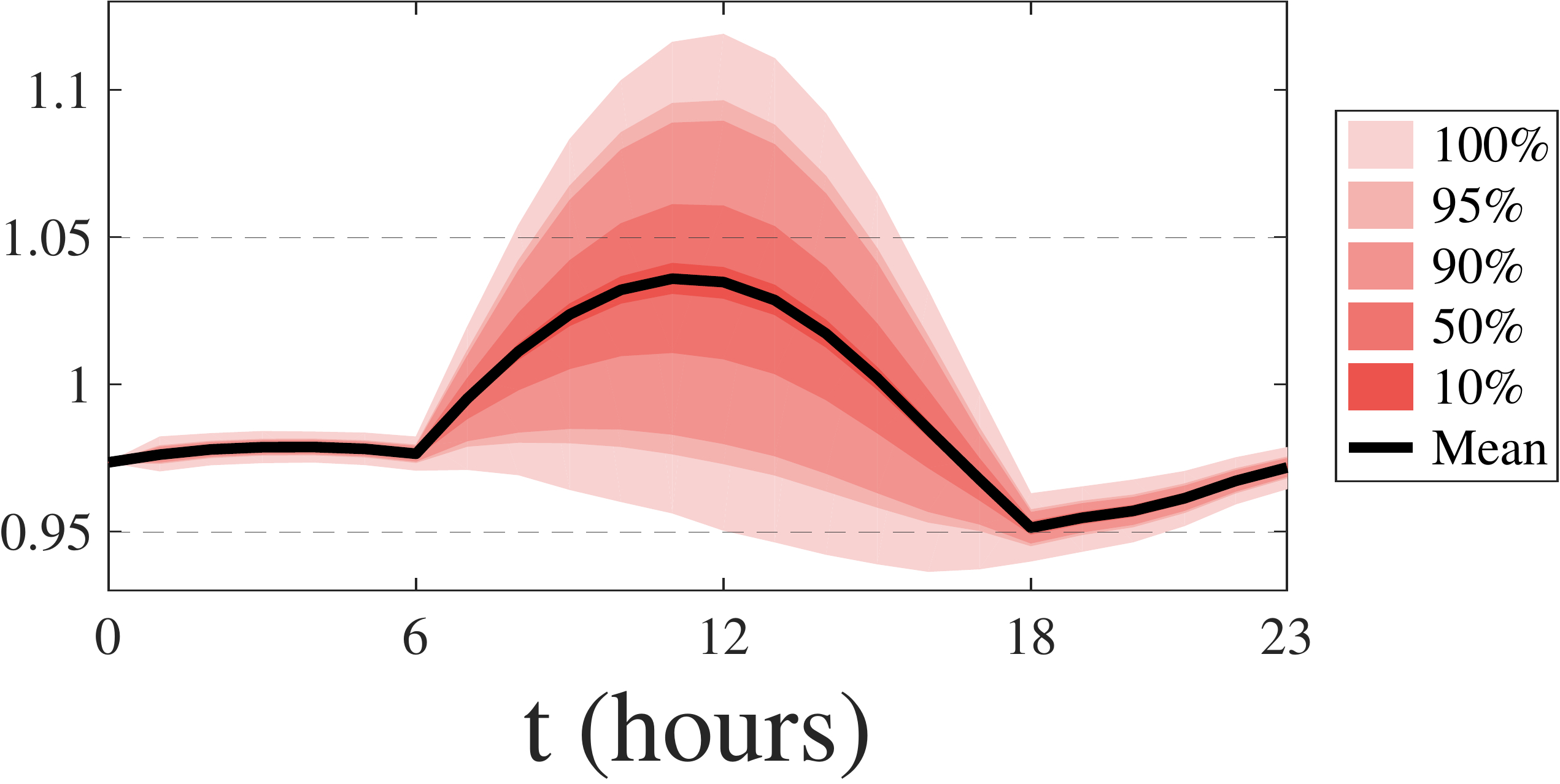}
\end{minipage}

\caption{Bus voltages in an \emph{uncontrolled distribution system}.} \label{fig:volt_uncontrolled}
\end{subfigure}
\caption{\rone{The above figures depict independent realizations of bus voltage magnitude trajectories and their empirical confidence intervals for (a) a \emph{controlled distribution system} operated under the decentralized affine controller, and (b) an \emph{uncontrolled distribution system}.  The empirical confidence intervals were estimated using $3\times10^5$ independent realizations of the disturbance trajectories. The dashed black lines indicate the range of allowable voltage magnitudes.}}
\label{fig:volt}
\end{figure*}

In Fig. \ref{fig:control_traj}, we illustrate the behavior of input and state trajectories generated by the decentralized affine controller computed according to Proposition \ref{prop:primal_affine_conic}. We consider the case of high PV penetration at a level of $\theta =4$ MW. 
In the first and third columns of Fig. \ref{fig:control_traj}, we plot several independent realizations of disturbance, input, and state trajectories associated with bus 4 and 8, respectively.
In the second and fourth columns, we plot the corresponding empirical confidence intervals.\footnote{The empirical confidence intervals were estimated using $3\times10^5$ independent realizations of the disturbance, input, and state trajectories.}
First, notice that both the sequence of reactive power injections from PV inverters and the sequence of active power extractions from storage exhibit large fluctuations during  daytime hours.  These fluctuations are due in large part to the underlying variability in the active power supplied by the PV resources. In particular, a large excess of active power supply from PV can manifest in overvoltage in the distribution network.  In order to ensure that voltage magnitude constraints
are not violated, the proposed control policy induces reactive power injections from PV inverters that are negatively correlated with their own active power supply. Clearly, in the absence of such a feedback control mechanism, certain realizations of the disturbance trajectory would have resulted in the violation of the voltage magnitude constraints at certain buses in the distribution system.

\rone{
In Fig. \ref{fig:volt}, we illustrate the effectiveness of the proposed decentralized affine controller in maintaining bus voltage magnitudes within their allowable range. In particular, we compare the behavior of  bus voltage magnitudes that occur in the distribution  system \emph{with} and \emph{without control}. In the first and second columns of Fig. \ref{fig:volt}, we illustrate the behavior of voltage magnitude trajectories that materialize in the \emph{controlled distribution system} operated under the proposed decentralized affine controller. In the third  and fourth columns of Fig. \ref{fig:volt}, we illustrate the behavior of voltage magnitude trajectories that materialize in the \emph{uncontrolled distribution system}, i.e., under the control policy $\gamma = 0$. Notice that, in the absence of  control, the distribution system may realize bus voltage magnitudes that substantially deviate from their allowable range. In particular, the distribution system appears to suffer from overvoltage when there is an overabundance of active power supply from PV, and undervoltage during hours of peak demand. However, when operated under the decentralized affine controller, the distribution system is guaranteed to satisfy the bus voltage magnitude constraints for any possible realization of the disturbance trajectory.
}

\section{Conclusion} \label{sec:Conclusions}

There are several interesting directions for future work. 
For example, one potential drawback of the approach considered in this paper is the explicit reliance of the control policy on the entire disturbance history.
Such dependency may result in the computational intractability 
of calculating control policies for problems with a long horizons $T$.
Accordingly, it will be of interest to extend the techniques developed in this paper to accommodate fixed-memory constraints on the control policy.

It is also worth noting that the class of controllers considered in this paper are \emph{fully decentralized}, in that explicit communication between subsystems is not permitted. It would be of theoretical and practical interest to investigate the extent to which the introduction of
additional communication links between subsystems might improve system performance. In particular, it would be of interest to explore the problem of designing a communication  topology between subsystems, in order to minimize the optimal control cost, subject to a constraint on the maximum number of allowable communication links.
While such problems are inherently combinatorial in nature, it is conceivable that regularization techniques, similar to those proposed in \cite{lin2013design, matni2014regularization}, might yield good approximations.

Finally, all of our results rely on the assumption that the distribution system is three-phase balanced. Such an assumption will not always hold in practice.
It would be of interest to extend the techniques in this paper to accommodate the possibility of imbalance in three-phase distribution systems.

\section{Acknowledgments}
The authors would like to thank the Editor and four anonymous reviewers for their helpful comments and feedback.
This work was supported in part by the US National Science Foundation (Contracts: ECCS-1351621, CNS-1239178, IIP-1632124), the US Department of Energy under the CERTS initiative, and the Simons Institute for the Theory of Computing. 

\begin{appendices}

\section{Matrix Definitions}  \label{app:matrix}

The block matrices $(\Ambb, \Bmbb)$ used in Eq. \eqref{eq:x_full} are given by
\begin{align*}
\Ambb = \bone_{(T+1) \times 1} \otimes I_n , \qquad  \Bmbb  = \bmat{ 0 \\  B &0 \\  B &  B & 0  \\  \vdots & \vdots & \ddots & \ddots \\ \vdots & \vdots & &\ddots &0  \\  B &  B  & \cdots & \cdots & B },
\end{align*}
where $\bone_{(T+1) \times 1} $ is a vector of all ones in $\RR^{T+1}$. The vector $c$ used in Eq. \eqref{eq:terminal_storage} is given by
\begin{align*}
c = \bmat{0_{nT \times 1} \\ \bone_{n \times 1}}.
\end{align*}

To define the matrices $L_u, \ L_\xi$, and $\Sigma$ used in Eq. \eqref{eq:loss}, one can first construct matrices $L_u^0, \ L_\xi^0 \in \RR^{2n \times 3n}$, and a positive definite diagonal matrix $\Sigma^0 \in \RR^{2n \times 2n}$, such that
\begin{align*}
& \sum_{(i,j) \in \Ecal} \frac{r_{ij} \left( p_{ij} (t)^2 + q_{ij} (t) ^2 \right)}{v_0^2} \\
& \hspace{0.9in} =\Big( L_u^0 u (t) + L_\xi^0 \xi (t) \Big)' \Sigma^0 \Big(L_u^0 u (t) + L_\xi^0 \xi (t) \Big)
\end{align*}
for $t = 0, \dots, T-1$. Then one can define matrices $\Sigma, L_u,$ and  $L_\xi$ according to 
\begin{align*}
\Sigma  = I_T \otimes \Sigma^0, \quad L_u  = I_T \otimes L_u^0,  \quad L_\xi = \bmat{0 & I_{T} \otimes L_\xi^0}.
\end{align*}

\section{Fast Time-Scale Controller Implementation} \label{app:continuous_time}

We now describe a method to enable the implementation of the decentralized control policy computed according to Proposition \ref{prop:primal_affine_conic} over a more finely grained (i.e., fast) time-scale.  The method we propose is simple. First, we compute  an affine control policy for the original (i.e., slow) time-scale according to Proposition \ref{prop:primal_affine_conic}. Via a suitable rescaling of the resulting feedback control gains,  we construct an affine control policy that can be implemented over a more finely grained time-scale. An attractive feature of the proposed implementation is that,  under a mild assumption on the quasi-stationarity of the support of the underlying disturbance trajectory, the  affine control policy we construct is guaranteed to yield 
state and input trajectories that are feasible on this more finely grained time-scale. In what follows, we provide a precise specification of this fast time-scale controller, and discuss its theoretical guarantees.

\subsection{Fast Time-Scale Processes}
We begin with a description of  the state, input, and disturbance processes on the more finely grained time-scale by
dividing each original time period $t$ into $K$ shorter time periods. It will be convenient to  index the original time periods by $t$, and  the more finely grained time periods by $(k, t)$, for $k = 0, \dots, K-1$, and $t = 0, \dots, T$.
In particular, each time period $(k,t)$ is defined over a time interval of length $\Delta / K$, where recall that each original time period $t$ is of length $\Delta$. For the remainder of this section, we will refer to the original and the more finely grained time-scales as the \emph{slow} and \emph{fast time-scales}, respectively.

\begin{figure}[http]
\centering

\tikzstyle{Line_thick} = [line width = 2pt]
\tikzstyle{Line_semi_thick} = [line width = 1.5pt]
\tikzstyle{Line} = [line width = 1pt]
\tikzstyle{Arrow} = [-angle 60, line width = 1pt]
\tikzstyle{Double_arrow} = [angle 60 - angle 60, line width = 0.75pt]
\tikzstyle{Line_dashed} = [line width = 1pt, densely dashed]
\tikzstyle{blank} = [circle, inner sep = 0pt, minimum size = 0mm]
\tikzstyle{Line_r_thick} = [draw = red, line width = 2.5pt]
\tikzstyle{Line_b_thick} = [draw = blue, line width = 2.5pt]
\tikzstyle{Line_r_thick_dashed} = [draw = red, line width = 2.5pt, densely dashed]
\tikzstyle{Line_b_thick_dashed} = [draw = blue, line width = 2.5pt, densely dashed]
\tikzstyle{textbox} = [below, text centered]
\tikzstyle{textbox_above} = [above, text centered]

\begin{tikzpicture}

\newcommand*{\axup}{3}
\newcommand*{\axdown}{0}
\newcommand*{\axleft}{-4.5}
\newcommand*{\length}{7.2}
\newcommand*{\tfast}{0.8}
\newcommand*{\tslow}{2.4}
\newcommand*{\vscale}{1}
\newcommand*{\lslow}{0.5}
\newcommand*{\lfast}{0.2}
\newcommand*{\edge}{0.7}

\coordinate (axup_left) at (\axleft, \axup);
\coordinate (axup_right) at ({\axleft + \length + \edge}, \axup);

\coordinate (axdown_left) at (\axleft, \axdown);
\coordinate (axdown_right) at ({\axleft + \length + \edge}, \axdown);

\draw [Arrow] ([xshift = {-\edge cm}, yshift = 0]axup_left.center) -- (axup_right);
\draw [Arrow] ([xshift = {-\edge cm}, yshift = 0]axdown_left.center) -- (axdown_right);


\coordinate (xi00a) at (\axleft, {\axup + 1.5 * \vscale});
\coordinate (xi00b) at ({\axleft + \tfast}, {\axup + 1.5 * \vscale});
\coordinate (xi10a) at ({\axleft + \tfast}, {\axup + 1.1 * \vscale});
\coordinate (xi10b) at ({\axleft + 2*\tfast}, {\axup + 1.1 * \vscale});
\coordinate (xi20a) at ({\axleft + 2*\tfast}, {\axup + 0.7 * \vscale});
\coordinate (xi20b) at ({\axleft + 3*\tfast}, {\axup + 0.7 * \vscale});
\coordinate (xi01a) at ({\axleft + 3*\tfast}, {\axup + 1.3 * \vscale});
\coordinate (xi01b) at ({\axleft + 4*\tfast}, {\axup + 1.3 * \vscale});
\coordinate (xi11a) at ({\axleft + 4*\tfast}, {\axup + 1 * \vscale});
\coordinate (xi11b) at ({\axleft + 5*\tfast}, {\axup + 1 * \vscale});
\coordinate (xi21a) at ({\axleft + 5*\tfast}, {\axup + 1.7 * \vscale});
\coordinate (xi21b) at ({\axleft + 6*\tfast}, {\axup + 1.7 * \vscale});
\coordinate (xi02a) at ({\axleft + 6*\tfast}, {\axup + 1.4 * \vscale});
\coordinate (xi02b) at ({\axleft + 7*\tfast}, {\axup + 1.4 * \vscale});
\coordinate (xi12a) at ({\axleft + 7*\tfast}, {\axup + 0.4 * \vscale});
\coordinate (xi12b) at ({\axleft + 8*\tfast}, {\axup + 0.4 * \vscale});
\coordinate (xi22a) at ({\axleft + 8*\tfast}, {\axup + 0.8 * \vscale});
\coordinate (xi22b) at ({\axleft + 9*\tfast}, {\axup + 0.8 * \vscale});

\coordinate (xi0a) at (\axleft, {3.3/3 * \vscale});
\coordinate (xi0b) at ({\axleft + \tslow}, {3.3/3 * \vscale});
\coordinate (xi1a) at ({\axleft + \tslow}, {4/3 * \vscale});
\coordinate (xi1b) at ({\axleft + 2*\tslow}, {4/3 * \vscale});
\coordinate (xi2a) at ({\axleft + 2*\tslow}, {2.6/3 * \vscale});
\coordinate (xi2b) at ({\axleft + 3*\tslow}, {2.6/3 * \vscale});


\coordinate (t0_upb) at (\axleft, {\axup + \lslow});
\coordinate (t1_upb) at ({\axleft+\tslow}, {\axup + \lslow});
\coordinate (t2_upb) at ({\axleft+2*\tslow}, {\axup + \lslow});
\coordinate (t3_upb) at ({\axleft+3*\tslow}, {\axup + \lslow});

\coordinate (t0_upc) at (\axleft, {\axup });
\coordinate (t1_upc) at ({\axleft+\tslow}, {\axup });
\coordinate (t2_upc) at ({\axleft+2*\tslow}, {\axup });
\coordinate (t3_upc) at ({\axleft+3*\tslow}, {\axup });

\draw [Line_semi_thick] (t0_upb) -- (t0_upc);
\draw [Line_semi_thick] (t1_upb) -- (t1_upc);
\draw [Line_semi_thick] (t2_upb) -- (t2_upc);
\draw [Line_semi_thick] (t3_upb) -- (t3_upc);

\coordinate (t0_downb) at (\axleft, {\axdown + \lslow});
\coordinate (t1_downb) at ({\axleft+\tslow}, {\axdown + \lslow});
\coordinate (t2_downb) at ({\axleft+2*\tslow}, {\axdown + \lslow});
\coordinate (t3_downb) at ({\axleft+3*\tslow}, {\axdown + \lslow});

\coordinate (t0_downc) at (\axleft, {\axdown });
\coordinate (t1_downc) at ({\axleft+\tslow}, {\axdown });
\coordinate (t2_downc) at ({\axleft+2*\tslow}, {\axdown });
\coordinate (t3_downc) at ({\axleft+3*\tslow}, {\axdown });

\draw [Line_semi_thick] (t0_downb) -- (t0_downc);
\draw [Line_semi_thick] (t1_downb) -- (t1_downc);
\draw [Line_semi_thick] (t2_downb) -- (t2_downc);
\draw [Line_semi_thick] (t3_downb) -- (t3_downc);


\coordinate (t10a) at ({\axleft + \tfast}, {\axup + \lfast});
\coordinate (t20a) at ({\axleft + 2*\tfast}, {\axup + \lfast});
\coordinate (t11a) at ({\axleft + \tfast + \tslow}, {\axup + \lfast});
\coordinate (t21a) at ({\axleft + 2*\tfast + \tslow}, {\axup + \lfast});
\coordinate (t12a) at ({\axleft + \tfast + 2*\tslow}, {\axup + \lfast});
\coordinate (t22a) at ({\axleft + 2*\tfast + 2*\tslow}, {\axup + \lfast});

\coordinate (t10b) at ({\axleft + \tfast}, {\axup});
\coordinate (t20b) at ({\axleft + 2*\tfast}, {\axup});
\coordinate (t11b) at ({\axleft + \tfast + \tslow}, {\axup});
\coordinate (t21b) at ({\axleft + 2*\tfast + \tslow}, {\axup});
\coordinate (t12b) at ({\axleft + \tfast + 2*\tslow}, {\axup});
\coordinate (t22b) at ({\axleft + 2*\tfast + 2*\tslow}, {\axup});

\draw [Line] (t10a) -- (t10b);
\draw [Line] (t20a) -- (t20b);
\draw [Line] (t11a) -- (t11b);
\draw [Line] (t21a) -- (t21b);
\draw [Line] (t12a) -- (t12b);
\draw [Line] (t22a) -- (t22b);

\node (t0) at ({\axleft + 0.5 * \tslow}, 0) [textbox] {$t-1$};
\node (t1) at ({\axleft + 1.5 * \tslow}, 0) [textbox] {$t$};
\node (t2) at ({\axleft + 2.5 * \tslow}, 0) [textbox] {$t+1$};

\node (t01) at ({\axleft + \tslow + 0.5 * \tfast}, {\axup}) [textbox] {\scriptsize{$(0, t)$}};
\node (t11) at ({\axleft + \tslow + 1.5 * \tfast}, {\axup}) [textbox] {\scriptsize{$(1, t)$}};
\node (t21) at ({\axleft + \tslow + 2.5 * \tfast}, {\axup}) [textbox] {\scriptsize{$(2, t)$}};


\coordinate (label_up_left) at ({\axleft + \tslow }, {\axup + 1.7});
\coordinate (label_up_right) at ({\axleft + \tslow + \tfast}, {\axup + 1.7});

\draw [Double_arrow] (label_up_left) -- (label_up_right);
\draw [Line_dashed] ([xshift = 0, yshift = 0cm]label_up_right) -- (xi01b);
\draw [Line_dashed] ([xshift = 0, yshift = 0cm]label_up_left) -- (xi01a);

\node (label_up) at  ({\axleft + \tslow + 0.5 * \tfast}, {\axup + 1.7}) [textbox_above] {$\frac{\Delta}{ K}$};

\coordinate (label_down_left) at ({\axleft + \tslow }, {\axdown + 1.7});
\coordinate (label_down_right) at ({\axleft + 2* \tslow }, {\axdown + 1.7});

\draw [Double_arrow] (label_down_left) -- (label_down_right);
\draw [Line_dashed] ([xshift = 0, yshift = 0cm]label_down_right) -- (xi1b);
\draw [Line_dashed] ([xshift = 0, yshift = 0cm]label_down_left) -- (xi1a);

\node (label_down) at ({\axleft + 1.5 * \tslow}, {\axdown + 1.7}) [textbox_above] {$\Delta$};


\draw [Line_r_thick] (xi00a) -- (xi00b) -- (xi10a) -- (xi10b) -- (xi20a) -- (xi20b) --
	(xi01a) -- (xi01b) -- (xi11a) -- (xi11b) -- (xi21a) -- (xi21b) --
	(xi02a) -- (xi02b) -- (xi12a) -- (xi12b) -- (xi22a) -- (xi22b);

\draw [Line_b_thick] (xi0a) -- (xi0b) -- (xi1a) -- (xi1b) -- (xi2a) -- (xi2b);

\draw [Line_r_thick_dashed] (xi00a) -- ([xshift = {-\edge cm + 0.1cm}, yshift = 0]xi00a);
\draw [Line_r_thick_dashed] (xi22b) -- ([xshift = {\edge cm - 0.1cm}, yshift = 0]xi22b);

\draw [Line_b_thick_dashed] (xi0a) -- ([xshift = {-\edge cm + 0.1cm}, yshift = 0]xi0a);
\draw [Line_b_thick_dashed] (xi2b) -- ([xshift = {\edge cm - 0.1cm}, yshift = 0]xi2b);


\coordinate (legend_up_left) at ({\axleft + 0.75 * \tslow}, {\axdown -1});
\coordinate (legend_up_right) at ([xshift = 0.5cm, yshift = 0]legend_up_left);

\draw [Line_r_thick] (legend_up_left) -- (legend_up_right);

\node at ([xshift = 0.1cm, yshift = 0cm]legend_up_right) [text centered, right] {$\xi(k,t)$};


\coordinate (legend_down_left) at ({\axleft + 1.75 * \tslow}, {\axdown -1});
\coordinate (legend_down_right) at ([xshift = 0.5cm, yshift = 0]legend_down_left);

\draw [Line_b_thick] (legend_down_left) -- (legend_down_right);

\node at ([xshift = 0.1cm, yshift = 0cm]legend_down_right) [text centered, right] {$\overline{\xi}(t)$};

\end{tikzpicture}
\caption{The above plots depicts (one component of)  the  fast time-scale disturbance process $\xi (k, t)$ and its slow time-scale average $\overline{\xi} (t)$ for $K=3$.}
\label{fig:fast}
\end{figure}
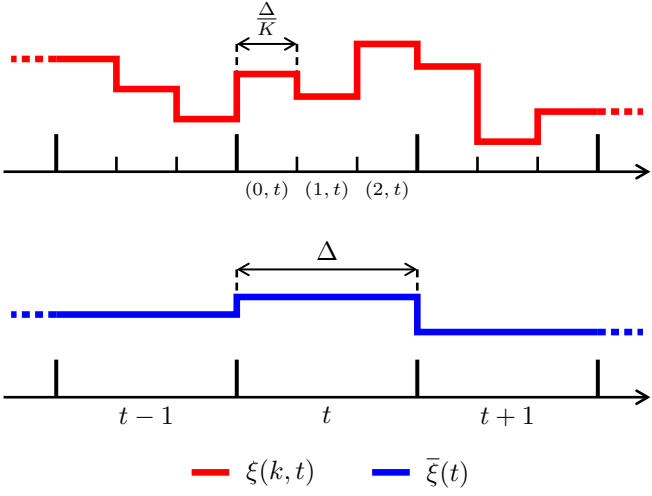

We denote the fast time-scale state, input, and disturbance processes  by $x(k, t)$, $u(k, t)$, and $\xi (k,t)$, respectively.
We emphasize that all the fast time-scale quantities have the same units as their slow time-scale counterparts. It will prove useful to define a \emph{slow time-scale average} of the fast time-scale disturbance process according to 
\begin{align} \label{eq:average}
\overline{\xi} (t) = \frac{1}{K} \left(  \sum_{k=0}^{K-1} \xi \left(k, t \right) \right),
\end{align}
for $t = 0, \dots, T-1$. 
We refer the reader to Fig. \ref{fig:fast}, which offers a graphical illustration comparing the  fast time-scale disturbance process $\xi (k, t)$ against its slow time-scale average $\overline{\xi} (t)$.

We describe the evolution of the \emph{fast time-scale state process} over each time period $t$ according to the state equation:
\begin{align*}
x(k+1, t) &= x(k, t) + \frac{1}{K} B u (k,t)
\end{align*}
for  $k = 0, \dots, K-1$. We link this process across the slow time-scale periods by enforcing the boundary conditions
$$x(K,t) = x(0, t+1)$$ for $t = 0, \dots, T-1$. We initialize the fast time-scale state process according to $x(0,0) = x(0)$.

\subsection{Fast Time-Scale Controller}
In what follows, we construct a fast time-scale controller based on the slow time-scale controller computed according to Proposition \ref{prop:primal_affine_conic}. 
More specifically, 
let $Q^*$ denote the optimal solution to problem \eqref{eqn:primal_affine_PN_conic}. 
Since $Q^* \in S$, it follows that $Q^*$ is a block lower-triangular matrix of the form
$$Q^* = \bmat{\overline{u}^* (0) & Q^* (0,0) \\
\vdots & \vdots  & \ddots \\
\overline{u}^* (T-1) & Q^* (T-1,0) & \cdots & Q^* (T-1, T-1) },$$
where each matrix $Q^* (t,s)$ is block diagonal of the form 
\begin{align*}
Q^* (t,s) = \bmat{ Q_1^* (t,s) & & \\ \vspace{-.1in} & \ddots & \\ & & Q_n^* (t,s)}
\end{align*}
for each $t = 0, \dots, T-1$ and $s = 0, \dots, t$.
Using these feedback control gains embedded in the matrix $Q^*$, we define the \emph{fast time-scale control input} at  period $(k,t)$  as
\begin{align}
u(k, t) = \overline{u}^* (t ) + Q^* (t, t) \xi \left( k , t \right) +  \sum_{s=0}^{t-1} Q^* (t,s ) \overline{\xi} (s), \label{eq:input_fast}
\end{align}
for all $t = 0, \dots, T-1$ and  $k = 0, \dots, K-1$. Recall from Eq. \eqref{eq:average} that $\overline{\xi} (t)$ denotes the average of the fast time-scale disturbance process over the period $t$.

\subsection{Constraint Satisfaction Guarantees}
The decentralized affine controller defined according to Eq. \eqref{eq:input_fast} is said to be \emph{feasible} if
it induces voltage, input, and state trajectories that  are guaranteed to satisfy  their respective constraints at the fast time-scale for all possible realizations of the fast time-scale disturbance process. That is to say, for each subsystem $i \in \{1, \dots, n\}$, it must hold that
\begin{align}
\underline{v}_i \leq  v_i  \left( k , t \right) &\leq \overline{v}_i,  \label{eq:con1}\\
-\sqrt{{s_i^{I}}^2 - {\xi_i^I(k,t)}^2}  \leq q_i^I \left( k, t \right)  &\leq \sqrt{{s_i^{I}}^2 - {\xi_i^I(k,t)}^2} ,  \label{eq:con2}\\
\underline{p}_i^S \leq p_i^S \left( k, t \right) &\leq \overline{p}_i^S,  \label{eq:con3} 
\end{align}
for all time periods $t= 0, \dots, T-1$, $k= 0, \dots, K-1$, and
\begin{align}
0 \leq x_i\left( k , t\right) &\leq b_i,    \label{eq:con4}
\end{align}
for all time periods $t= 0, \dots, T-1$, $k= 0, \dots, K$; and all possible realizations of the fast time-scale disturbance process.

We now make a \emph{mild} assumption on the support of the fast time-scale disturbance process, which ensures that the fast time-scale controller defined according to Eq. \eqref{eq:input_fast} is feasible. 
\begin{assumptio}[Quasi-Stationarity] \label{ass:fast_time}
We assume that
\begin{align*}
\left( 1, \xi \left(k_0, 0 \right), \xi \left( k_1, 1 \right), \dots, \xi \left( k_{T-1}, T-1 \right) \right) \in \Xi,
\end{align*}
for all $k_t \in \{0, \dots, K-1\}$ and $t = 0, \dots, T-1$.
\end{assumptio}
Assumption \ref{ass:fast_time} can be interpreted as an assumption on the \emph{quasi-stationarity} of the support of the fast time-scale disturbance process $\xi (k,t)$. Moreover,  Assumption \ref{ass:fast_time} is reasonable,  as it is always possible to construct a set $\Xi$ such this assumption is satisfied, given a  characterization of the set of all possible realizations taken by the fast time-scale disturbance process.

\begin{propositio}[Fast Time-Scale Feasibility]  \label{prop:fast}
Let Assumption \ref{ass:fast_time} hold. The fast time-scale controller specified according to Equation \eqref{eq:input_fast} is feasible. 
\end{propositio}
Proposition \ref{prop:fast} reveals that the \emph{slow time-scale} controller computed according to Proposition \ref{prop:primal_affine_conic} can be implemented as a feasible \emph{fast time-scale} controller.
We present its proof in  Appendix \ref{pf:propofast}.

\section{Proof of Proposition \ref{prop:fast}} \label{pf:propofast}
Let the system control input be specified according to Eq. \eqref{eq:input_fast}.
The proof consists of two parts.
In Part 1, we show that for any realization of the fast time-scale disturbance process, the input constraints specified in inequalities \eqref{eq:con1}-\eqref{eq:con3} are all satisfied.
In Part 2, we show that for any realization of the fast time-scale disturbance process, the state constraint specified in inequality \eqref{eq:con4} is satisfied.

\

\noindent \textbf{Part 1:}
We will only show that for any realization of the fast time-scale disturbance process, the voltage magnitude constraint specified in inequality \eqref{eq:con1} is satisfied. The proof of the satisfaction of the input constraints specified in inequalities \eqref{eq:con2} and \eqref{eq:con3} is analogous. It is thus omitted for the sake of brevity.

It will be convenient to work with the vector of squared voltage magnitudes $v \left( k , t \right)^2 = (v_1(k,t)^2, \dots, v_n (k,t)^2)$ for the remainder of the proof. We will show that
\begin{align*}
\underline{v}^2 \leq v \left( k , t \right)^2 \leq \overline{v}^2,
\end{align*}
where $\underline{v}^2 = (\underline{v}^2_1, \dots, \underline{v}^2_n)$, and $\overline{v}^2 = (\overline{v}_1^2, \dots, \overline{v}_n^2)$. 
It follows from the linearized branch flow model \eqref{eq:v} that for each time period $(k, t)$, the vector of squared voltage magnitudes is given by
$$v  \left( k , t \right)^2 = V_u u  \left( k , t \right) + V_{\xi} \xi  \left( k , t \right) + v_0^2 \bone ,$$
 where the matrices $V_u$ and $V_{\xi}$  are defined according to
\begin{align*}
V_u &= R \otimes \bmat{1 &0} + X \otimes \bmat{0 &1}, \\
V_{\xi} &= R \otimes \bmat{-1 &0 &1} - X \otimes \bmat{0 &1 &0}.
\end{align*}
Given the specification of the fast time-scale control input $u  \left( k , t \right)$  according to Eq. \eqref{eq:input_fast}, we have that
\begin{equation}
\begin{split}
v  \left( k , t  \right)^2 = &V_u \left( \overline{u}^* (t)   +  \sum_{s=0}^{t-1} Q^* (t,s ) \overline{\xi} (s)\right) \\
& + \left( V_u Q^* (t, t) + V_{\xi} \right) \xi \left(k , t \right) + v_0^2 \bone.
\end{split} \label{eq:v2}
\end{equation}
Given Assumption \ref{ass:fast_time} and the convexity of the set $\Xi$, it holds that
\begin{align}
\left( 1, \overline{\xi} ( 0 ) , .., \overline{\xi} ( t-1  ), \xi  ( k , t) , .., \xi ( k, T -1 )\right) \in \Xi. \label{eq:inclusion1}
\end{align}
Condition \eqref{eq:inclusion1}, in combination with the guaranteed feasibility of the control policy $Q^*$ for the original \emph{slow time-scale} decentralized control design problem \eqref{opt:decent_full}, implies that
\begin{align*}
\underline{v}^2 \leq &V_u \left( \overline{u}^* (t)   +  \sum_{s=0}^{t-1} Q^* (t,s ) \overline{\xi} (s)\right) \\
& + \left( V_u Q^* (t, t) + V_{\xi} \right) \xi \left(k , t \right) + v_0^2 \bone \leq \overline{v}^2.
\end{align*}
It immediately follows from Eq. \eqref{eq:v2} that 
$\underline{v}^2 \leq v (k, t)^2 \leq \overline{v}^2$. This completes Part 1 of the proof.

\

\noindent \textbf{Part 2:} We show that for any realization of the fast time-scale disturbance process, the system state satisfies $0 \leq x(k,t) \leq b$ for $t = 0, \dots, T-1$, $k = 0, \dots, K$. Here, $b = (b_1, \dots, b_n)$ is the vector of energy storage capacities. We fix an arbitrary realization of the fast time-scale disturbance process throughout this part of the proof.

We first consider the case of $k = 0$, and show that $0 \leq x( 0 ,t ) \leq b$ for $t = 0, \dots, T$. It holds that
\begin{align}
&x \left( 0 , t \right)  = x(0, 0) + \frac{1}{K} B \left( \sum_{s=0}^{t-1} \sum_{\ell = 0}^{K-1} u \left( \ell, s \right) \right) \label{eq:x0t_1}\\
 =  &x(0) + B \left( \sum_{s=0}^{t-1} \left( \overline{u}^* (s)  + \sum_{r=0}^s Q^* (s,r) \overline{\xi} (r) \right) \right), \label{eq:x0t_2}
\end{align}
where Eq. \eqref{eq:x0t_2} follows from the specification of the fast time-scale system control input according to Eq. \eqref{eq:input_fast}.

Condition \eqref{eq:inclusion1}, in combination with the guaranteed feasibility of the control policy $Q^*$ for the slow time-scale decentralized control design problem \eqref{opt:decent_full}, implies that
$$0 \leq x(0) + B \left( \sum_{s=0}^{t-1} \left( \overline{u}^* (s)  + \sum_{r=0}^s Q^* (s,r) \overline{\xi} (r) \right) \right) \leq b.$$
It follows from Eq. \eqref{eq:x0t_2} that  $0 \leq x(0, t) \leq b$ for $t = 0, \dots, T$.
The enforcement of the boundary condition of the fast time-scale state equation requires that
$x (K, t) = x(0, t+1)$
for $ t = 0, \dots, T-1$. This, in combination with the above inequality, implies that $0 \leq x(k, t) \leq b$ for $t = 0, \dots, T-1$, and $k = 0$ and $ K$.

Next, we show that $0 \leq x( k ,t ) \leq b$ for 
$k = 1, \dots, K-1$, $t = 0, \dots, T-1$. 
We first write $x(k,t)$ as
\begin{align*}
&x \left( k , t \right) 
=  x(0,t) + \frac{k}{K} B \left( \frac{1}{k} \sum_{\ell' = 0}^{k-1} u (\ell', t) \right) \\
=  &\frac{K - k}{K}  x(0, t) + \frac{k}{K} \left( x(0,t) +  B \left( \frac{1}{k} \sum_{\ell' = 0}^{k-1} u (\ell', t) \right) \right) , 
\end{align*}
where $x(0,t)$ is specified according to Eq. \eqref{eq:x0t_1}. 
Recall that we previously established that $0 \leq x(0, t) \leq b$. Thus, to show that $0 \leq x(k, t) \leq b$, it suffices to show that
\begin{align}
0 \leq x(0,t) + B \left( \frac{1}{k} \sum_{\ell' = 0}^{k-1} u (\ell', t) \right) \leq b. \label{eq:xkt_bound2}
\end{align}
First notice that under the fast time-scale control policy specified by Eq. \eqref{eq:input_fast}, we have that
\begin{align}
\frac{1}{k} \sum_{\ell' = 0}^{k-1} u (\ell', t) = \overline{u}^* (t) + Q^* (t,t) \widetilde{\xi} (k, t) + \sum_{s' = 0}^{t-1} Q^* (t, s') \overline{\xi} (s'), \label{eq:avg_u}
\end{align}
where the vector $\widetilde{\xi} (k, t)$ is defined according to 
\begin{align*}
\widetilde{\xi}(k, t) &= \frac{1}{k}  \sum_{\ell' = 0}^{k-1} \xi \left( \ell', t \right).
\end{align*}

Given Assumption \ref{ass:fast_time} and the convexity of $\Xi$, it holds that
\begin{align}
\left( 1, \overline{\xi} (0), .., \overline{\xi} (t-1), \widetilde{\xi} (k, t),  \xi \left(0,  t+1  \right), .., \xi \left(0,  T-1 \right) \right) \in \Xi. \label{eq:inclusion3}
\end{align}
Condition \eqref{eq:inclusion3}, in combination with the guaranteed feasibility of the control policy $Q^*$ for the slow time-scale decentralized control design problem \eqref{opt:decent_full}, implies that
\begin{align*}
0 &\leq x(0, t ) + B \left( \overline{u}^* (t) + Q^* (t,t) \widetilde{\xi} (k, t) + \sum_{s' = 0}^{t-1} Q^* (t, s') \overline{\xi} (s') \right) \\
&\leq b,
\end{align*}
where $x(0, t)$ is specified according to Eq. \eqref{eq:x0t_2}. It follows from Eq. \eqref{eq:avg_u} that inequality \eqref{eq:xkt_bound2} is satisfied. 
This completes Part 2 of the proof.

\end{appendices}

\bibliographystyle{IEEEtran}
\bibliography{decentralized_bib}{\markboth{References}{References}}

\begin{IEEEbiography}[{\includegraphics[width=1in,height=1.25in,clip,keepaspectratio]{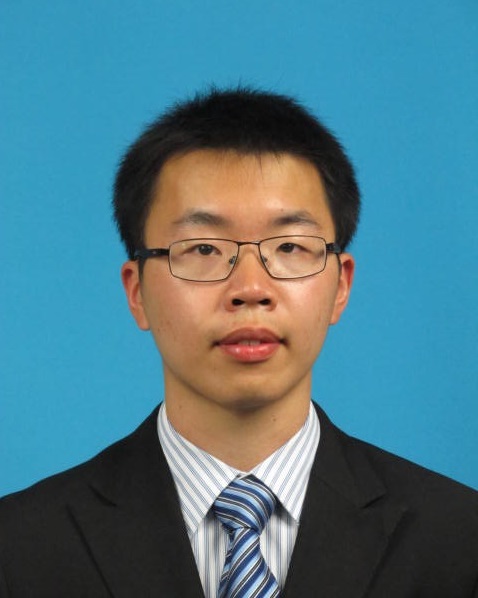}}]{Weixuan Lin} has been pursuing the Ph.D. degree in Electrical and Computer Engineering from Cornell University, Ithaca, NY, USA, since 2013. He received the B.S. degree in Electrical Engineering with high honors from Tsinghua University, Beijing, China, in 2013. 
His research interests include stochastic control, algorithmic game theory and online algorithms, with particular applications to the management of uncertainty in distributed energy resources and the design and analysis of electricity markets.
He is a recipient of the Jacobs Fellowship and the Hewlett Packard Fellowship. 
\end{IEEEbiography}

\begin{IEEEbiography}[{\includegraphics[width=1in,height=1.25in,clip,keepaspectratio]{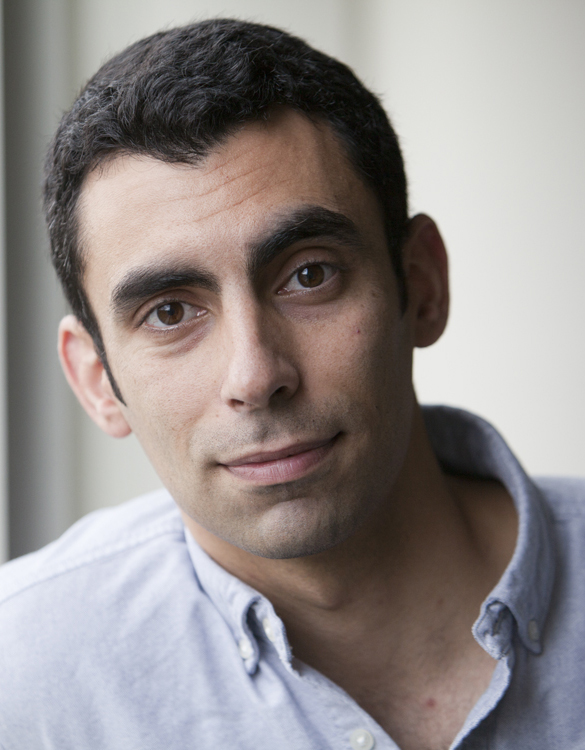}}]{Eilyan Bitar}
currently serves as an Assistant Professor and the David D. Croll Sesquicentennial Faculty Fellow in the School of Electrical and Computer Engineering at Cornell University, Ithaca, NY, USA.  Prior to joining Cornell in the Fall of 2012, he was engaged as a Postdoctoral Fellow in the department of Computing and Mathematical Sciences at the California Institute of Technology and at the University of California, Berkeley in Electrical Engineering and Computer Science, during the 2011-2012 academic year. His current research examines the operation and economics of modern power systems, with an emphasis on the design of markets and optimization methods to manage uncertainty in renewable and distributed energy resources. He received the B.S. and Ph.D. degrees in Mechanical Engineering from the University of California at Berkeley in 2006 and 2011, respectively.

Dr. Bitar is a recipient of the NSF Faculty Early Career Development Award (CAREER), the John and Janet McMurtry Fellowship, the John G. Maurer Fellowship, and the Robert F. Steidel Jr. Fellowship.
\end{IEEEbiography}

\end{document}